\documentclass[11pt,a4paper]{article}

\usepackage{slashed}
\usepackage{xkeyval}
\usepackage{url,amsmath,amssymb,latexsym,pstricks,mathrsfs,comment,amsthm,graphicx,tikz,tikz-cd,enumerate,accents,pgffor,cite,wrapfig,multicol,float,cases}
\usepackage{bibspacing}

\usepackage[colorlinks]{hyperref}

\newcommand{\nc}{\newcommand}
\nc{\rnc}{\renewcommand}

\nc{\OEIS}{}

\nc{\PT}{\mathcal{PT}}
\nc{\X}{\mathcal X}
\nc{\Y}{\mathcal Y}
\nc{\cL}{\mathrel{\mathcal L}}
\nc{\RC}{R_\leftrightarrow}
\nc{\bQ}{{\bf Q}}
\nc{\bp}{{\bf p}}
\nc{\bq}{{\bf q}}
\nc{\yb}{\overline{y}}
\newcommand{\hack}{\smash{\left\{\vphantom{\begin{aligned} \\[0.5ex] \\[0.5ex] \\[0.5ex] \end{aligned}}\right.}\!}

\nc{\pfitem}[1]{\bigskip\noindent (#1).}
\nc{\pfcase}[1]{\bigskip\noindent {\bf Case #1.}}

\usepackage{geometry} 

\let\oldproofname=\proofname
\renewcommand{\proofname}{\rm\bf{\oldproofname}}

\allowdisplaybreaks

\makeatletter

\DeclareMathSymbol{\widehatsym}{\mathord}{largesymbols}{"62}
\newcommand\lowerwidehatsym{%
  \text{\smash{\raisebox{-1.3ex}{%
    $\widehatsym$}}}}
\newcommand\fixwidehat[1]{%
  \mathchoice
    {\accentset{\displaystyle\lowerwidehatsym}{#1}}
    {\accentset{\textstyle\lowerwidehatsym}{#1}}
    {\accentset{\scriptstyle\lowerwidehatsym}{#1}}
    {\accentset{\scriptscriptstyle\lowerwidehatsym}{#1}}
}
\rnc{\widehat}{\fixwidehat}

\begin{document}

\nc{\boldsection}[1]{\section[#1]{\boldmath #1}}
\nc{\boldsubsection}[1]{\subsection[#1]{\boldmath #1}}

\nc{\ubluebox}[2]{\bluebox{#1}{1.7}{#2}2\udotted{#1}{#2}}
\nc{\lbluebox}[2]{\bluebox{#1}0{#2}{.3}\ldotted{#1}{#2}}
\nc{\ublueboxes}[1]{{
\foreach \x/\y in {#1}
{ \ubluebox{\x}{\y}}}
}
\nc{\lblueboxes}[1]{{
\foreach \x/\y in {#1}
{ \lbluebox{\x}{\y}}}
}

\nc{\bluebox}[4]{
\draw[color=blue!20, fill=blue!20] (#1,#2)--(#3,#2)--(#3,#4)--(#1,#4)--(#1,#2);
}
\nc{\redbox}[4]{
\draw[color=red!20, fill=red!20] (#1,#2)--(#3,#2)--(#3,#4)--(#1,#4)--(#1,#2);
}

\nc{\bluetrap}[8]{
\draw[color=blue!20, fill=blue!20] (#1,#2)--(#3,#4)--(#5,#6)--(#7,#8)--(#1,#2);
}
\nc{\redtrap}[8]{
\draw[color=red!20, fill=red!20] (#1,#2)--(#3,#4)--(#5,#6)--(#7,#8)--(#1,#2);
}

\usetikzlibrary{decorations.markings}
\usetikzlibrary{positioning}
\usetikzlibrary{arrows,matrix}
\usepgflibrary{arrows}
\tikzset{->-/.style={decoration={
  markings,
  mark=at position #1 with {\arrow{>}}},postaction={decorate}}}
\tikzset{-<-/.style={decoration={
  markings,
  mark=at position #1 with {\arrow{<}}},postaction={decorate}}}
\nc{\Unode}[1]{\draw(#1,-2)node{$U$};}
\nc{\Dnode}[1]{\draw(#1,-2)node{$D$};}
\nc{\Fnode}[1]{\draw(#1,-2)node{$F$};}
\nc{\Cnode}[1]{\draw(#1-.1,-2)node{$\phantom{+0},$};}
\nc{\Unodes}[1]{\foreach \x in {#1}{ \Unode{\x} }}
\nc{\Dnodes}[1]{\foreach \x in {#1}{ \Dnode{\x} }}
\nc{\Fnodes}[1]{\foreach \x in {#1}{ \Fnode{\x} }}
\nc{\Cnodes}[1]{\foreach \x in {#1}{ \Cnode{\x} }}
\nc{\Uedge}[2]{\draw[->-=0.6,line width=.3mm](#1,#2-9)--(#1+1,#2+1-9); \vertsm{#1}{#2-9} \vertsm{#1+1}{#2+1-9}}
\nc{\Dedge}[2]{\draw[->-=0.6,line width=.3mm](#1,#2-9)--(#1+1,#2-1-9); \vertsm{#1}{#2-9} \vertsm{#1+1}{#2-1-9}}
\nc{\Fedge}[2]{\draw[->-=0.6,line width=.3mm](#1,#2-9)--(#1+1,#2-9); \vertsm{#1}{#2-9} \vertsm{#1+1}{#2-9}}
\nc{\Uedges}[1]{\foreach \x/\y in {#1}{\Uedge{\x}{\y}}}
\nc{\Dedges}[1]{\foreach \x/\y in {#1}{\Dedge{\x}{\y}}}
\nc{\Fedges}[1]{\foreach \x/\y in {#1}{\Fedge{\x}{\y}}}
\nc{\xvertlabel}[1]{\draw(#1,-10+.6)node{{\tiny $#1$}};}
\nc{\yvertlabel}[1]{\draw(0-.4,-9+#1)node{{\tiny $#1$}};}
\nc{\xvertlabels}[1]{\foreach \x in {#1}{ \xvertlabel{\x} }}
\nc{\yvertlabels}[1]{\foreach \x in {#1}{ \yvertlabel{\x} }}

\nc{\bbE}{\mathbb E}
\nc{\floorn}{\lfloor\tfrac n2\rfloor}
\rnc{\sp}{\supseteq}
\rnc{\arraystretch}{1.2}

\nc{\bn}{\mathbf{n}} \nc{\bt}{\mathbf{t}} \nc{\ba}{\mathbf{a}} \nc{\bb}{\mathbf{b}} \nc{\bc}{\mathbf{c}} \nc{\bl}{\mathbf{l}} \nc{\bm}{\mathbf{m}} \nc{\bk}{\mathbf{k}} \nc{\br}{\mathbf{r}} \nc{\bs}{{\mathbf s}} \nc{\bnf}{\bnf}\nc{\bone}{\mathbf{1}}

\nc{\M}{\mathcal M}
\nc{\G}{\mathbb G}
\nc{\F}{\mathbb F}
\nc{\MnJ}{\mathcal M_n^J}
\nc{\EnJ}{\mathcal E_n^J}
\nc{\Mat}{\operatorname{Mat}}
\nc{\RegMnJ}{\Reg(\MnJ)}
\nc{\row}{\mathfrak r}
\nc{\col}{\mathfrak c}
\nc{\Row}{\operatorname{Row}}
\nc{\Col}{\operatorname{Col}}
\nc{\Span}{\operatorname{span}}
\nc{\mat}[4]{\left[\begin{matrix}#1&#2\\#3&#4\end{matrix}\right]}
\nc{\tmat}[4]{\left[\begin{smallmatrix}#1&#2\\#3&#4\end{smallmatrix}\right]}
\nc{\ttmat}[4]{{\tiny \left[\begin{smallmatrix}#1&#2\\#3&#4\end{smallmatrix}\right]}}
\nc{\tmatt}[9]{\left[\begin{smallmatrix}#1&#2&#3\\#4&#5&#6\\#7&#8&#9\end{smallmatrix}\right]}
\nc{\ttmatt}[9]{{\tiny \left[\begin{smallmatrix}#1&#2&#3\\#4&#5&#6\\#7&#8&#9\end{smallmatrix}\right]}}
\nc{\MnGn}{\M_n\sm\G_n}
\nc{\MrGr}{\M_r\sm\G_r}
\nc{\qbin}[2]{\left[\begin{matrix}#1\\#2\end{matrix}\right]_q}
\nc{\tqbin}[2]{\left[\begin{smallmatrix}#1\\#2\end{smallmatrix}\right]_q}
\nc{\qbinx}[3]{\left[\begin{matrix}#1\\#2\end{matrix}\right]_{#3}}
\nc{\tqbinx}[3]{\left[\begin{smallmatrix}#1\\#2\end{smallmatrix}\right]_{#3}}
\nc{\MNJ}{\M_nJ}
\nc{\JMN}{J\M_n}
\nc{\RegMNJ}{\Reg(\MNJ)}
\nc{\RegJMN}{\Reg(\JMN)}
\nc{\RegMMNJ}{\Reg(\MMNJ)}
\nc{\RegJMMN}{\Reg(\JMMN)}
\nc{\Wb}{\overline{W}}
\nc{\Xb}{\overline{X}}
\nc{\Yb}{\overline{Y}}
\nc{\Zb}{\overline{Z}}
\nc{\Sib}{\overline{\Si}}
\nc{\Om}{\Omega}
\nc{\Omb}{\overline{\Om}}
\nc{\Gab}{\overline{\Ga}}
\nc{\qfact}[1]{[#1]_q!}
\nc{\smat}[2]{\left[\begin{matrix}#1&#2\end{matrix}\right]}
\nc{\tsmat}[2]{\left[\begin{smallmatrix}#1&#2\end{smallmatrix}\right]}
\nc{\hmat}[2]{\left[\begin{matrix}#1\\#2\end{matrix}\right]}
\nc{\thmat}[2]{\left[\begin{smallmatrix}#1\\#2\end{smallmatrix}\right]}
\nc{\LVW}{\mathcal L(V,W)}
\nc{\KVW}{\mathcal K(V,W)}
\nc{\LV}{\mathcal L(V)}
\nc{\RegLVW}{\Reg(\LVW)}
\nc{\sM}{\mathscr M}
\nc{\sN}{\mathscr N}
\rnc{\iff}{\ \Leftrightarrow\ }
\nc{\Hom}{\operatorname{Hom}}
\nc{\End}{\operatorname{End}}
\nc{\Aut}{\operatorname{Aut}}
\nc{\Lin}{\mathcal L}
\nc{\Hommn}{\Hom(V_m,V_n)}
\nc{\Homnm}{\Hom(V_n,V_m)}
\nc{\Homnl}{\Hom(V_n,V_l)}
\nc{\Homkm}{\Hom(V_k,V_m)}
\nc{\Endm}{\End(V_m)}
\nc{\Endn}{\End(V_n)}
\nc{\Endr}{\End(V_r)}
\nc{\Autm}{\Aut(V_m)}
\nc{\Autn}{\Aut(V_n)}
\nc{\MmnJ}{\M_{mn}^J}
\nc{\MmnA}{\M_{mn}^A}
\nc{\MmnB}{\M_{mn}^B}
\nc{\Mmn}{\M_{mn}}
\nc{\Mkl}{\M_{kl}}
\nc{\Mnm}{\M_{nm}}
\nc{\EmnJ}{\mathcal E_{mn}^J}
\nc{\MmGm}{\M_m\sm\G_m}
\nc{\RegMmnJ}{\Reg(\MmnJ)}
\rnc{\implies}{\ \Rightarrow\ }
\nc{\DMmn}[1]{D_{#1}(\Mmn)}
\nc{\DMmnJ}[1]{D_{#1}(\MmnJ)}
\nc{\MMNJ}{\Mmn J}
\nc{\JMMN}{J\Mmn}
\nc{\JMMNJ}{J\Mmn J}
\nc{\Inr}{\mathcal I(V_n,W_r)}
\nc{\Lnr}{\mathcal L(V_n,W_r)}
\nc{\Knr}{\mathcal K(V_n,W_r)}
\nc{\Imr}{\mathcal I(V_m,W_r)}
\nc{\Kmr}{\mathcal K(V_m,W_r)}
\nc{\Lmr}{\mathcal L(V_m,W_r)}
\nc{\Kmmr}{\mathcal K(V_m,W_{m-r})}
\nc{\tr}{{\operatorname{T}}}
\nc{\MMN}{\MmnA(\F_1)}
\nc{\MKL}{\Mkl^B(\F_2)}
\nc{\RegMMN}{\Reg(\MmnA(\F_1))}
\nc{\RegMKL}{\Reg(\Mkl^B(\F_2))}
\nc{\gRhA}{\widehat{\mathscr R}^A}
\nc{\gRhB}{\widehat{\mathscr R}^B}
\nc{\gLhA}{\widehat{\mathscr L}^A}
\nc{\gLhB}{\widehat{\mathscr L}^B}
\nc{\timplies}{\Rightarrow}
\nc{\tiff}{\Leftrightarrow}
\nc{\Sija}{S_{ij}^a}
\nc{\dmat}[8]{\draw(#1*1.5,#2)node{$\left[\begin{smallmatrix}#3&#4&#5\\#6&#7&#8\end{smallmatrix}\right]$};}
\nc{\bdmat}[8]{\draw(#1*1.5,#2)node{${\mathbf{\left[\begin{smallmatrix}#3&#4&#5\\#6&#7&#8\end{smallmatrix}\right]}}$};}
\nc{\rdmat}[8]{\draw(#1*1.5,#2)node{\rotatebox{90}{$\left[\begin{smallmatrix}#3&#4&#5\\#6&#7&#8\end{smallmatrix}\right]$}};}
\nc{\rldmat}[8]{\draw(#1*1.5-0.375,#2)node{\rotatebox{90}{$\left[\begin{smallmatrix}#3&#4&#5\\#6&#7&#8\end{smallmatrix}\right]$}};}
\nc{\rrdmat}[8]{\draw(#1*1.5+.375,#2)node{\rotatebox{90}{$\left[\begin{smallmatrix}#3&#4&#5\\#6&#7&#8\end{smallmatrix}\right]$}};}
\nc{\rfldmat}[8]{\draw(#1*1.5-0.375+.15,#2)node{\rotatebox{90}{$\left[\begin{smallmatrix}#3&#4&#5\\#6&#7&#8\end{smallmatrix}\right]$}};}
\nc{\rfrdmat}[8]{\draw(#1*1.5+.375-.15,#2)node{\rotatebox{90}{$\left[\begin{smallmatrix}#3&#4&#5\\#6&#7&#8\end{smallmatrix}\right]$}};}
\nc{\xL}{[x]_{\! _\gL}}\nc{\yL}{[y]_{\! _\gL}}\nc{\xR}{[x]_{\! _\gR}}\nc{\yR}{[y]_{\! _\gR}}\nc{\xH}{[x]_{\! _\gH}}\nc{\yH}{[y]_{\! _\gH}}\nc{\XK}{[X]_{\! _\gK}}\nc{\xK}{[x]_{\! _\gK}}
\nc{\RegSija}{\Reg(\Sija)}
\nc{\MnmK}{\M_{nm}^K}
\nc{\cC}{\mathcal C}
\nc{\cR}{\mathcal R}
\nc{\Ckl}{\cC_k(l)}
\nc{\Rkl}{\cR_k(l)}
\nc{\Cmr}{\cC_m(r)}
\nc{\Rmr}{\cR_m(r)}
\nc{\Cnr}{\cC_n(r)}
\nc{\Rnr}{\cR_n(r)}
\nc{\Z}{\mathcal Z}

\nc{\Reg}{\operatorname{Reg}}
\nc{\RP}{\operatorname{RP}}
\nc{\TXa}{\T_X^a}
\nc{\TXA}{\T(X,A)}
\nc{\TXal}{\T(X,\al)}
\nc{\RegTXa}{\Reg(\TXa)}
\nc{\RegTXA}{\Reg(\TXA)}
\nc{\RegTXal}{\Reg(\TXal)}
\nc{\PalX}{\P_\al(X)}
\nc{\EAX}{\E_A(X)}
\nc{\Bb}{\overline{B}}
\nc{\bw}{{\bf w}}
\nc{\bz}{{\bf z}}
\nc{\TASA}{\T_A\sm\S_A}
\nc{\Ub}{\overline{U}}
\nc{\Vb}{\overline{V}}
\nc{\eb}{\overline{e}}
\nc{\EXa}{\E_X^a}
\nc{\oijr}{1\leq i<j\leq r}
\nc{\veb}{\overline{\ve}}
\nc{\bbT}{\mathbb T}
\nc{\Surj}{\operatorname{Surj}}
\nc{\Sone}{S^{(1)}}
\nc{\fillbox}[2]{\draw[fill=gray!30](#1,#2)--(#1+1,#2)--(#1+1,#2+1)--(#1,#2+1)--(#1,#2);}
\nc{\raa}{\rangle_J}
\nc{\raJ}{\rangle_J}
\nc{\Ea}{E_J}
\nc{\EJ}{E_J}
\nc{\ep}{\epsilon} \nc{\ve}{\varepsilon}
\nc{\IXa}{\I_X^a}
\nc{\RegIXa}{\Reg(\IXa)}
\nc{\JXa}{\J_X^a}
\nc{\RegJXa}{\Reg(\JXa)}
\nc{\IXA}{\I(X,A)}
\nc{\IAX}{\I(A,X)}
\nc{\RegIXA}{\Reg(\IXA)}
\nc{\RegIAX}{\Reg(\IAX)}
\nc{\trans}[2]{\left(\begin{smallmatrix} #1 \\ #2 \end{smallmatrix}\right)}
\nc{\bigtrans}[2]{\left(\begin{matrix} #1 \\ #2 \end{matrix}\right)}
\nc{\lmap}[1]{\mapstochar \xrightarrow {\ #1\ }}
\nc{\EaTXa}{E}

\nc{\gL}{\mathscr L}
\nc{\gR}{\mathscr R}
\nc{\gH}{\mathscr H}
\nc{\gJ}{\mathscr J}
\nc{\gD}{\mathscr D}
\nc{\gK}{\mathscr K}
\nc{\gLa}{\mathscr L^a}
\nc{\gRa}{\mathscr R^a}
\nc{\gHa}{\mathscr H^a}
\nc{\gJa}{\mathscr J^a}
\nc{\gDa}{\mathscr D^a}
\nc{\gKa}{\mathscr K^a}
\nc{\gLJ}{\mathscr L^J}
\nc{\gRJ}{\mathscr R^J}
\nc{\gHJ}{\mathscr H^J}
\nc{\gJJ}{\mathscr J^J}
\nc{\gDJ}{\mathscr D^J}
\nc{\gKJ}{\mathscr K^J}
\nc{\gLh}{\widehat{\mathscr L}^J}
\nc{\gRh}{\widehat{\mathscr R}^J}
\nc{\gHh}{\widehat{\mathscr H}^J}
\nc{\gJh}{\widehat{\mathscr J}^J}
\nc{\gDh}{\widehat{\mathscr D}^J}
\nc{\gKh}{\widehat{\mathscr K}^J}
\nc{\Lh}{\widehat{L}^J}
\nc{\Rh}{\widehat{R}^J}
\nc{\Hh}{\widehat{H}^J}
\nc{\Jh}{\widehat{J}^J}
\nc{\Dh}{\widehat{D}^J}
\nc{\Kh}{\widehat{K}^J}
\nc{\gLb}{\widehat{\mathscr L}}
\nc{\gRb}{\widehat{\mathscr R}}
\nc{\gHb}{\widehat{\mathscr H}}
\nc{\gJb}{\widehat{\mathscr J}}
\nc{\gDb}{\widehat{\mathscr D}}
\nc{\gKb}{\widehat{\mathscr K}}
\nc{\Lb}{\widehat{L}^J}
\nc{\Rb}{\widehat{R}^J}
\nc{\Hb}{\widehat{H}^J}
\nc{\Jb}{\widehat{J}^J}
\nc{\Db}{\overline{D}}
\nc{\Kb}{\widehat{K}}

\hyphenation{mon-oid mon-oids}

\nc{\itemit}[1]{\item[\emph{(#1)}]}
\nc{\E}{\mathcal E}
\nc{\TX}{\T(X)}
\nc{\TXP}{\T(X,\P)}
\nc{\EX}{\E(X)}
\nc{\EXP}{\E(X,\P)}
\nc{\SX}{\S(X)}
\nc{\SXP}{\S(X,\P)}
\nc{\Sing}{\operatorname{Sing}}
\nc{\idrank}{\operatorname{idrank}}
\nc{\SingXP}{\Sing(X,\P)}
\nc{\De}{\Delta}
\nc{\sgp}{\operatorname{sgp}}
\nc{\mon}{\operatorname{mon}}
\nc{\Dn}{\mathcal D_n}
\nc{\Dm}{\mathcal D_m}

\nc{\lline}[1]{\draw(3*#1,0)--(3*#1+2,0);}
\nc{\uline}[1]{\draw(3*#1,5)--(3*#1+2,5);}
\nc{\thickline}[2]{\draw(3*#1,5)--(3*#2,0); \draw(3*#1+2,5)--(3*#2+2,0) ;}
\nc{\thicklabel}[3]{\draw(3*#1+1+3*#2*0.15-3*#1*0.15,4.25)node{{\tiny $#3$}};}

\nc{\slline}[3]{\draw(3*#1+#3,0+#2)--(3*#1+2+#3,0+#2);}
\nc{\suline}[3]{\draw(3*#1+#3,5+#2)--(3*#1+2+#3,5+#2);}
\nc{\sthickline}[4]{\draw(3*#1+#4,5+#3)--(3*#2+#4,0+#3); \draw(3*#1+2+#4,5+#3)--(3*#2+2+#4,0+#3) ;}
\nc{\sthicklabel}[5]{\draw(3*#1+1+3*#2*0.15-3*#1*0.15+#5,4.25+#4)node{{\tiny $#3$}};}

\nc{\stll}[5]{\sthickline{#1}{#2}{#4}{#5} \sthicklabel{#1}{#2}{#3}{#4}{#5}}
\nc{\tll}[3]{\stll{#1}{#2}{#3}00}

\nc{\mfourpic}[9]{
\slline1{#9}0
\slline3{#9}0
\slline4{#9}0
\slline5{#9}0
\suline1{#9}0
\suline3{#9}0
\suline4{#9}0
\suline5{#9}0
\stll1{#1}{#5}{#9}{0}
\stll3{#2}{#6}{#9}{0}
\stll4{#3}{#7}{#9}{0}
\stll5{#4}{#8}{#9}{0}
\draw[dotted](6,0+#9)--(8,0+#9);
\draw[dotted](6,5+#9)--(8,5+#9);
}
\nc{\vdotted}[1]{
\draw[dotted](3*#1,10)--(3*#1,15);
\draw[dotted](3*#1+2,10)--(3*#1+2,15);
}

\nc{\Clab}[2]{
\sthicklabel{#1}{#1}{{}_{\phantom{#1}}C_{#1}}{1.25+5*#2}0
}
\nc{\sClab}[3]{
\sthicklabel{#1}{#1}{{}_{\phantom{#1}}C_{#1}}{1.25+5*#2}{#3}
}
\nc{\Clabl}[3]{
\sthicklabel{#1}{#1}{{}_{\phantom{#3}}C_{#3}}{1.25+5*#2}0
}
\nc{\sClabl}[4]{
\sthicklabel{#1}{#1}{{}_{\phantom{#4}}C_{#4}}{1.25+5*#2}{#3}
}
\nc{\Clabll}[3]{
\sthicklabel{#1}{#1}{C_{#3}}{1.25+5*#2}0
}
\nc{\sClabll}[4]{
\sthicklabel{#1}{#1}{C_{#3}}{1.25+5*#2}{#3}
}

\nc{\mtwopic}[6]{
\slline1{#6*5}{#5}
\slline2{#6*5}{#5}
\suline1{#6*5}{#5}
\suline2{#6*5}{#5}
\stll1{#1}{#3}{#6*5}{#5}
\stll2{#2}{#4}{#6*5}{#5}
}
\nc{\mtwopicl}[6]{
\slline1{#6*5}{#5}
\slline2{#6*5}{#5}
\suline1{#6*5}{#5}
\suline2{#6*5}{#5}
\stll1{#1}{#3}{#6*5}{#5}
\stll2{#2}{#4}{#6*5}{#5}
\sClabl1{#6}{#5}{i}
\sClabl2{#6}{#5}{j}
}

\nc{\keru}{\operatorname{ker}^\wedge} \nc{\kerl}{\operatorname{ker}_\vee}

\nc{\coker}{\operatorname{coker}}
\nc{\KER}{\ker}
\nc{\N}{\mathbb N}
\nc{\LaBn}{L_\al(\B_n)}
\nc{\RaBn}{R_\al(\B_n)}
\nc{\LaPBn}{L_\al(\PB_n)}
\nc{\RaPBn}{R_\al(\PB_n)}
\nc{\rhorBn}{\rho_r(\B_n)}
\nc{\DrBn}{D_r(\B_n)}
\nc{\DrPn}{D_r(\P_n)}
\nc{\DrPBn}{D_r(\PB_n)}
\nc{\DrKn}{D_r(\K_n)}
\nc{\alb}{\al_{\vee}}
\nc{\beb}{\be^{\wedge}}
\nc{\Bal}{\operatorname{Bal}}
\nc{\Red}{\operatorname{Red}}
\nc{\Pnxi}{\P_n^\xi}
\nc{\Bnxi}{\B_n^\xi}
\nc{\PBnxi}{\PB_n^\xi}
\nc{\Knxi}{\K_n^\xi}
\nc{\C}{\mathbb C}
\nc{\exi}{e^\xi}
\nc{\Exi}{E^\xi}
\nc{\eximu}{e^\xi_\mu}
\nc{\Eximu}{E^\xi_\mu}
\nc{\REF}{ {\red [Ref?]} }
\nc{\GL}{\operatorname{GL}}
\rnc{\O}{\mathcal O}

\nc{\vtx}[2]{\fill (#1,#2)circle(.2);}
\nc{\lvtx}[2]{\fill (#1,0)circle(.2);}
\nc{\uvtx}[2]{\fill (#1,1.5)circle(.2);}

\nc{\Eq}{\mathfrak{Eq}}
\nc{\Gau}{\Ga^\wedge} \nc{\Gal}{\Ga_\vee}
\nc{\Lamu}{\Lam^\wedge} \nc{\Laml}{\Lam_\vee}
\nc{\bX}{{\bf X}}
\nc{\bY}{{\bf Y}}
\nc{\ds}{\displaystyle}

\nc{\uuvert}[1]{\fill (#1,3)circle(.2);}
\nc{\uuuvert}[1]{\fill (#1,4.5)circle(.2);}
\nc{\overt}[1]{\fill (#1,0)circle(.1);}
\nc{\overtl}[3]{\node[vertex] (#3) at (#1,0) {  {\tiny $#2$} };}
\nc{\cv}[2]{\draw(#1,1.5) to [out=270,in=90] (#2,0);}
\nc{\cvs}[2]{\draw(#1,1.5) to [out=270+30,in=90+30] (#2,0);}
\nc{\ucv}[2]{\draw(#1,3) to [out=270,in=90] (#2,1.5);}
\nc{\uucv}[2]{\draw(#1,4.5) to [out=270,in=90] (#2,3);}
\nc{\textpartn}[1]{{\lower1.0 ex\hbox{\begin{tikzpicture}[xscale=.3,yscale=0.3] #1 \end{tikzpicture}}}}
\nc{\textpartnx}[2]{{\lower1.0 ex\hbox{\begin{tikzpicture}[xscale=.3,yscale=0.3] 
\foreach \x in {1,...,#1}
{ \uvert{\x} \lvert{\x} }
#2 \end{tikzpicture}}}}
\nc{\disppartnx}[2]{{\lower1.0 ex\hbox{\begin{tikzpicture}[scale=0.3] 
\foreach \x in {1,...,#1}
{ \uvert{\x} \lvert{\x} }
#2 \end{tikzpicture}}}}
\nc{\disppartnxd}[2]{{\lower2.1 ex\hbox{\begin{tikzpicture}[scale=0.3] 
\foreach \x in {1,...,#1}
{ \uuvert{\x} \uvert{\x} \lvert{\x} }
#2 \end{tikzpicture}}}}
\nc{\disppartnxdn}[2]{{\lower2.1 ex\hbox{\begin{tikzpicture}[scale=0.3] 
\foreach \x in {1,...,#1}
{ \uuvert{\x} \lvert{\x} }
#2 \end{tikzpicture}}}}
\nc{\disppartnxdd}[2]{{\lower3.6 ex\hbox{\begin{tikzpicture}[scale=0.3] 
\foreach \x in {1,...,#1}
{ \uuuvert{\x} \uuvert{\x} \uvert{\x} \lvert{\x} }
#2 \end{tikzpicture}}}}

\nc{\dispgax}[2]{{\lower0.0 ex\hbox{\begin{tikzpicture}[scale=0.3] 
#2
\foreach \x in {1,...,#1}
{\lvert{\x} }
 \end{tikzpicture}}}}
\nc{\textgax}[2]{{\lower0.4 ex\hbox{\begin{tikzpicture}[scale=0.3] 
#2
\foreach \x in {1,...,#1}
{\lvert{\x} }
 \end{tikzpicture}}}}
\nc{\textlinegraph}[2]{{\raise#1 ex\hbox{\begin{tikzpicture}[scale=0.8] 
#2
 \end{tikzpicture}}}}
\nc{\textlinegraphl}[2]{{\raise#1 ex\hbox{\begin{tikzpicture}[scale=0.8] 
\tikzstyle{vertex}=[circle,draw=black, fill=white, inner sep = 0.07cm]
#2
 \end{tikzpicture}}}}
\nc{\displinegraph}[1]{{\lower0.0 ex\hbox{\begin{tikzpicture}[scale=0.6] 
#1
 \end{tikzpicture}}}}
 
\nc{\disppartnthreeone}[1]{{\lower1.0 ex\hbox{\begin{tikzpicture}[scale=0.3] 
\foreach \x in {1,2,3,5,6}
{ \uvert{\x} }
\foreach \x in {1,2,4,5,6}
{ \lvert{\x} }
\draw[dotted] (3.5,1.5)--(4.5,1.5);
\draw[dotted] (2.5,0)--(3.5,0);
#1 \end{tikzpicture}}}}

\nc{\partn}[4]{\left( \begin{array}{c|c} 
#1 \ & \ #3 \ \ \\ \cline{2-2}
#2 \ & \ #4 \ \
\end{array} \!\!\! \right)}
\nc{\partnlong}[6]{\partn{#1}{#2}{#3,\ #4}{#5,\ #6}} 
\nc{\partnsh}[2]{\left( \begin{array}{c} 
#1 \\
#2 
\end{array} \right)}
\nc{\partncodefz}[3]{\partn{#1}{#2}{#3}{\emptyset}}
\nc{\partndefz}[3]{{\partn{#1}{#2}{\emptyset}{#3}}}
\nc{\partnlast}[2]{\left( \begin{array}{c|c}
#1 \ &  \ #2 \\
#1 \ &  \ #2
\end{array} \right)}

\nc{\uuarcx}[3]{\draw(#1,3)arc(180:270:#3) (#1+#3,3-#3)--(#2-#3,3-#3) (#2-#3,3-#3) arc(270:360:#3);}
\nc{\uuarc}[2]{\uuarcx{#1}{#2}{.4}}
\nc{\uuuarcx}[3]{\draw(#1,4.5)arc(180:270:#3) (#1+#3,4.5-#3)--(#2-#3,4.5-#3) (#2-#3,4.5-#3) arc(270:360:#3);}
\nc{\uuuarc}[2]{\uuuarcx{#1}{#2}{.4}}
\nc{\udarcx}[3]{\draw(#1,1.5)arc(180:90:#3) (#1+#3,1.5+#3)--(#2-#3,1.5+#3) (#2-#3,1.5+#3) arc(90:0:#3);}
\nc{\udarc}[2]{\udarcx{#1}{#2}{.4}}
\nc{\uudarcx}[3]{\draw(#1,3)arc(180:90:#3) (#1+#3,3+#3)--(#2-#3,3+#3) (#2-#3,3+#3) arc(90:0:#3);}
\nc{\uudarc}[2]{\uudarcx{#1}{#2}{.4}}
\nc{\darcxhalf}[3]{\draw(#1,0)arc(180:90:#3) (#1+#3,#3)--(#2,#3) ;}
\nc{\darchalf}[2]{\darcxhalf{#1}{#2}{.4}}
\nc{\uarcxhalf}[3]{\draw(#1,1.5)arc(180:270:#3) (#1+#3,1.5-#3)--(#2,1.5-#3) ;}
\nc{\uarchalf}[2]{\uarcxhalf{#1}{#2}{.4}}
\nc{\uarcxhalfr}[3]{\draw (#1+#3,1.5-#3)--(#2-#3,1.5-#3) (#2-#3,1.5-#3) arc(270:360:#3);}
\nc{\uarchalfr}[2]{\uarcxhalfr{#1}{#2}{.4}}

\nc{\bdarcx}[3]{\draw[blue](#1,0)arc(180:90:#3) (#1+#3,#3)--(#2-#3,#3) (#2-#3,#3) arc(90:0:#3);}
\nc{\bdarc}[2]{\darcx{#1}{#2}{.4}}
\nc{\rduarcx}[3]{\draw[red](#1,0)arc(180:270:#3) (#1+#3,0-#3)--(#2-#3,0-#3) (#2-#3,0-#3) arc(270:360:#3);}
\nc{\rduarc}[2]{\uarcx{#1}{#2}{.4}}
\nc{\duarcx}[3]{\draw(#1,0)arc(180:270:#3) (#1+#3,0-#3)--(#2-#3,0-#3) (#2-#3,0-#3) arc(270:360:#3);}
\nc{\duarc}[2]{\uarcx{#1}{#2}{.4}}

\nc{\uuv}[1]{\fill (#1,4)circle(.1);}
\nc{\uv}[1]{\fill (#1,2)circle(.1);}
\nc{\lv}[1]{\fill (#1,0)circle(.1);}
\nc{\uvred}[1]{\fill[red] (#1,2)circle(.1);}
\nc{\lvred}[1]{\fill[red] (#1,0)circle(.1);}
\nc{\lvwhite}[1]{\fill[white] (#1,0)circle(.1);}

\nc{\uvs}[1]{{
\foreach \x in {#1}
{ \uv{\x}}
}}
\nc{\uuvs}[1]{{
\foreach \x in {#1}
{ \uuv{\x}}
}}
\nc{\lvs}[1]{{
\foreach \x in {#1}
{ \lv{\x}}
}}

\nc{\uvreds}[1]{{
\foreach \x in {#1}
{ \uvred{\x}}
}}
\nc{\lvreds}[1]{{
\foreach \x in {#1}
{ \lvred{\x}}
}}

\nc{\uudotted}[2]{\draw [dotted] (#1,4)--(#2,4);}
\nc{\uudotteds}[1]{{
\foreach \x/\y in {#1}
{ \uudotted{\x}{\y}}
}}
\nc{\uudottedsm}[2]{\draw [dotted] (#1+.4,4)--(#2-.4,4);}
\nc{\uudottedsms}[1]{{
\foreach \x/\y in {#1}
{ \uudottedsm{\x}{\y}}
}}
\nc{\udottedsm}[2]{\draw [dotted] (#1+.4,2)--(#2-.4,2);}
\nc{\udottedsms}[1]{{
\foreach \x/\y in {#1}
{ \udottedsm{\x}{\y}}
}}
\nc{\udotted}[2]{\draw [dotted] (#1,2)--(#2,2);}
\nc{\udotteds}[1]{{
\foreach \x/\y in {#1}
{ \udotted{\x}{\y}}
}}
\nc{\ldotted}[2]{\draw [dotted] (#1,0)--(#2,0);}
\nc{\ldotteds}[1]{{
\foreach \x/\y in {#1}
{ \ldotted{\x}{\y}}
}}
\nc{\ldottedsm}[2]{\draw [dotted] (#1+.4,0)--(#2-.4,0);}
\nc{\ldottedsms}[1]{{
\foreach \x/\y in {#1}
{ \ldottedsm{\x}{\y}}
}}

\nc{\stlinest}[2]{\draw(#1,4)--(#2,0);}

\nc{\stlined}[2]{\draw[dotted](#1,2)--(#2,0);}

\nc{\tlab}[2]{\draw(#1,2)node[above]{\tiny $#2$};}
\nc{\tudots}[1]{\draw(#1,2)node{$\cdots$};}
\nc{\tldots}[1]{\draw(#1,0)node{$\cdots$};}

\nc{\uvw}[1]{\fill[white] (#1,2)circle(.1);}
\nc{\huv}[1]{\fill (#1,1)circle(.1);}
\nc{\llv}[1]{\fill (#1,-2)circle(.1);}
\nc{\arcup}[2]{
\draw(#1,2)arc(180:270:.4) (#1+.4,1.6)--(#2-.4,1.6) (#2-.4,1.6) arc(270:360:.4);
}
\nc{\harcup}[2]{
\draw(#1,1)arc(180:270:.4) (#1+.4,.6)--(#2-.4,.6) (#2-.4,.6) arc(270:360:.4);
}
\nc{\arcdn}[2]{
\draw(#1,0)arc(180:90:.4) (#1+.4,.4)--(#2-.4,.4) (#2-.4,.4) arc(90:0:.4);
}
\nc{\cve}[2]{
\draw(#1,2) to [out=270,in=90] (#2,0);
}
\nc{\hcve}[2]{
\draw(#1,1) to [out=270,in=90] (#2,0);
}
\nc{\catarc}[3]{
\draw(#1,2)arc(180:270:#3) (#1+#3,2-#3)--(#2-#3,2-#3) (#2-#3,2-#3) arc(270:360:#3);
}

\nc{\arcr}[2]{
\draw[red](#1,0)arc(180:90:.4) (#1+.4,.4)--(#2-.4,.4) (#2-.4,.4) arc(90:0:.4);
}
\nc{\arcb}[2]{
\draw[blue](#1,2-2)arc(180:270:.4) (#1+.4,1.6-2)--(#2-.4,1.6-2) (#2-.4,1.6-2) arc(270:360:.4);
}
\nc{\loopr}[1]{
\draw[blue](#1,-2) edge [out=130,in=50,loop] ();
}
\nc{\loopb}[1]{
\draw[red](#1,-2) edge [out=180+130,in=180+50,loop] ();
}
\nc{\redto}[2]{\draw[red](#1,0)--(#2,0);}
\nc{\bluto}[2]{\draw[blue](#1,0)--(#2,0);}
\nc{\dotto}[2]{\draw[dotted](#1,0)--(#2,0);}
\nc{\lloopr}[2]{\draw[red](#1,0)arc(0:360:#2);}
\nc{\lloopb}[2]{\draw[blue](#1,0)arc(0:360:#2);}
\nc{\rloopr}[2]{\draw[red](#1,0)arc(-180:180:#2);}
\nc{\rloopb}[2]{\draw[blue](#1,0)arc(-180:180:#2);}
\nc{\uloopr}[2]{\draw[red](#1,0)arc(-270:270:#2);}
\nc{\uloopb}[2]{\draw[blue](#1,0)arc(-270:270:#2);}
\nc{\dloopr}[2]{\draw[red](#1,0)arc(-90:270:#2);}
\nc{\dloopb}[2]{\draw[blue](#1,0)arc(-90:270:#2);}
\nc{\llloopr}[2]{\draw[red](#1,0-2)arc(0:360:#2);}
\nc{\llloopb}[2]{\draw[blue](#1,0-2)arc(0:360:#2);}
\nc{\lrloopr}[2]{\draw[red](#1,0-2)arc(-180:180:#2);}
\nc{\lrloopb}[2]{\draw[blue](#1,0-2)arc(-180:180:#2);}
\nc{\ldloopr}[2]{\draw[red](#1,0-2)arc(-270:270:#2);}
\nc{\ldloopb}[2]{\draw[blue](#1,0-2)arc(-270:270:#2);}
\nc{\luloopr}[2]{\draw[red](#1,0-2)arc(-90:270:#2);}
\nc{\luloopb}[2]{\draw[blue](#1,0-2)arc(-90:270:#2);}

\nc{\larcb}[2]{
\draw[blue](#1,0-2)arc(180:90:.4) (#1+.4,.4-2)--(#2-.4,.4-2) (#2-.4,.4-2) arc(90:0:.4);
}
\nc{\larcr}[2]{
\draw[red](#1,2-2-2)arc(180:270:.4) (#1+.4,1.6-2-2)--(#2-.4,1.6-2-2) (#2-.4,1.6-2-2) arc(270:360:.4);
}

\rnc{\H}{\mathrel{\mathscr H}}
\rnc{\L}{\mathrel{\mathscr L}}
\nc{\R}{\mathrel{\mathscr R}}
\nc{\D}{\mathrel{\mathscr D}}
\nc{\J}{\mathrel{\mathscr J}}

\nc{\ssim}{\mathrel{\raise0.25 ex\hbox{\oalign{$\approx$\crcr\noalign{\kern-0.84 ex}$\approx$}}}}
\nc{\POI}{\mathcal{O}}
\nc{\wb}{\overline{w}}
\nc{\ub}{\overline{u}}
\nc{\vb}{\overline{v}}
\nc{\fb}{\overline{f}}
\nc{\gb}{\overline{g}}
\nc{\hb}{\overline{h}}
\nc{\pb}{\overline{p}}
\rnc{\sb}{\overline{s}}
\nc{\XR}{\pres{X}{R\,}}
\nc{\YQ}{\pres{Y}{Q}}
\nc{\ZP}{\pres{Z}{P\,}}
\nc{\XRone}{\pres{X_1}{R_1}}
\nc{\XRtwo}{\pres{X_2}{R_2}}
\nc{\XRthree}{\pres{X_1\cup X_2}{R_1\cup R_2\cup R_3}}
\nc{\er}{\eqref}
\nc{\larr}{\mathrel{\hspace{-0.35 ex}>\hspace{-1.1ex}-}\hspace{-0.35 ex}}
\nc{\rarr}{\mathrel{\hspace{-0.35 ex}-\hspace{-0.5ex}-\hspace{-2.3ex}>\hspace{-0.35 ex}}}
\nc{\lrarr}{\mathrel{\hspace{-0.35 ex}>\hspace{-1.1ex}-\hspace{-0.5ex}-\hspace{-2.3ex}>\hspace{-0.35 ex}}}
\nc{\nn}{\tag*{}}
\nc{\tagd}[1]{\tag*{(#1)$'$}}
\nc{\ldb}{[\![}
\nc{\rdb}{]\!]}
\nc{\sm}{\setminus}
\nc{\I}{\mathcal I}
\nc{\InSn}{\I_n\setminus\S_n}
\nc{\dom}{\operatorname{dom}} \nc{\codom}{\operatorname{codom}}
\nc{\ojin}{1\leq j<i\leq n}
\nc{\eh}{\widehat{e}}
\nc{\wh}{\widehat{w}}
\nc{\uh}{\widehat{u}}
\nc{\vh}{\widehat{v}}
\nc{\sh}{\widehat{s}}
\nc{\fh}{\widehat{f}}
\nc{\textres}[1]{\text{\emph{#1}}}
\nc{\aand}{\emph{\ and \ }}
\nc{\iif}{\emph{\ if \ }}
\nc{\textlarr}{\ \larr\ }
\nc{\textrarr}{\ \rarr\ }
\nc{\textlrarr}{\ \lrarr\ }

\nc{\comma}{,\ }

\nc{\COMMA}{,\quad}
\nc{\COMMa}{,\ \ }
\nc{\TnSn}{\T_n\setminus\S_n} 
\nc{\TmSm}{\T_m\setminus\S_m} 
\nc{\TXSX}{\T_X\setminus\S_X} 
\rnc{\S}{\mathcal S}

\nc{\T}{\mathcal T} 
\nc{\A}{\mathscr A} 
\nc{\B}{\mathcal B} 
\rnc{\P}{\mathcal P} 
\nc{\K}{\mathscr K}
\nc{\PB}{\mathcal{PB}} 
\nc{\rank}{\operatorname{rank}}

\nc{\mtt}{\!\!\!\mt\!\!\!}

\nc{\sub}{\subseteq}
\nc{\la}{\langle}
\nc{\ra}{\rangle}
\nc{\mt}{\mapsto}
\nc{\im}{\mathrm{im}}
\nc{\id}{\mathrm{id}}
\nc{\al}{\alpha}
\nc{\be}{\beta}
\nc{\ga}{\gamma}
\nc{\Ga}{\Gamma}
\nc{\de}{\delta}
\nc{\ka}{\kappa}
\nc{\lam}{\lambda}
\nc{\Lam}{\Lambda}
\nc{\si}{\sigma}
\nc{\Si}{\Sigma}
\nc{\oijn}{1\leq i<j\leq n}
\nc{\oijm}{1\leq i<j\leq m}

\nc{\comm}{\rightleftharpoons}
\nc{\AND}{\qquad\text{and}\qquad}
\nc{\ANd}{\quad\text{and}\quad}

\nc{\bit}{\begin{itemize}}
\nc{\bitbmc}{\begin{itemize}\begin{multicols}}
\nc{\bmc}{\begin{itemize}\begin{multicols}}
\nc{\emc}{\end{multicols}\end{itemize}}
\nc{\eit}{\end{itemize}}
\nc{\ben}{\begin{enumerate}}
\nc{\een}{\end{enumerate}}
\nc{\eitres}{\end{itemize}}

\nc{\set}[2]{\{ {#1} : {#2} \}} 
\nc{\bigset}[2]{\big\{ {#1}: {#2} \big\}} 
\nc{\Bigset}[2]{\left\{ \,{#1} :{#2}\, \right\}}

\nc{\pres}[2]{\la {#1} : {#2} \ra} \nc{\bigpres}[2]{\big\la {#1} : {#2} \big\ra} \nc{\Bigpres}[2]{\Big\la \,{#1}\, : \,{#2}\, \Big\ra} \nc{\Biggpres}[2]{\Bigg\la {#1} : {#2} \Bigg\ra}

\nc{\pf}{\begin{proof}}
\nc{\epf}{\end{proof}}
\nc{\epfres}{\hfill\qed}
\nc{\pfnb}{\pf}
\nc{\epfnb}{\bigskip}
\nc{\pfthm}[1]{\bigskip \noindent{\bf Proof of Theorem \ref{#1}.}  } 
\nc{\pfprop}[1]{\bigskip \noindent{\bf Proof of Proposition \ref{#1}}  } 

\nc{\uvert}[1]{\fill (#1,2)circle(.2);}
\rnc{\lvert}[1]{\fill (#1,0)circle(.2);}
\nc{\guvert}[1]{\fill[lightgray] (#1,2)circle(.2);}
\nc{\glvert}[1]{\fill[lightgray] (#1,0)circle(.2);}
\nc{\uvertx}[2]{\fill (#1,#2)circle(.2);}
\nc{\guvertx}[2]{\fill[lightgray] (#1,#2)circle(.2);}
\nc{\uvertxs}[2]{
\foreach \x in {#1}
{ \uvertx{\x}{#2}  }
}
\nc{\guvertxs}[2]{
\foreach \x in {#1}
{ \guvertx{\x}{#2}  }
}

\nc{\uvertth}[2]{\fill (#1,2)circle(#2);}
\nc{\lvertth}[2]{\fill (#1,0)circle(#2);}
\nc{\uvertths}[2]{
\foreach \x in {#1}
{ \uvertth{\x}{#2}  }
}
\nc{\lvertths}[2]{
\foreach \x in {#1}
{ \lvertth{\x}{#2}  }
}

\nc{\vertlabel}[2]{\draw(#1,2+.3)node{{\tiny $#2$}};}
\nc{\vertlabelh}[2]{\draw(#1,2+.4)node{{\tiny $#2$}};}
\nc{\vertlabelhh}[2]{\draw(#1,2+.6)node{{\tiny $#2$}};}
\nc{\vertlabelhhh}[2]{\draw(#1,2+.64)node{{\tiny $#2$}};}
\nc{\vertlabelup}[2]{\draw(#1,4+.6)node{{\tiny $#2$}};}
\nc{\vertlabels}[1]{
{\foreach \x/\y in {#1}
{ \vertlabel{\x}{\y} }
}
}

\nc{\dvertlabel}[2]{\draw(#1,-.4)node{{\tiny $#2$}};}
\nc{\dvertlabels}[1]{
{\foreach \x/\y in {#1}
{ \dvertlabel{\x}{\y} }
}
}
\nc{\vertlabelsh}[1]{
{\foreach \x/\y in {#1}
{ \vertlabelh{\x}{\y} }
}
}
\nc{\vertlabelshh}[1]{
{\foreach \x/\y in {#1}
{ \vertlabelhh{\x}{\y} }
}
}
\nc{\vertlabelshhh}[1]{
{\foreach \x/\y in {#1}
{ \vertlabelhhh{\x}{\y} }
}
}

\nc{\vertlabelx}[3]{\draw(#1,2+#3+.6)node{{\tiny $#2$}};}
\nc{\vertlabelxs}[2]{
{\foreach \x/\y in {#1}
{ \vertlabelx{\x}{\y}{#2} }
}
}

\nc{\vertlabelupdash}[2]{\draw(#1,2.7)node{{\tiny $\phantom{'}#2'$}};}
\nc{\vertlabelupdashess}[1]{
{\foreach \x/\y in {#1}
{\vertlabelupdash{\x}{\y}}
}
}

\nc{\vertlabeldn}[2]{\draw(#1,0-.6)node{{\tiny $\phantom{'}#2'$}};}
\nc{\vertlabeldnph}[2]{\draw(#1,0-.6)node{{\tiny $\phantom{'#2'}$}};}

\nc{\vertlabelups}[1]{
{\foreach \x in {#1}
{\vertlabel{\x}{\x}}
}
}
\nc{\vertlabeldns}[1]{
{\foreach \x in {#1}
{\vertlabeldn{\x}{\x}}
}
}
\nc{\vertlabeldnsph}[1]{
{\foreach \x in {#1}
{\vertlabeldnph{\x}{\x}}
}
}

\nc{\dotsup}[2]{\draw [dotted] (#1+.6,2)--(#2-.6,2);}
\nc{\dotsupx}[3]{\draw [dotted] (#1+.6,#3)--(#2-.6,#3);}
\nc{\dotsdn}[2]{\draw [dotted] (#1+.6,0)--(#2-.6,0);}
\nc{\dotsups}[1]{\foreach \x/\y in {#1}
{ \dotsup{\x}{\y} }
}
\nc{\dotsupxs}[2]{\foreach \x/\y in {#1}
{ \dotsupx{\x}{\y}{#2} }
}
\nc{\dotsdns}[1]{\foreach \x/\y in {#1}
{ \dotsdn{\x}{\y} }
}

\nc{\nodropcustpartn}[3]{
\begin{tikzpicture}[scale=.3]
\foreach \x in {#1}
{ \uvert{\x}  }
\foreach \x in {#2}
{ \lvert{\x}  }
#3 \end{tikzpicture}
}

\nc{\custpartn}[3]{{\lower1.4 ex\hbox{
\begin{tikzpicture}[scale=.3]
\foreach \x in {#1}
{ \uvert{\x}  }
\foreach \x in {#2}
{ \lvert{\x}  }
#3 \end{tikzpicture}
}}}

\nc{\smcustpartn}[3]{{\lower0.7 ex\hbox{
\begin{tikzpicture}[scale=.2]
\foreach \x in {#1}
{ \uvert{\x}  }
\foreach \x in {#2}
{ \lvert{\x}  }
#3 \end{tikzpicture}
}}}

\nc{\dropcustpartn}[3]{{\lower5.2 ex\hbox{
\begin{tikzpicture}[scale=.3]
\foreach \x in {#1}
{ \uvert{\x}  }
\foreach \x in {#2}
{ \lvert{\x}  }
#3 \end{tikzpicture}
}}}

\nc{\dropcustpartnx}[4]{{\lower#4 ex\hbox{
\begin{tikzpicture}[scale=.4]
\foreach \x in {#1}
{ \uvert{\x}  }
\foreach \x in {#2}
{ \lvert{\x}  }
#3 \end{tikzpicture}
}}}

\nc{\dropcustpartnxy}[3]{\dropcustpartnx{#1}{#2}{#3}{4.6}}

\nc{\uvertsm}[1]{\fill (#1,2)circle(.15);}
\nc{\lvertsm}[1]{\fill (#1,0)circle(.15);}
\nc{\vertsm}[2]{\fill (#1,#2)circle(.15);}

\nc{\bigdropcustpartn}[3]{{\lower6.93 ex\hbox{
\begin{tikzpicture}[scale=.6]
\foreach \x in {#1}
{ \uvertsm{\x}  }
\foreach \x in {#2}
{ \lvertsm{\x}  }
#3 \end{tikzpicture}
}}}

\nc{\gcustpartn}[5]{{\lower1.4 ex\hbox{
\begin{tikzpicture}[scale=.3]
\foreach \x in {#1}
{ \uvert{\x}  }
\foreach \x in {#2}
{ \guvert{\x}  }
\foreach \x in {#3}
{ \lvert{\x}  }
\foreach \x in {#4}
{ \glvert{\x}  }
#5 \end{tikzpicture}
}}}

\nc{\gcustpartndash}[5]{{\lower6.97 ex\hbox{
\begin{tikzpicture}[scale=.3]
\foreach \x in {#1}
{ \uvert{\x}  }
\foreach \x in {#2}
{ \guvert{\x}  }
\foreach \x in {#3}
{ \lvert{\x}  }
\foreach \x in {#4}
{ \glvert{\x}  }
#5 \end{tikzpicture}
}}}

\nc{\stline}[2]{\draw(#1,2)--(#2,0);}
\nc{\stlines}[1]{
{\foreach \x/\y in {#1}
{ \stline{\x}{\y} }
}
}

\nc{\stlinests}[1]{
{\foreach \x/\y in {#1}
{ \stlinest{\x}{\y} }
}
}

\nc{\stlineds}[1]{
{\foreach \x/\y in {#1}
{ \stlined{\x}{\y} }
}
}

\nc{\gstline}[2]{\draw[lightgray](#1,2)--(#2,0);}
\nc{\gstlines}[1]{
{\foreach \x/\y in {#1}
{ \gstline{\x}{\y} }
}
}

\nc{\gstlinex}[3]{\draw[lightgray](#1,2+#3)--(#2,0+#3);}
\nc{\gstlinexs}[2]{
{\foreach \x/\y in {#1}
{ \gstlinex{\x}{\y}{#2} }
}
}

\nc{\stlinex}[3]{\draw(#1,2+#3)--(#2,0+#3);}
\nc{\stlinexs}[2]{
{\foreach \x/\y in {#1}
{ \stlinex{\x}{\y}{#2} }
}
}

\nc{\darcx}[3]{\draw(#1,0)arc(180:90:#3) (#1+#3,#3)--(#2-#3,#3) (#2-#3,#3) arc(90:0:#3);}
\nc{\darc}[2]{\darcx{#1}{#2}{.4}}
\nc{\uarcx}[3]{\draw(#1,2)arc(180:270:#3) (#1+#3,2-#3)--(#2-#3,2-#3) (#2-#3,2-#3) arc(270:360:#3);}
\nc{\uarc}[2]{\uarcx{#1}{#2}{.4}}

\nc{\darcxx}[4]{\draw(#1,0+#4)arc(180:90:#3) (#1+#3,#3+#4)--(#2-#3,#3+#4) (#2-#3,#3+#4) arc(90:0:#3);}
\nc{\uarcxx}[4]{\draw(#1,2+#4)arc(180:270:#3) (#1+#3,2-#3+#4)--(#2-#3,2-#3+#4) (#2-#3,2-#3+#4) arc(270:360:#3);}

\makeatletter
\newcommand\footnoteref[1]{\protected@xdef\@thefnmark{\ref{#1}}\@footnotemark}
\makeatother

\numberwithin{equation}{section}

\newtheorem{thm}[equation]{Theorem}
\newtheorem{lemma}[equation]{Lemma}
\newtheorem{cor}[equation]{Corollary}
\newtheorem{prop}[equation]{Proposition}
\newtheorem{conj}[equation]{Conjecture}

\theoremstyle{definition}

\newtheorem{rem}[equation]{Remark}
\newtheorem{defn}[equation]{Definition}
\newtheorem{eg}[equation]{Example}
\newtheorem{ass}[equation]{Assumption}
\newtheorem{prob}{Problem}

%
%
%
%

\title{Presentations for singular wreath products~\vspace{-5ex}}
\author{}
\date{}

\maketitle
\begin{center}
{\large 
Ying-Ying Feng,%
\hspace{-.2em}\footnote{School of Mathematics and Big Data, Foshan University, Guangdong, P. R. China 528000. {\it Email:} {\tt rickyfungyy@fosu.edu.cn}}
Asawer Al-Aadhami,%
\hspace{-.2em}\footnote{Department of Mathematics, College of Science, University of Baghdad, Baghdad, Iraq.  Department of Mathematics, University of York, Heslington, York YO10 5DD, UK. {\it Email:} {\tt adaa501@york.ac.uk}, {\tt asawerduraid@gmail.com}}
Igor Dolinka,%
\hspace{-.2em}\footnote{Department of Mathematics and Informatics, University of Novi Sad, Trg Dositeja Obradovi\'ca 4, 21101 Novi Sad, Serbia. {\it Email:} {\tt dockie@dmi.uns.ac.rs}}
James East,%
\hspace{-.2em}\footnote{Centre for Research in Mathematics, School of Computing, Engineering and Mathematics, Western Sydney University, Locked Bag 1797, Penrith NSW 2751, Australia. {\it Email:} {\tt j.east@westernsydney.edu.au}}
Victoria Gould
\hspace{-.2em}\footnote{Department of Mathematics, University of York, Heslington, York YO10 5DD, UK. {\it Email:} {\tt victoria.gould@york.ac.uk}}
}
\end{center}

\vspace{-0.5cm}

\begin{abstract}
For a monoid $M$ and a subsemigroup $S$ of the full transformation semigroup $\T_n$, the wreath product $M\wr S$ is defined to be the semidirect product $M^n\rtimes S$, with the coordinatewise action of $S$ on $M^n$.  
The full wreath product $M\wr\T_n$ is isomorphic to the endomorphism monoid of the free $M$-act on $n$ generators.
Here, we are particularly interested in the case that $S=\Sing_n$ is the singular part of $\T_n$, consisting of all non-invertible transformations.  
Our main results are presentations for $M\wr\Sing_n$ in terms of certain natural generating sets, and we prove these via general results on semidirect products and wreath products.
We re-prove a classical result of Bulman-Fleming that $M\wr\Sing_n$ is idempotent-generated if and only if the set $M/{\L}$ of $\L$-classes of $M$ forms a chain under the usual ordering of $\L$-classes, and we give a presentation for $M\wr\Sing_n$ in terms of idempotent generators for such a monoid $M$.
Among other results, we also give estimates for the minimal size of a generating set for $M\wr\Sing_n$, as well as exact values in some cases (including the case that $M$ is finite and $M/{\L}$ is a chain, in which case we also calculate the minimal size of an idempotent generating set).
As an application of our results, we obtain a presentation (with idempotent generators) for the idempotent-generated subsemigroup of the endomorphism monoid of a uniform partition of a finite set.

{\it Keywords}: Wreath products, semidirect products, transformation semigroups, presentations, rank, idempotent rank.

MSC: 20M05; 20M20.
\end{abstract}

\tableofcontents

\section{Introduction}\label{sect:intro}

We draw our main inspiration from three closely related themes in semigroup theory: idempotent-generation, endomorphism monoids, and presentations.
Idempotents have long played an important role in algebraic and combinatorial semigroup theory.  On the one hand, Howie's 1966 article \cite{Howie1966} on singular transformation semigroups initiated a vast research theme in idempotent-generated semigroups, with much attention focused on (partial) endomorphism monoids of sets, vector spaces, modules, independence algebras, free acts and other structures \cite{Ballantine1978,Ruitenburg1993, BR2009, EastGray,GH1987,FG2007,Gray2008,Gray2007,Howie1990,Howie1978,
LX2008,SZ2014,AJLL2014, JEpnsn,EF2012,BFF1997,BF1995,RS1985,Erdos1967,Dawlings81/82,FL1992, FL1993,Laffey1983,Fountain1991}.  
Among other things, Howie's article \cite{Howie1966} demonstrated a universal property: every semigroup $S$ embeds in an idempotent-generated (singular transformation) semigroup that may be taken to be finite if $S$ is finite.  
On the other hand, every idempotent-generated semigroup $T$ is a homomorphic image of a so-called \emph{free idempotent-generated semigroup} that has the same \emph{biordered set} of idempotents as $T$.  These semigroups (henceforth referred to as FIGSs) are defined by means of a presentation (by generators and relations) derived from abstract biordered sets \cite{Nambooripad1979,Easdown1985}.  Among other things, these FIGSs encode combinatorial and topological properties of semigroups, as their maximal subgroups are isomorphic to the fundamental groups of certain complexes arising from the biordered set \cite{BMM2009}.  These maximal subgroups have therefore been the subject of intense study, particularly in the last decade or so.  Early results \cite{NP1980,Pastijn1980,McElwee2002} led to a conjecture that these maximal subgroups were always free, but this turned out to be false.  The first counter-example was provided in \cite{BMM2009} and, soon after, it was shown by Gray and Ru\v skuc \cite{GR2012_2} that \emph{every} group occurs as a maximal subgroup of some FIGS.  A recent focus has, therefore, been to describe the maximal subgroups of FIGSs arising from biordered sets of well-known semigroups \cite{Dandan2016,GR2012,Dolinka2013,DG2014,DGY2015,BMM2010}, while relatively fewer studies of the global structure of a FIGS have been carried out \cite{Dandan2016,DG2016_2,DGR2014}.
The Gray-Ru\v skuc result \cite{GR2012_2} has now been proved in a number of different ways \cite{DR2013,GY2014,DGY2015}, with endomorphism monoids of free $G$-acts playing a key role in some of these later proofs, and providing a natural biordered set for the construction.  It has long been known that the endomorphism monoid of a free $G$-act of (finite) rank $n$ is isomorphic to a wreath product~$G\wr\T_n$; in fact, this is true in the more general case of free $M$-acts, where $M$ is a monoid \cite{Skornjakov1979,KKM2000}.  It is therefore very natural to study the structure of wreath products $M\wr\T_n$ or, more generally, $M\wr S$ for an arbitrary subsemigroup $S$ of $\T_n$.  Further motivation for studying such wreath products comes from the fact that the wreath product $\T_m\wr\T_n$ of two transformation semigroups is isomorphic to the endomorphism monoid of a uniform partition with $n$ blocks of size $m$ \cite{Pei1994,DE2016,AS2009}.

Close to the time that Howie's article \cite{Howie1966} appeared, a number of authors obtained presentations for various monoids of (partial) transformations \cite{Aizenstat1958,Aizenstat1962,Popova1962,Popova1961}, sparking a field of research that continues to this day; see for example \cite{EL2004,KM2006,JEgrpm, Fernandes2001,Lavers1997}, and especially the survey  \cite{Fernandes2001_2} and references therein.  Presentations have recently been obtained for certain singular semigroups of transformations and related structures \cite{JEtnsn2,JEptnsn2,JEpnsn,MM2007,JEinsn2}.  In particular, a presentation for $\Sing_n$, the singular part of $\T_n$, is given in \cite{JEtnsn2} in terms of the generating set consisting of all idempotents of rank $n-1$.  To the authors' knowledge, a general presentation for the endomorphism monoid $\End(A)$ of an arbitrary independence algebra $A$ is not currently known.  But for a special subclass of such algebras, the above-mentioned free $G$-acts of finite rank, such a presentation can be described using results of Lavers \cite{Lavers1998} on general products of \emph{monoids}, since (as noted above) these endomorphism monoids are isomorphic to wreath products of the form $G\wr\T_n$.  Here, we are interested in the more general problem of finding presentations for wreath products $M\wr S$ for an arbitrary monoid $M$ and an arbitrary \emph{subsemigroup} $S\sub\T_n$, particularly in the case that $S=\Sing_n$.  This kind of problem can be quite difficult in the case that $S$ does not contain the identity transformation (as happens when $S=\Sing_n$, for example), since many articles on presentations for semigroup constructions (including wreath and semidirect products) focus on the case of monoids \cite{RRT2003,Lavers1998,FQ2011,HR1994}; notable exceptions that are not restricted to monoids have concentrated on constructions that do not capture the kind of wreath and semidirect products that arise from endomorphisms of $M$-acts \cite{DR2009,RRT2002}.  As such, to achieve our main aim of finding presentations for wreath products $M\wr\Sing_n$, we first prove general results on presentations for arbitrary semidirect products $M\rtimes S$ (of which the wreath product is a special case) where $M$ is a monoid and $S$ a semigroup.  From these, we are able to deduce a number of presentations for $M\wr\Sing_n$ that extend the presentation for $\Sing_n$ from \cite{JEtnsn2}.  Along the way, we obtain several other results of independent interest, as we describe below.

The article is organised as follows.  In Section \ref{sect:prelim}, we establish notation and gather some background results on (transformation) semigroups and presentations.   In Section \ref{sect:semidirect}, we give a general presentation for a semidirect product $M\rtimes S$ of a monoid $M$ and semigroup $S$ (Theorem \ref{thm:semidirect}), including wreath products as a special case.  In Section \ref{sect:gen}, we classify and enumerate the idempotents of $M\wr\T_n$ (Proposition \ref{prop:idempotents} and Corollary~\ref{cor:idempotents}), describe canonical generating sets for $M\wr\Sing_n$ (Theorem \ref{thm:IG}), give necessary and sufficient conditions for $M\wr\Sing_n$ to be idempotent-generated (also contained in Theorem \ref{thm:IG}, and originally proved by Bulman-Fleming \cite{BF1995}), and give bounds (and exact values in the idempotent-generated case) for the minimal size of a generating set for finite $M\wr\Sing_n$ (Proposition \ref{prop:rank1} and Theorem~\ref{thm:rank2}; see also Theorem \ref{thm:infinite}).  In Section~\ref{sect:pres}, we give a number of presentations for $M\wr\Sing_n$ with respect to the canonical generating sets alluded to above (Corollary \ref{cor:Xn} and Theorem \ref{thm:X2}); in the case that $M\wr\Sing_n$ is idempotent-generated, we give a presentation in terms of a particularly natural idempotent generating set (Theorems~\ref{thm:X1} and \ref{thm:X1group}).  Finally, in Section~\ref{sect:EndP}, we apply our results to obtain presentations for the idempotent-generated subsemigroups of a class of monoids including the endomorphism monoid of a uniform partition of a finite set (see Theorem~\ref{thm:IGMwrTn}).

\section{Preliminaries}\label{sect:prelim}

Let $S$ be a semigroup, and write $S^1$ for the monoid obtained by adjoining an identity $1$ to $S$ if necessary.  Unless otherwise specified, we will generally write $1$ for the identity element of any monoid.  Recall that \emph{Green's preorders} are defined, for $a,b\in S$, by
\[
a\leq_{\R} b \iff a\in bS^1 \COMMA a\leq_{\L} b \iff a\in S^1b \COMMA a\leq_{\J} b \iff a\in S^1bS^1,
\]
and that \emph{Green's relations} are defined by
\[
{\R}={\leq_{\R}}\cap{\geq_{\R}} \COMMA {\L}={\leq_{\L}}\cap{\geq_{\L}} \COMMA {\J}={\leq_{\J}}\cap{\geq_{\J}} \COMMA {\H}={\R}\cap{\L} \COMMA {\D}={\R}\circ{\L}.
\]
It can easily be proved that ${\D}={\R}\circ{\L}={\L}\circ{\R}$ (so that $\D$ is the least upper bound of $\R$ and $\L$ in the lattice of equivalence relations on $S$).
If $S$ is finite, then~${\D}={\J}$.  If $a\in S$, and $\K$ is one of Green's relations, we write $K_a$ for the $\K$-class of $S$ containing $a$.  If a subset $A\sub S$ is a union of $\K$-classes, we write $A/\K=\set{K_a}{a\in A}$ for the set of all $\K$-classes contained in $A$.  If $\K$ is one of $\R$, $\L$ or $\J$, then $\leq_\K$ induces a partial order on $S/\K$; as usual, we will generally write $\leq_\K$ for this induced order, so that $a\leq_\K b$ in $S$ $\iff$ $K_a\leq_\K K_b$ in $S/\K$.  For any subset $A\sub S$, we write $E(A)=\set{a\in A}{a=a^2}$ for the set of idempotents in $A$, and $\la A\ra$ for the subsemigroup of $S$ generated by $A$.  The \emph{rank} of a semigroup $S$, denoted $\rank(S)$, is the smallest size of a generating set for $S$.  If $S$ is idempotent-generated, then the \emph{idempotent rank}, denoted $\idrank(S)$, is defined analogously with respect to generating sets consisting of idempotents.  For more background on semigroups, see \cite{Hig,Howie}.


Let $X$ be an \emph{alphabet} (a set whose elements are called \emph{letters}), and denote by $X^+$ (resp., $X^*$) the free semigroup (resp., free monoid) on~$X$.  We denote the \emph{empty word} (over any alphabet) by $1$, so $X^*=X^+\cup\{1\}$.  If~$R\sub X^+\times X^+$ (resp., $R\sub X^*\times X^*$), we denote by $R^\sharp$ the congruence on $X^+$ (resp., $X^*$) generated by~$R$.  We say a semigroup (resp., monoid) $S$ has \emph{semigroup} (resp., \emph{monoid}) \emph{presentation} $\pres{X\!}{\!R}$ if~${S\cong X^+/R^\sharp}$ (resp., $S\cong X^*/R^\sharp$) or, equivalently, if there is an epimorphism ${X^{+}}\to S$ (resp., $X^*\to S$) with kernel $R^\sharp$.  If $\phi$ is such an epimorphism, we say $S$ has \emph{presentation $\pres{X\!}{\!R}$ via $\phi$}.  The elements of $R$ are generally referred to as \emph{relations}, and a relation $(w_1,w_2)\in R$ will usually be displayed as an equation: $w_1=w_2$.  
Unless otherwise specified (and this will only occur in Section \ref{sect:EndP}), all presentations will be semigroup presentations.
%
%

For an integer $n\geq0$, we write $\bn=\{1,\ldots,n\}$ and $\T_n$ for the \emph{full transformation semigroup of degree $n$}, which consists of all transformations of $\bn$ (i.e., all maps $\bn\to\bn$) under composition.  (When $n=0$, $\bn=\emptyset$ and $\T_0$ consists only of the empty function $\emptyset\to\emptyset$.)  For $\al\in\T_n$ and $i\in\bn$, we write $i\al$ for the image of $i$ under $\al$, so that transformations are composed left-to-right.  For $\al\in\T_n$, define
\[
\im(\al)=\set{i\al}{i\in\bn} \COMMA \ker(\al)=\set{(i,j)\in\bn\times\bn}{i\al=j\al} \COMMA \rank(\al)=|\im(\al)|.
\]
It is well known (see \cite[Exercise 16]{Howie}) that, for $\al,\be\in\T_n$,
\[
\al\leq_{\R}\be \iff \ker(\al)\supseteq\ker(\be) \COMMA \al\leq_{\L}\be \iff \im(\al)\sub\im(\be) \COMMA \al\leq_{\J}\be \iff \rank(\al)\leq\rank(\be).
\]
Since $\T_n$ is finite (of size $n^n$), ${\D}={\J}$.  The next result is well known, and is easily proved.

\begin{prop}\label{prop:ETn}
A transformation $\al\in\T_n$ is an idempotent if and only if the restriction $\al|_{\im(\al)}$ of $\al$ to $\im(\al)$ is the identity map.  Consequently, $|E(\T_n)| = \sum_{k=1}^n\binom nk k^{n-k}$. \epfres
\end{prop}

The group of units of $\T_n$ is the \emph{symmetric group} $\S_n=\set{\al\in\T_n}{\rank(\al)=n}$.  We write $\Sing_n=\TnSn$ for the \emph{singular part} of $\T_n$, which consists of all non-invertible (i.e., \emph{singular}) transformations of $\bn$.  A famous result of Howie \cite{Howie1966} states that $\Sing_n$ is generated by its idempotents: in fact, by its idempotents of rank~$n-1$.
The latter are precisely the maps $\ve_{ij}$ (for $i,j\in\bn$ with $i\not=j$) defined by
\[
k\ve_{ij} = \begin{cases}
k &\text{if $k\not=j$}\\
i &\text{if $k=j$.}
\end{cases}
\]
As usual, these idempotents may be represented diagrammatically, for $\oijn$, by
\[
\ve_{ij} = 
\custpartn{1,3,4,5,7,8,9,11}{1,3,4,5,7,8,9,11}{\dotsups{1/3,5/7,9/11}\dotsdns{1/3,5/7,9/11}\stlines{1/1,3/3,4/4,5/5,7/7,8/4,9/9,11/11}\vertlabelshh{1/1,4/i,8/j,11/n}}
\AND
\ve_{ji} = 
\custpartn{1,3,4,5,7,8,9,11}{1,3,4,5,7,8,9,11}{\dotsups{1/3,5/7,9/11}\dotsdns{1/3,5/7,9/11}\stlines{1/1,3/3,4/8,5/5,7/7,8/8,9/9,11/11}\vertlabelshh{1/1,4/i,8/j,11/n}}
.
\]
We will write $\X=\set{\ve_{ij}}{i,j\in\bn,\ i\not=j}$ for the set of all rank $n-1$ idempotents from $\T_n$.  Note that $|\X|=2\binom n2=n(n-1)$.  It is easy to check (and follows from facts mentioned above) that for all $i,j,k,l\in\bn$ with $i\not=j$ and $k\not=l$,
\[
\ve_{ij}\R\ve_{kl} \iff \{i,j\}=\{k,l\} \AND \ve_{ij}\L\ve_{kl} \iff j=l.
\]
The next result is \cite[Theorem I]{Howie1966}.

\begin{thm}\label{thm:howie}
If $n\geq2$, then $\Sing_n=\la\X\ra$. \epfres
\end{thm}

Note that $\Sing_n=\emptyset$ if $n\leq1$.  Note also that $\Sing_2=\X=\{\ve_{12},\ve_{21}\}$ is a right-zero semigroup, so the smallest (idempotent) generating set for $\Sing_2$ has size $2=|\X|$.  For $n\geq3$, $\Sing_n$ may be generated by $\frac12|\X|$ elements, by further results of Howie, as we now explain.  For a subset $F\sub\X$, let $\Ga_F$ be the directed graph with vertex set $\bn$ and edge set $\set{(i,j)}{\ve_{ij}\in F}$.  We say a directed graph $\Ga$ is \emph{complete} if its associated undirected graph (obtained by changing directed edges to undirected edges, and converting any resulting parallel edges to single edges) is the complete graph (with the same vertex set as $\Ga$).  Part (i) of the next result is \cite[Theorem 2.1]{GH1987}, and Part (ii) is \cite[Theorem 1]{Howie1978}.  Recall that a directed graph is \emph{strongly connected} if, for any pair of vertices $x,y$, there is a (possibly empty) directed path from $x$ to $y$.

\begin{thm}\label{thm:howie2}
Let $n\geq3$.  Then 
\bit
\itemit{i} $\rank(\Sing_n)=\idrank(\Sing_n)=\binom n2=\frac12n(n-1)$\emph{;}
\itemit{ii} if $F\sub\X$, then $\Sing_n=\la F\ra$ if and only if $\Ga_F$ is strongly connected and complete. \epfres
\eit
\end{thm}


A presentation for $\Sing_n$ was given in \cite{JEtnsn2}, in terms of the generating set $\X$.  Define an alphabet
\[
X = \set{e_{ij}}{i,j\in\bn,\ i\not=j},
\]
an epimorphism
\[
\phi:X^+\to\Sing_n:e_{ij}\mt\ve_{ij},
\]
and let $R$ be the set of relations
\begin{align}
\label{R1}\tag{R1} e_{ij}^2 = e_{ij} &= e_{ji}e_{ij} &&\text{for distinct $i,j$}\\
\label{R2}\tag{R2} e_{ij}e_{kl} &= e_{kl}e_{ij} &&\text{for distinct $i,j,k,l$}\\
\label{R3}\tag{R3} e_{ik}e_{jk} &= e_{ik} &&\text{for distinct $i,j,k$}\\
\label{R4}\tag{R4} e_{ij}e_{ik} = e_{ik}e_{ij} &= e_{jk}e_{ij} &&\text{for distinct $i,j,k$}\\
\label{R5}\tag{R5} e_{ki}e_{ij}e_{jk} &= e_{ik}e_{kj}e_{ji}e_{ik} &&\text{for distinct $i,j,k$}\\
\label{R6}\tag{R6} e_{ki}e_{ij}e_{jk}e_{kl} &= e_{ik}e_{kl}e_{li}e_{ij}e_{jl} &&\text{for distinct $i,j,k,l$.}
\end{align}
The next result is \cite[Theorem 6]{JEtnsn2}.

\begin{thm}\label{thm:X}
For $n\geq2$, the semigroup $\Sing_n$ has presentation $\pres{X\!}{\!R}$ via $\phi$. \epfres
\end{thm}


\section{Semidirect products and wreath products}\label{sect:semidirect}

In this section, we prove some general results about semidirect products $M\rtimes S$ and wreath products $M\wr S$, where $M$ is a monoid and $S$ a semigroup ($S$ will be a subsemigroup of some $\T_n$ in the case of $M\wr S$).

\subsection{Semidirect products}\label{subsect:semidirect}

Let $S$ be a semigroup and $M$ a monoid with identity $1$.\footnote{\label{fn:1}In all that follows, $1$ need not be a two-sided identity of $M$, but could instead be a distinguished right identity.}  (Note that $S$ might also be a monoid.)  Suppose $S$ has a left action on $M$ by monoid endomorphisms;  that is, there is a homomorphism $\varphi:S\to\End^*(M):s\mt\varphi_s$, where $\End^*(M)$ denotes the monoid of endomorphisms of $M$ with right-to-left composition.\footnote{\label{fn:2}In the case that $1$ is a distinguished right identity, we still need to assume that every endomorphism $\varphi_s$ fixes $1$.}  For $s\in S$ and $a\in M$, we write $s\cdot a=\varphi_s(a)$.
So
\[
s\cdot1=1 \COMMA s\cdot(t\cdot a) = (st)\cdot a \COMMA s\cdot(ab) = (s\cdot a)(s\cdot b) \qquad\text{for all $s,t\in S$ and $a,b\in M$.}
\]
(Note that if $S$ happens to be a monoid, we do not assume that it acts \emph{monoidally} on $M$: i.e., we do not assume that the identity of $S$ acts as the identity automorphism of $M$.)
The \emph{semidirect product} $M\rtimes S=M\rtimes_\varphi S$ has underlying set $M\times S=\set{(a,s)}{a\in M,\ s\in S}$, and product defined by
\[
(a,s)(b,t) = (a(s\cdot b),st)  \qquad\text{for all $s,t\in S$ and $a,b\in M$.}
\]
The fact that $S$ acts by \emph{monoid} endomorphisms ensures that $S$ may be identified with the subsemigroup $\set{(1,s)}{s\in S}$ of $M\rtimes S$.  If $S$ is a monoid acting monoidally on $M$ (i.e., $1\cdot a=a$ for all $a\in M$), then $\set{(a,1)}{a\in M}$ is an isomorphic copy of $M$ inside $M\rtimes S$.  However, this article is mostly concerned with the case that $S$ is not a monoid, in which case $M\rtimes S$ does not contain such a canonical copy of $M$.
%
%
A motivating example of the semidirect product are the \emph{wreath products}, which are the subject of Section \ref{subsect:wreath}.

Suppose now that $S$ has semigroup presentation $\pres{X\!}{\!R}$ via $\phi:X^+\to S$.  Define a new alphabet $X_M=\set{x_a}{x\in X,\ a\in M}$.  We regard $X$ as a subset of $X_M$ by identifying $x\in X$ with $x_1\in X_M$.  For a word $w=x_1\cdots x_k\in X^+$, and for an element $a\in M$, we define the word $w_a=(x_1)_a x_2\cdots x_k\in X_M^+$.  Consider the set $R_M=R_M^1\cup R_M^2$ of relations over $X_M$, where $R_M^1$ and $R_M^2$ are defined by
\[
R_M^1=\set{(u_a,v_a)}{(u,v)\in R,\ a\in M} \AND R_M^2 = \set{(x_ay_b,x_{a(x\phi\cdot b)}y)}{x,y\in X,\ a,b\in M}.
\]
Note that, by identifying $X\sub X_M$ as above, we also have $R\sub R_M$, via $(u,v)\equiv(u_1,v_1)$.
Define a map
\[
\phi_M:X_M^+\to M\rtimes S \qquad\text{by}\qquad x_a\phi_M = (a,x\phi) \qquad\text{for all $x\in X$ and $a\in M$.}
\]
It is easy to check that $w_a\phi_M=(a,w\phi)$ for all $a\in M$ and $w\in X^+$.  It quickly follows (from the surjectivity of $\phi:X^+\to S$) that $\phi_M$ is surjective.
For convenience, in the following proof, even though $R$ and $R_M$ may not be symmetric, we will say ``$(u,v)\in R$'' to mean ``$(u,v)\in R$ or $(v,u)\in R$'' (and similarly for $R_M$).

\begin{thm}\label{thm:semidirect}
With the above notation, $M\rtimes S$ has semigroup presentation $\pres{X_M}{R_M}$ via $\phi_M$.
\end{thm}

\pf We showed above that $\phi_M$ is surjective, so we just need to show that $\ker(\phi_M)=R_M^\sharp$.
%
%
First note that for any $(u,v)\in R$ and $a\in M$, 
$
u_a\phi_M=(a,u\phi)=(a,v\phi)=v_a\phi_M,
$
while for any $x,y\in X$ and $a,b\in M$,
$
(x_ay_b)\phi_M = (a,x\phi)(b,y\phi) = (a(x\phi\cdot b),(xy)\phi) = (a(x\phi\cdot b),x\phi)(1,y\phi) = (x_{a(x\phi\cdot b)}y)\phi_M,
$
showing that $R_M\sub\ker(\phi_M)$.  

Conversely, suppose $u = (x_1)_{a_1}\cdots (x_k)_{a_k},v = (y_1)_{b_1}\cdots(y_l)_{b_l}\in X_M^+$ are such that $u\phi_M=v\phi_M$.  For the remainder of the proof, write ${\approx}=R_M^\sharp$.  Using relations from $R_M^2$, we have
\[
u \approx (x_1)_ax_2\cdots x_k=(x_1\cdots x_k)_a \AND v\approx (y_1)_by_2\cdots y_l=(y_1\cdots y_l)_b \qquad\text{for some $a,b\in M$.}
\]
Since ${\approx}\sub\ker(\phi_M)$, we have
\[
(a,(x_1\cdots x_k)\phi)=(x_1\cdots x_k)_a\phi_M=u\phi_M=v\phi_M=(y_1\cdots y_l)_b\phi_M = (b,(y_1\cdots y_l)\phi).
\]
It follows that $a=b$ and $(x_1\cdots x_k)\phi=(y_1\cdots y_l)\phi$.  Since $\ker(\phi)=R^\sharp$, it follows that there is a sequence of words
$
x_1\cdots x_k=w_0, w_1,\ldots,w_r=y_1\cdots y_l
$
such that, for each $0\leq i\leq r-1$, $w_i=w_i'uw_i''$ and $w_{i+1}=w_i'vw_i''$ for some $w_i',w_i''\in X^*$ and $(u,v)\in R$.  But then we see that $(w_i)_a\approx (w_{i+1})_a$, using either $(u,v)\in R\sub R_M$ (if $w_i'$ is non-empty) or $(u_a,v_a)\in R_M$ (if $w_i'$ is empty).  But then
\[
u\approx (x_1\cdots x_k)_a =(w_0)_a \approx (w_1)_a \approx\cdots\approx (w_r)_a=(y_1\cdots y_l)_a=(y_1\cdots y_l)_b\approx v,
\]
completing the proof. \epf

%


The next result follows immediately from Theorem \ref{thm:semidirect}

\begin{cor}\label{cor:semidirect}
If $S$ is finitely presented and $M$ is finite, then $M\rtimes S$ is finitely presented.~\epfres
\end{cor}

Of course the converse of Corollary \ref{cor:semidirect} is not true: for example, the semidirect product $M\rtimes S$ of finitely presented infinite \emph{monoids} (with a monoidal action of $S$ on $M$) is finitely presented \cite[Corollary 2]{Lavers1998}.  One might hope to improve Corollary \ref{cor:semidirect} by assuming only that $M$ is finitely presented (rather than the stronger assumption of being finite).  But this is not the case in general, as the following examples show.

\begin{eg}\label{eg:1}
Let $\N$ be the additive monoid of natural numbers (including $0$), and let $S=\{x,x^2\}$ with $x^3=x^2\not=x$.  It is easy to show that any generating set for the \emph{direct} product $\N\times S$ must contain $\N\times\{x\}$.
In particular, $\N\times S$ is not finitely generated (and, hence, not finitely presented).
\end{eg}

\begin{eg}\label{eg:2}
Our second example is a semidirect product that is not direct.  Let $M=\N\times\N$, and let $S=\{\ve\}$ be a trivial semigroup.  Define an action of $S$ on $M$ by $\ve\cdot(a,b)=(a,a)$.  Since $|S|=1$, we may identify the element $((a,b),\ve)$ of $M\rtimes S$ with $(a,b)$: in this way, the operation on $M\rtimes S$ obeys $(a,b)(c,d)=(a+c,b+c)$.  
It is again easy to show that any generating set for $M\rtimes S$ must contain $\N\times\{0\}$.  (In fact, this example can be viewed as a \emph{wreath product}; see Section \ref{subsect:wreath} for more details.)
\end{eg}

%

%

\subsection{Wreath products}\label{subsect:wreath}

Let $S$ be a subsemigroup of the full transformation semigroup $\T_n$, and let $M$ be an arbitrary monoid.  Then~$S$ has a natural left action on $M^n$ (the direct product of $n$ copies of $M$) given by
\[
\al\cdot(a_1,\ldots,a_n)=(a_{1\al},\ldots,a_{n\al}) \qquad\text{for $\al\in S$ and $a_1,\ldots,a_n\in M$,}
\]
and the identity element $(1,\ldots,1)$ of $M^n$ is clearly fixed by every $\al\in S$.\footnote{Note that if $1$ is \emph{any} right identity of $M$, then $(1,\ldots,1)$ is a right identity of $M^n$ and is still fixed under the action of $S$; cf.~Footnotes \ref{fn:1} and \ref{fn:2}.}
The resulting semidirect product $M^n\rtimes S$ is the \emph{wreath product} of $M$ by $S$, denoted $M\wr S$.  
Multiplication in $M\wr S$ obeys the rule
\[
((a_1,\ldots,a_n),\al)((b_1,\ldots,b_n),\be) = ((a_1b_{1\al},\ldots,a_nb_{n\al}),\al\be).
\]
For example, writing $\ve=\ve_{12}\in\T_2$, if $M=\N$ and $S=\{\ve\}$, then $M\wr S$ is the semigroup from Example \ref{eg:2}.  When $S=\T_n$, we obtain the \emph{full} wreath product $M\wr\T_n$.  When $S=\Sing_n=\TnSn$, we obtain the \emph{singular} wreath product $M\wr\Sing_n$; these singular wreath products are the main focus of this article.  If $M=\{1\}$, then $M\wr S\cong S$ for any $S$.  On the other hand, if $S=\{1\}$, where $1\in\T_n$ denotes the identity map, then $M\wr S\cong M^n$.

\begin{figure}[h]
\begin{center}
\scalebox{.9}{
\begin{tikzpicture}[scale=1]
\tikzstyle{vertex}=[circle,draw=black, fill=white, inner sep = 0.06cm]
\draw[-{latex}](6+1,-.5)--(10-1,-.5);
\draw[-{latex}](12.5,-2.5-1)--(12.5,-5.5+1);
\draw[{latex}-](6+1,-7.5)--(10-1,-7.5);
\begin{scope}[shift={(0,0)}]	
\lvs{1,...,6}
\stlines{1/2,2/1,3/3,4/2,5/6,6/6}
\draw(.5,1)node[left]{$(\ba,\al)=$};
\foreach \x in {1,...,6} {\node[vertex] () at (\x,2){$a_{\x}$};}
\end{scope}
\begin{scope}[shift={(0,-3)}]	
\lvs{1,...,6}
\stlines{1/1,2/1,3/4,4/4,5/5,6/4}
\draw(.5,1)node[left]{$(\bb,\be)=$};
\foreach \x in {1,...,6} {\node[vertex] () at (\x,2){$b_{\x}$};}
\end{scope}
\begin{scope}[shift={(9,-.5)}]	
\lvs{1,...,6}
\stlines{1/2,2/1,3/3,4/2,5/6,6/6}
\foreach \x in {1,...,6} {\node[vertex] () at (\x,2){$a_{\x}$};}
\end{scope}
\begin{scope}[shift={(9,-2.5)}]	
\lvs{1,...,6}
\stlines{1/1,2/1,3/4,4/4,5/5,6/4}
\foreach \x in {1,...,6} {\node[vertex] () at (\x,2){$b_{\x}$};}
\end{scope}
\begin{scope}[shift={(9,-7-.5)}]	
\stlines{1/2,2/1,3/3,4/2,5/6,6/6}
\foreach \x/\y in {1/2,2/1,3/3,4/2,5/6,6/6} {\node[vertex] () at (\x,2){{\footnotesize $a_{\x}b_{\y}$}};}
\end{scope}
\begin{scope}[shift={(9,-7-2.5)}]	
\lvs{1,...,6}
\stlines{1/1,2/1,3/4,4/4,5/5,6/4}
\draw(6.5,1)node[right]{$\phantom{=(\ba,\al)(\bb,\be)}$};
\end{scope}
\begin{scope}[shift={(0,-8.5)}]	
\lvs{1,...,6}
\stlines{1/1,2/1,3/4,4/1,5/4,6/4}
\draw(.5,1)node[left]{$(\ba,\al)(\bb,\be)=$};
\foreach \x/\y in {1/2,2/1,3/3,4/2,5/6,6/6} {\node[vertex] () at (\x,2){{\footnotesize $a_{\x}b_{\y}$}};}
\end{scope}
\end{tikzpicture}
}
\end{center}
\vspace{-5mm}
\caption{Elements of $M\wr\T_6$ (top left) and their product (bottom left).}
\label{fig:MwrT6}
\end{figure}
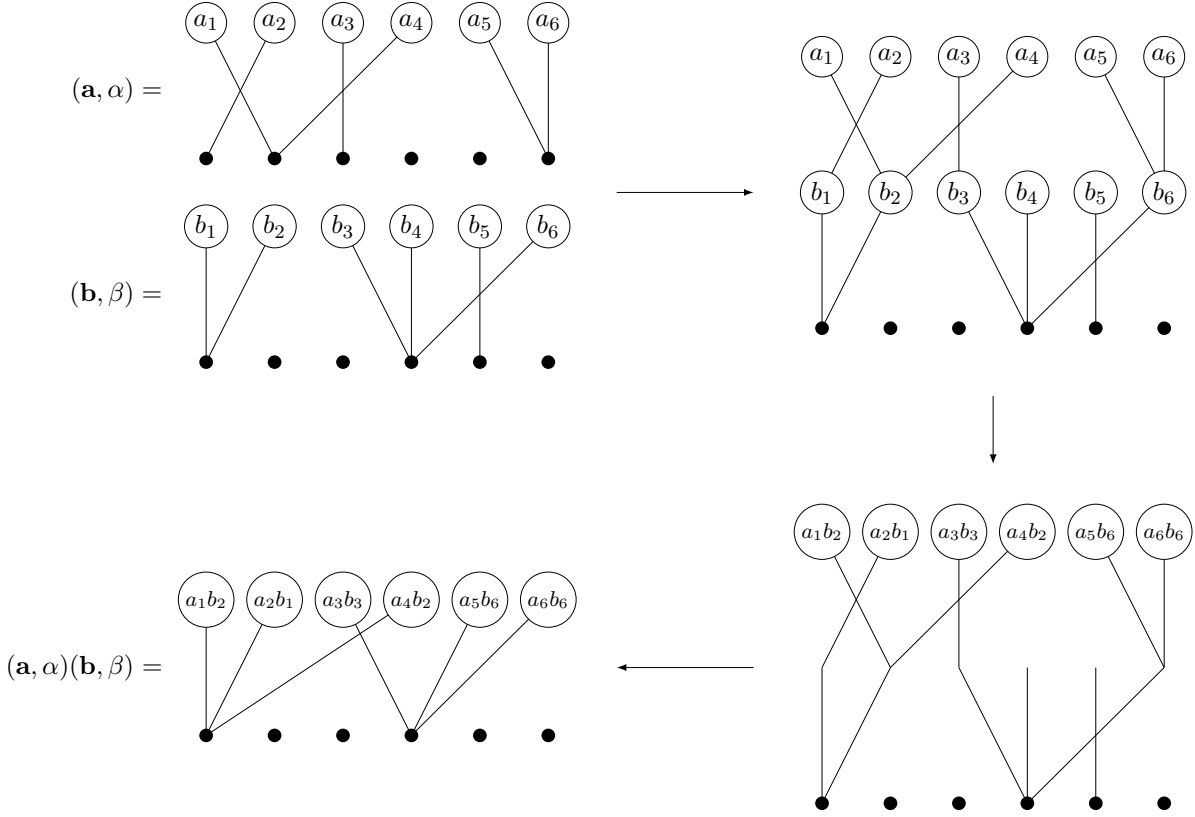

There is a useful way to picture an element $(\ba,\al)$ of $M\wr\T_n$.  We draw two parallel rows of vertices, both labelled $1,\ldots,n$ (and assumed to be increasing from left to right, unless otherwise specified); for each $i\in\bn$, we draw an edge between upper vertex $i$ and lower vertex $i\al$; and we decorate upper vertex $i$ with the monoid element $a_i$, where $\ba=(a_1,\ldots,a_n)$.  For example, two elements $(\ba,\al),(\bb,\be)\in M\wr\T_6$ are pictured in Figure~\ref{fig:MwrT6}.  To calculate the product $(\ba,\al)(\bb,\be)=(\ba(\al\cdot\bb),\al\be)$ diagrammatically: we first stack $(\ba,\al)$ above $(\bb,\be)$, identifying lower vertex $i$ of $(\ba,\al)$ with upper vertex $i$ of $(\bb,\be)$ for each $i$; we then ``slide'' the $b_i$ decorations up the edges of $(\ba,\al)$ and form the appropriate products $a_jb_i$; finally, we straighten the remaining edges, and remove incomplete edges.  An example calculation is also given in Figure \ref{fig:MwrT6}.  By convention, we will often omit the decoration on upper vertex $i$ of $(\ba,\al)$ if $a_i=1$.

Since $M\wr S=M^n\rtimes S$ is a semidirect product, Theorem \ref{thm:semidirect} directly leads to a general presentation for~$M\wr S$, modulo a presentation $\pres{X\!}{\!R}$ for $S$, but we will not state this explicitly.

The remainder of the article almost exclusively concerns singular wreath products $M\wr\Sing_n$.  Because $\Sing_n$ is empty for $n\leq1$, we will assume that $n\geq2$ whenever we make a statement about $\Sing_n$.


\boldsection{Idempotents and generators for $M\wr\Sing_n$}\label{sect:gen}

In this section, we extend Proposition \ref{prop:ETn} to full wreath products $M\wr\T_n$ and Theorems \ref{thm:howie} and \ref{thm:howie2} to singular wreath products $M\wr\Sing_n$.  In particular, we characterise and enumerate the idempotents of~${M\wr\T_n}$ (Lemma~\ref{lem:idempotents}, Proposition \ref{prop:idempotents} and Corollary \ref{cor:idempotents}); we give canonical generators and re-prove Bulman-Fleming's~\cite{BF1995} necessary and sufficient conditions for $M\wr\Sing_n$ to be idempotent-generated (Theorem \ref{thm:IG}); and we give bounds (and exact values in some cases) on the minimal size of (idempotent) generating sets for $M\wr\Sing_n$ (Proposition~\ref{prop:rank1} and Theorems \ref{thm:rank2} and \ref{thm:infinite}).

For $i,j\in\bn$ with $i\not=j$, and for $\ba\in M^n$, we define $\ve_{ij;\ba}=(\ba,\ve_{ij})\in M\wr\Sing_n$.  
%
As a special case, for $a,b\in M$, we define
\[
\ve_{ij;ab} = \ve_{ij;\ba} \qquad\text{where $\ba\in M^n$ is defined by $a_k = \begin{cases}
a &\text{if $k=i$}\\
b &\text{if $k=j$}\\
1 &\text{otherwise.}
\end{cases}$}
\]
As a special case of the latter, we define $\ve_{ij;a}=\ve_{ij;1a}$, for $a\in M$.  We gather these elements into the sets
%
\begin{align*}
\X_n &= \set{\ve_{ij;\ba}}{i,j\in\bn,\ i\not=j,\ \ba\in M^n}, \\
\X_2 &= \set{\ve_{ij;ab}}{i,j\in\bn,\ i\not=j,\ a,b\in M}, \\
\X_1 &= \set{\ve_{ij;a}}{i,j\in\bn,\ i\not=j,\ a\in M}.
\end{align*}
As usual, we also identify $\ve_{ij}\in\Sing_n$ with $\ve_{ij;1}\in M\wr\Sing_n$.  So we have $\X\sub\X_1\sub\X_2\sub\X_n$.  Note that
\[
|\X|=2\tbinom n2 \COMMA |\X_1|=2|M|\times \tbinom n2 \COMMA |\X_2|=2|M|^2\times \tbinom n2 \COMMA |\X_n|=2|M|^n\times \tbinom n2 .
\]
Note also that $\X_1\sub E(M\wr\Sing_n)$, as we show in Figure \ref{fig:eija} (and also follows from Lemma \ref{lem:idempotents} below), but that elements of $\X_2$ need not be idempotent in general.

\begin{figure}[h]
\begin{center}
\begin{tikzpicture}[scale=1]
\tikzstyle{vertex}=[circle,draw=black, fill=white, inner sep = 0.06cm]
\begin{scope}[shift={(0,0)}]	
\uvs{1,2,3,5}
\lvs{1,2,3,5}
\stlines{1/1,2/1,3/3,5/5}
\dotsups{3/5}
\dotsdns{3/5}
\vertlabelsh{1/i,2/j,3/i_1,5/i_{n-2}}
\draw(0,1)node[left]{$\ve_{ij;a}$};
\draw[|-|] (0,0)--(0,2);
\node[vertex] () at (2,2){$a$};
\end{scope}
\begin{scope}[shift={(0,-2)}]	
\uvs{1,2,3,5}
\lvs{1,2,3,5}
\stlines{1/1,2/1,3/3,5/5}
\dotsups{3/5}
\dotsdns{3/5}
\draw(0,1)node[left]{$\ve_{ij;a}$};
\draw[|-|] (0,0)--(0,2);
\node[vertex] () at (2,2){$a$};
\end{scope}
\begin{scope}[shift={(8,-1)}]	
\uvs{1,2,3,5}
\lvs{1,2,3,5}
\stlines{1/1,2/1,3/3,5/5}
\dotsups{3/5}
\dotsdns{3/5}
\vertlabelsh{1/i,2/j,3/i_1,5/i_{n-2}}
\draw(6,1)node[right]{$\ve_{ij;a}$};
\draw[|-|] (6,0)--(6,2);
\node[vertex] () at (2,2){$a$};
\draw(-1,1)node{$=$};
\end{scope}
\end{tikzpicture}
\end{center}
\vspace{-5mm}
\caption{Diagrammatic proof that $\ve_{ij;a}$ is an idempotent.  Vertices are drawn in the order $i,j,i_1,\ldots,i_{n-2}$, where $\bn\sm\{i,j\}=\{i_1,\ldots,i_{n-2}\}$.}
\label{fig:eija}
\end{figure}
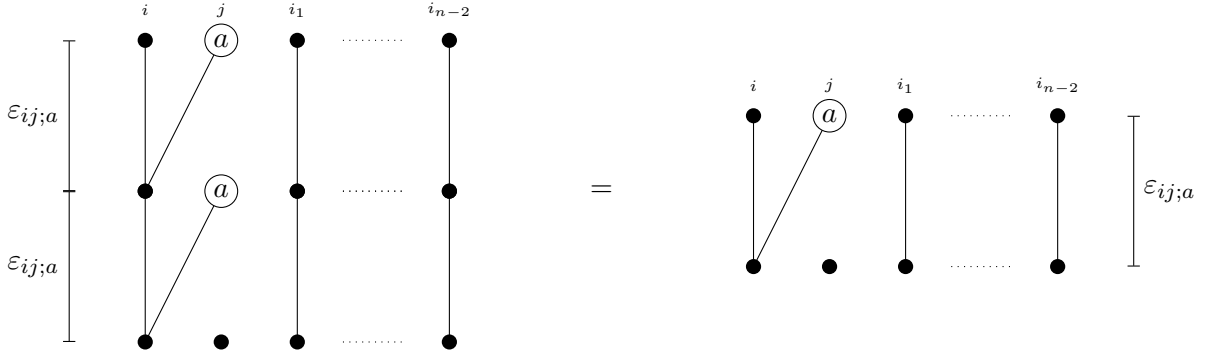

\boldsubsection{Idempotents in $M\wr\T_n$}

We begin with a simple description of the idempotents of $M\wr\T_n$.

\begin{lemma}\label{lem:idempotents}
We have $E(M\wr\T_n) = \set{(\ba,\al)\in M\wr\T_n}{\al\in E(\T_n) \text{ and } a_ia_{i\al}=a_i \text{ for all } i\in\bn}$.
\end{lemma}

\pf This follows immediately from the fact that $(\ba,\al)^2=((a_1a_{1\al},\ldots,a_na_{n\al}),\al^2)$. \epf

\begin{rem}\label{rem:idempotents}
If $(\ba,\al)\in E(M\wr\T_n)$, then $a_i\in E(M)$ for all $i\in\im(\al)$, since $\al\in E(\T_n)$ acts as the identity on $\im(\al)$.
\end{rem}

If $M$ is infinite, then $|E(M\wr\T_n)|=|M|=|M\wr\T_n|$, since $|\X_1|=|M|$ and $\X_1\sub E(M\wr\T_n)$.  The next result calculates the number of idempotents in $M\wr\T_n$ when $M$ is finite.  

\begin{prop}\label{prop:idempotents}
If $M$ is a finite monoid, then
\[
|E(M\wr\T_n)|=
\sum_{k=1}^n \binom nk \sum_{e_1,\ldots,e_k \atop \in E(M)} \left( \sum_{i=1}^k|Me_i|\right)^{n-k}.
\]
\end{prop}

%

\pf To specify an idempotent $(\ba,\al)\in E(M\wr\T_n)$, we first choose $\im(\al)$, of size $k$ (say), in $\binom nk$ ways.  Write $\im(\al)=\{j_1,\ldots,j_k\}$, where $j_1<\cdots<j_k$.  We then choose idempotents $a_{j_1}=e_1,\ldots,a_{j_k}=e_k\in E(M)$; see Remark \ref{rem:idempotents}.  We then choose the preimages $j_1\al^{-1},\ldots,j_k\al^{-1}$, of sizes $l_1+1,\ldots,l_k+1$ (say), in $\binom{n-k}{l_1,\ldots,l_k}$ ways, noting that $j_i\al=j_i$ for each $i$.  For each $1\leq i\leq k$, and for each $q\in j_i\al^{-1}\sm\{j_i\}$, we must choose $a_q\in M$ so that $a_qe_i=a_q$: i.e., $a_q\in Me_i$.  So there are $|Me_i|^{l_i}$ ways to choose $\set{a_q}{q\in j_i\al^{-1}\sm\{j_i\}}$.  
Multiplying and adding these values, as appropriate, gives 
\[
|E(M\wr\T_n)|=\sum_{k=1}^n \binom nk \sum_{e_1,\ldots,e_k \atop \in E(M)} \sum_{(l_1,\ldots,l_k)\in\N^k\atop l_1+\cdots+l_k=n-k}\binom{n-k}{l_1,\ldots,l_k}\prod_{i=1}^k|Me_i|^{l_i}.
\]
The multinomial formula then completes the proof. \epf


\begin{rem}
In the case that $M=\T_m$ for some $m$, Proposition \ref{prop:idempotents} specialises to a much simpler version of \cite[Proposition 3.1]{DE2016}; the latter gives a formula for $|E(\T_m\wr\T_n)|$ that bears more similarity to the displayed equation in the proof of Proposition \ref{prop:idempotents}.
A formula for $|E(M\wr\Sing_n)|$ may be obtained from $|E(M\wr\T_n)|$ by subtracting $|E(M)|^n$ (i.e., the $k=n$ term in the sum from Proposition \ref{prop:idempotents}).%
\end{rem}

Proposition \ref{prop:idempotents} simplifies substantially in the case that the finite monoid $M$ has only a single idempotent; this occurs precisely when $M$ is a finite group.

\begin{cor}\label{cor:idempotents}
If $G$ is a finite group, then $\ds{|E(G\wr\T_n)|= \sum_{k=1}^n \binom nk  k^{n-k}|G|^{n-k}}$. \epfres
\end{cor}

%

\begin{rem}
%
In the case that $|G|=1$, $G\wr\T_n\cong\T_n$, and Corollary \ref{cor:idempotents} reduces to Proposition~\ref{prop:ETn}.
\end{rem}


\boldsubsection{Generation and idempotent-generation in $M\wr\Sing_n$}

The statement and proof of the next result refer to the sets $\X\sub\X_1\sub\X_2\sub\X_n\sub M\wr\Sing_n$ defined at the beginning of Section \ref{sect:gen}, and to the preorder $\leq_{\L}$ and relation $\L$ defined (on $M$) as in Section \ref{sect:prelim}.  Recall that $M/{\L}$ denotes the partially ordered set of all $\L$-classes of $M$.  We denote by $\de$ the diagonal equivalence $\set{(i,i)}{i\in\bn}$.  Parts (ii) and (iii) of the next result are contained in \cite[Theorems 5.5 and 5.7]{BF1995}, but they are crucial in what follows, and we provide short proofs for convenience.  We note that the corresponding statements from \cite{BF1995} are dual to ours, because we compose maps from left to right.  Note also that 
\begin{align*}
\text{$M/{\L}$ is a chain} &\iff \text{the principal left ideals of $M$ form a chain under inclusion}\\
&\iff \text{all finitely generated left ideals of $M$ are principal.}
\end{align*}

\begin{thm}\label{thm:IG}
If $M$ is any monoid, then
\bit
\itemit{i} $M\wr\Sing_n = \la\X_n\ra=\la\X_2\ra$\emph{;}
\itemit{ii} $\la E(M\wr\Sing_n)\ra = \la\X_1\ra = \set{(\ba,\al)\in M\wr\Sing_n}{a_i\leq_{\L} a_j\text{ for some } (i,j)\in\ker(\al)\sm\de}$\emph{;} and
\itemit{iii} $M\wr\Sing_n$ is idempotent-generated if and only if $M/{\L}$ is a chain.
\end{itemize}
\end{thm}

\begin{figure}[h]
\begin{center}
\begin{tikzpicture}[scale=1]
\tikzstyle{vertex}=[circle,draw=black, fill=white, inner sep = 0.06cm]
\begin{scope}[shift={(0,-3)}]	
\uvs{1,2,3,5}
\lvs{1,2,3,5}
\stlines{1/1,2/1,3/3,5/5}
\dotsups{3/5}
\dotsdns{3/5}
\draw(0,1)node[left]{$\ve_{12;\ba}$};
\draw[|-|] (0,0)--(0,2);
\foreach \x/\y in {1/1,2/2} {\node[vertex] () at (\x,2){$a_{\y}$};}
\foreach \x/\y in {3/3,5/n} {\node[vertex] () at (\x,2){$a_{\y}$};}
\draw(5.5+1.5,1)node{$=$};
\end{scope}
\begin{scope}[shift={(8,0)}]	
\uvs{1,2,3,5}
\lvs{1,2,3,5}
\stlines{1/1,2/1,3/3,5/5}
\dotsups{3/5}
\dotsdns{3/5}
\draw(6,1)node[right]{$\ve_{12;a_1,a_2}$};
\draw[|-|] (6,0)--(6,2);
\foreach \x/\y in {1/1,2/2} {\node[vertex] () at (\x,2){$a_{\y}$};}
\end{scope}
\begin{scope}[shift={(8,-2)}]	
\uvs{1,2,3,5}
\lvs{1,2,3,5}
\stlines{1/1,2/3,3/3,5/5}
\dotsups{3/5}
\dotsdns{3/5}
\draw(6,1)node[right]{$\ve_{32;a_3,1}$};
\draw[|-|] (6,0)--(6,2);
\foreach \x/\y in {3/3} {\node[vertex] () at (\x,2){$a_{\y}$};}
\foreach \x in {1,2,3,5} {\draw[dotted](\x,0)--(\x,-2);}
\end{scope}
\begin{scope}[shift={(8,-6)}]	
\uvs{1,2,3,5}
\lvs{1,2,3,5}
\stlines{1/1,2/5,3/3,5/5}
\dotsups{3/5}
\dotsdns{3/5}
\draw(6,1)node[right]{$\ve_{n2;a_n,1}$};
\draw[|-|] (6,0)--(6,2);
\foreach \x/\y in {5/n} {\node[vertex] () at (\x,2){$a_{\y}$};}
\end{scope}
\end{tikzpicture}
\end{center}
\vspace{-5mm}
\caption{Diagrammatic proof that  that $\ve_{12;\ba} = \ve_{12;a_1,a_2}\ve_{32;a_3,1}\cdots\ve_{n2;a_n,1}$.}
\label{fig:eijba}
\end{figure}
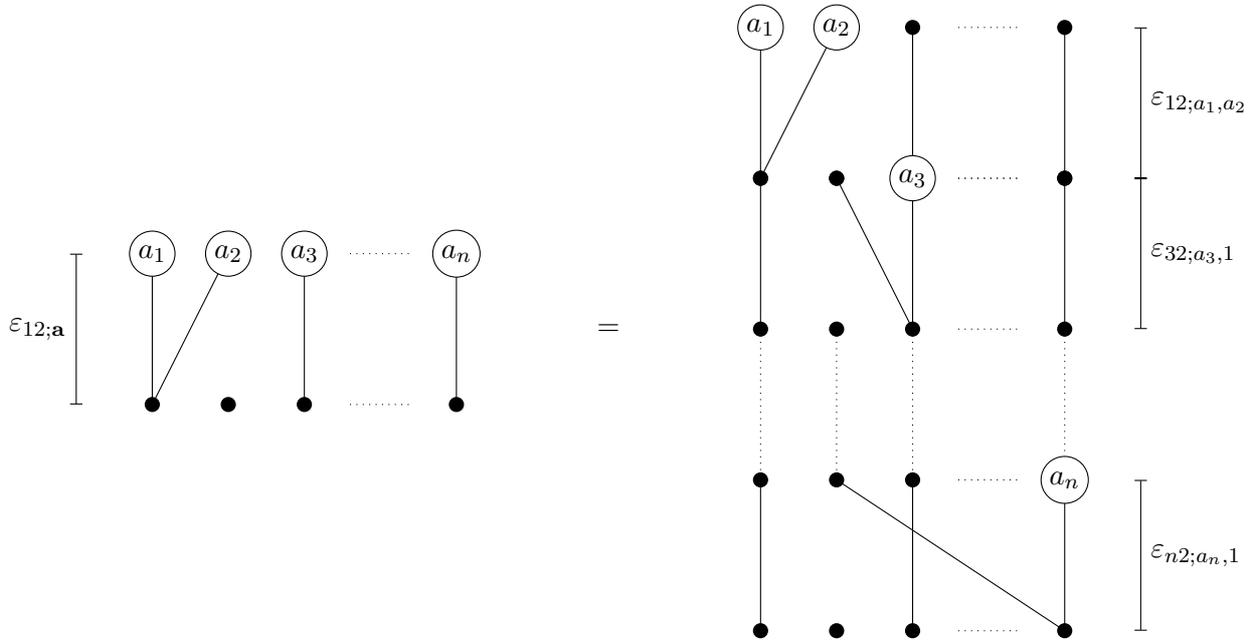

\pf (i).  Theorems \ref{thm:howie} and \ref{thm:semidirect} quickly give $M\wr\Sing_n=\la\X_n\ra$.  To complete the proof of (i), since $\X_2\sub\X_n$, it suffices to prove that $\X_n\sub\la\X_2\ra$, so let $\ve_{ij;\ba}\in\X_n$.  Relabelling the elements of $\bn$ if necessary, we may assume that $(i,j)=(1,2)$.  We show in Figure \ref{fig:eijba} that $\ve_{12;\ba} = \ve_{12;a_1,a_2}\ve_{32;a_3,1}\cdots\ve_{n2;a_n,1}$.

\pfitem{ii}  Put $\Si=\set{(\ba,\al)\in M\wr\Sing_n}{a_i\leq_{\L} a_j\text{ for some } (i,j)\in\ker(\al)\sm\de}$.  Since $\X_1\sub E(M\wr\Sing_n)$, it suffices to show that:
\bit\begin{multicols}{2}
\item[(a)] $\la E(M\wr\Sing_n)\ra\sub\Si$, and 
\item[(b)] $\Si\sub\la\X_1\ra$.
\end{multicols}\eit
For (a), let 
$(\ba,\al)\in\la E(M\wr\Sing_n)\ra$.  So $(\ba,\al)=(\bb,\be)(\ba,\al)=(\bb(\be\cdot\ba),\be\al)$ for some $(\bb,\be)\in E(M\wr\Sing_n)$.  Since $\be\in\Sing_n$, we may choose some $i\in\bn\sm\im(\be)$.  Put $j=i\be$, noting also that $j=j\be$ (since $\be$ is an idempotent).  Examining the $i$th coordinate of $\ba=\bb(\be\cdot\ba)$ gives $a_i=b_ia_{i\be}=b_ia_j\leq_{\L} a_j$.  Since $(i,j)\in\ker(\be)\sub\ker(\be\al)=\ker(\al)$, this completes the proof of (a).

For (b), let 
$(\ba,\al)\in\Si$.  Relabelling the elements of $\bn$ if necessary, we may assume that $(1,2)\in\ker(\al)$ and $a_1\leq_{\L} a_2$, so that $a_1=xa_2$ for some $x\in M$.  Note that $(\ba,\al)=\ve_{12;\ba}\al$ (thinking of $\al\in\Sing_n$ as an element of $M\wr\Sing_n$).  Since $\Sing_n=\la\X\ra\sub\la\X_1\ra$, it remains to show that $\ve_{12;\ba}\in\la\X_1\ra$.  As in Part (i), we have $\ve_{12;\ba}=\ve_{12;a_1,a_2}\ve_{32;a_3,1}\cdots\ve_{n2;a_n,1}$, and the proof of (b) is complete, after noting that
\[
\ve_{12;a_1,a_2} = \ve_{21;x}\ve_{12;a_2}  \AND  \ve_{k2;a_k,1} = \ve_{2k;a_k}\ve_{k2;1} \qquad\text{for all $3\leq k\leq n$.}
\]

\pfitem{iii}  Suppose first that $M/{\L}$ is a chain.  
Consider an arbitrary element $(\ba,\al)\in M\wr\Sing_n$, and choose any $(i,j)\in\ker(\al)\sm\de$.  Then either $a_i\leq_{\L} a_j$ or $a_j\leq_{\L} a_i$ (since $M/{\L}$ is a chain).  Part (ii) then gives $(\ba,\al)\in\la E(M\wr\Sing_n)\ra$, showing that $M\wr\Sing_n$ is idempotent-generated.  Conversely, suppose $M\wr\Sing_n$ is idempotent-generated, and let $a,b\in M$ be arbitrary.  Then $\ve_{12;ab}\in\Si$, by (ii).  As we clearly have $\ker(\ve_{12})\sm\de=\{(1,2),(2,1)\}$, it follows that $a\leq_{\L} b$ or $b\leq_{\L} a$.  This completes the proof of (iii). \epf

\begin{rem}
Recall from Theorem \ref{thm:howie} that $\Sing_n$ is generated by its idempotents of rank $n-1$.  Parts~(ii) and (iii) of Theorem \ref{thm:IG} show that when $M\wr\Sing_n$ is idempotent-generated (i.e., when $M/{\L}$ is a chain), it is generated by idempotents whose underlying (idempotent) transformation has rank $n-1$.  See also \cite[Theorem 4.2]{FG2007} and \cite[Theorem 2.1]{FL1992}.
\end{rem}

The elements of $\la E(M\wr\T_n)\ra$ were characterised in \cite[Theorem 5.6]{BF1995}; the description is a little more technical than Theorem \ref{thm:IG}(ii).  The next result (which was not proved in \cite{BF1995}) characterises the monoids $M$ for which $\la E(M\wr\T_n)\ra$ takes on a particularly simple form.  For its statement, recall that we write $1$ for the identity element of any monoid (including $M\wr\T_n$).

\begin{cor}\label{cor:IGMwrTn}
For any monoid $M$, the following are equivalent:
\bit
\itemit{i} $\la E(M\wr\T_n)\ra=\{1\}\cup(M\wr\Sing_n)$\emph{;} and
\itemit{ii} $M/{\L}$ is a chain and $E(M)=\{1\}$.
\eit
\end{cor}

\pf Suppose first that (i) holds.  If there was a non-identity idempotent $a\in E(M)\sm\{1\}$, then $((a,1,\ldots,1),1)\in E(M\wr\T_n)\sub\{1\}\cup(M\wr\Sing_n)$, a contradiction.  It follows that $E(M)=\{1\}$.  Now suppose $a,b\in M$ are arbitrary.  By assumption, we may write $\ve_{12;ab}=(\ba_1,\al_1)\cdots(\ba_k,\al_k)$ for some idempotents $(\ba_1,\al_1),\ldots,(\ba_k,\al_k)\in E(M\wr\T_n)$.  If any $\al_i=1$, then Lemma \ref{lem:idempotents} and $E(M)=\{1\}$ gives $\ba_i=(1,\ldots,1)$; in this case, the factor $(\ba_i,\al_i)$ is the identity of $M\wr\T_n$, and so may be cancelled from the product.  After cancelling all such trivial factors, we deduce that $\ve_{12;ab}\in\la E(M\wr\Sing_n)\ra$.  Since $\ker(\ve_{12})\sm\de=\{(1,2),(2,1)\}$, Theorem \ref{thm:IG}(ii) gives $a\leq_{\L} b$ or $b\leq_{\L} a$.  So $M/{\L}$ is a chain.

Conversely, suppose (ii) holds.  First note that ${\{1\}\cup(M\wr\Sing_n)=\{1\}\cup\la E(M\wr\Sing_n)\ra\sub\la E(M\wr\T_n)\ra}$, by Theorem \ref{thm:IG}(iii).  To establish the reverse inclusion, suppose $(\ba_1,\al_1),\ldots,(\ba_k,\al_k)\in E(M\wr\T_n)$, and put $(\ba,\al)=(\ba_1,\al_1)\cdots(\ba_k,\al_k)$.  Again, $E(M)=\{1\}$ tells us that any factor with $\al_i=1$ is the identity of $M\wr\T_n$.  If $\al_i=1$ for all $i$, then $(\ba,\al)=1$; otherwise, we cancel all such trivial factors and deduce that ${(\ba,\al)\in\la E(M\wr\Sing_n)\ra}$. \epf

\begin{rem}
If $M$ is a group, then clearly $M$ satisfies the conditions of Corollary \ref{cor:IGMwrTn}(ii).
\end{rem}

\boldsubsection{Rank and idempotent rank in $M\wr\Sing_n$}

Recall that the rank (resp., idempotent rank) of a semigroup (resp., idempotent-generated semigroup), denoted $\rank(S)$ (resp., $\idrank(S)$), is the minimum size of a generating set (resp., idempotent generating set) for $S$.  In this section, we prove a number of results concerning the rank and idempotent rank of $M\wr\Sing_n$.

For $\oijn$, write $\eta_{ij}$ for the equivalence relation on $\bn$ with unique non-trivial equivalence class $\{i,j\}$.  In particular, $\eta_{ij}=\ker(\ve_{ij})=\ker(\ve_{ji})$.  For several proofs in this section, recall that a directed graph~$\Ga$ with no loops is a \emph{tournament} if, for every pair of distinct vertices $x,y$, $\Ga$ has exactly one of the edges $(x,y)$ or $(y,x)$.  Note that a tournament on vertex set $\bn$ has $\binom n2$ edges if $n\geq3$ (and that there are no strongly connected tournaments on $2$ vertices).

\begin{prop}\label{prop:rank1}
If $M$ is a finite monoid with group of units $G$, then
\[
(2|M|-|G|)\times\tbinom n2 \leq \rank(M\wr\Sing_n) \leq \begin{cases}|M|^2+1 &\text{if $n=2$}\\|M|^2\times\binom n2&\text{if $n\geq3$.}\end{cases}
\]
\end{prop}

\pf Put
\[
\Lam=\begin{cases}
\set{\ve_{12;ab}}{a,b\in M}\cup\{\ve_{21;11}\} &\text{if $n=2$}\\
\set{\ve_{ij;ab}}{(i,j)\in E_\Ga,\ a,b\in M} &\text{if $n\geq3$,}
\end{cases}
\]
where $\Ga$ is a strongly connected tournament with vertex set $\bn$ and edge set $E_\Ga$ in the $n\geq3$ case.  To establish the stated upper bound, it suffices to show that $M\wr\Sing_n=\la\Lam\ra$; by Theorem \ref{thm:IG}(i), it is enough to show that $\X_2\sub\la\Lam\ra$.  With this in mind, let $i,j\in\bn$ with $i\not=j$, and let $a,b\in M$.  We must show that $\ve_{ij;ab},\ve_{ji;ab}\in\la\Lam\ra$.  Without loss of generality, we may assume that $\ve_{ij;ab}\in\Lam$ (so $(i,j)=(1,2)$ if $n=2$, or $(i,j)\in E_\Ga$ if $n\geq3$).  Note also that $\ve_{ij;ba}\in\Lam$.  Next, note that $\Lam\cap\Sing_n$ is a generating set for $\Sing_n$ (regarded as a subsemigroup of $M\wr\Sing_n$, as usual), by Theorem \ref{thm:howie2}(ii).  In particular, $\ve_{ji;11}\in\la\Lam\ra$.  But then $\ve_{ji;ab}=\ve_{ij;ba}\ve_{ji;11}\in\la\Lam\ra$, as required.

To establish the stated lower bound, suppose $M\wr\Sing_n=\la\Lam\ra$.  Let $\oijn$, let $a\in M$, and consider an expression $\ve_{ij;a} = (\ba,\al)(\ba_1,\al_1)\cdots(\ba_k,\al_k)$, where $k\geq0$ and $(\ba,\al),(\ba_1,\al_1),\ldots,(\ba_k,\al_k)\in\Lam$.  Since $\ve_{ij}=\al\al_1\cdots\al_k$, we see that $\ker(\al)=\eta_{ij}$.  Write $(\bb,\be)=(\ba_1,\al_1)\cdots(\ba_k,\al_k)$, so that $\ve_{ij;a}=(\ba(\al\cdot\bb),\al\be)$.  In particular, examining the $i$th and $j$th coordinates of $\ba(\al\cdot\bb)$, we see that $1=a_ib_{i\al}$ and $a=a_jb_{j\al}=a_jb_{i\al}$.  It follows that $a_i\in G$ (as $M$ is finite) and $a_i^{-1}=b_{i\al}$, so that $a=a_ja_i^{-1}$.    To summarise, $\Lam$ contains an element $(\ba,\al)$ with $\ker(\al)=\eta_{ij}$, $a_i\in G$, and $a_ja_i^{-1}=a$.  A similar argument (considering a factorisation of~$\ve_{ji;a}$) shows that $\Lam$ contains an element $(\bc,\ga)$ with $\ker(\ga)=\eta_{ij}$, $c_j\in G$, and $c_ic_j^{-1}=a$.  If $a\in M\sm G$, then $(\ba,\al)\not=(\bc,\ga)$, or else then $a=a_ja_i^{-1}=c_ja_i^{-1}\in G$ (though might possibly be the case that $(\ba,\al)=(\bc,\ga)$ if $a\in G$).  In particular, $\Lam$ contains at least $2|M\sm G|+|G|=2|M|-|G|$ elements whose underlying transformation has kernel $\eta_{ij}$.  Since this is true for all $\oijn$, it follows that $|\Lam|\geq(2|M|-|G|)\times\binom n2$.  Since this is true for any generating set, the proof is complete. \epf

Although the next result is a special case of the one immediately following it, it will be convenient to state and prove it separately.

\begin{prop}\label{prop:rank2}
If $G$ is a finite group, then 
\[
\rank(G\wr\Sing_n)=\idrank(G\wr\Sing_n)=\begin{cases}
2 &\text{if $n=2$ and $|G|=1$}\\
|G|\times\binom n2 &\text{otherwise.}
\end{cases}
\]
\end{prop}

\pf If $n=2$ and $|G|=1$, then $G\wr\Sing_n\cong\Sing_2$ is a right zero semigroup of size~$2$, so the result is clear.  For the rest of the proof, suppose at least one of $n\geq3$ or $|G|\geq2$ holds.  By Proposition \ref{prop:rank1}, the proof will be complete if we can produce an idempotent generating set for $G\wr\Sing_n$ of the specified size.  With this in mind, put
\[
\Lam=\begin{cases}
\set{\ve_{12;a}}{a\in G\sm\{1\}}\cup\{\ve_{21;1}\} &\text{if $n=2$}\\
\set{\ve_{ij;a}}{(i,j)\in E_\Ga,\ a\in G} &\text{if $n\geq3$,}
\end{cases}
\]
where $\Ga$ is a strongly connected tournament with vertex set $\bn$ and edge set $E_\Ga$ in the $n\geq3$ case.  Since~$\Lam$ has the required size, it suffices to show that $\X_1\sub\la\Lam\ra$, by Theorem \ref{thm:IG}(ii) and (iii).  As in the proof of Proposition \ref{prop:rank1}, we immediately see that $\la\Lam\ra$ contains $\X$ if $n\geq3$.  We claim that this is also the case for $n=2$.  Indeed, if $n=2$, then we already have $\ve_{21;1}\in\Lam$, while for any $a\in G\sm\{1\}$ with $a^k=1$ (and $k\geq1$), we have $\ve_{12;1}=(\ve_{12;a,a})^k=(\ve_{21;1}\ve_{12;a})^k\in\la\Lam\ra$.  So far we have shown that $\X\sub\la\Lam\ra$.  Let $i,j\in\bn$ with $i\not=j$, and let $a\in G\sm\{1\}$.  The proof will be complete if we can show that $\ve_{ij;a},\ve_{ji;a}\in\la\Lam\ra$.  Without loss of generality, we may assume that $\ve_{ij;a}\in\Lam$, which also implies that $\ve_{ij;a^{-1}}\in\Lam$ (by definition of $\Lam$).  But then $\ve_{ji;a}=\ve_{ij;a^{-1}}\ve_{ji;1}\ve_{ij;a}\ve_{ji;1}\in\la\Lam\ra$. \epf

The next result extends Proposition \ref{prop:rank2} to the case of an arbitrary finite idempotent-generated singular wreath product.
In parts of its proof, we use some general results of Gray \cite{Gray2008} on (idempotent) ranks of completely $0$-simple semigroups; see also \cite{EastGray,Gray2014}.
We also make use of Green's relations, as defined in Section~\ref{sect:prelim}.

\begin{thm}\label{thm:rank2}
If $M$ is a finite monoid with group of units $G$, and if $M/{\L}$ is a chain, then
\begin{align*}
\rank(M\wr\Sing_n) &= \begin{cases}
2 &\text{if $n=2$ and $|M|=1$}\\
(2|M|-|G|)\times\binom n2 &\text{otherwise,}
\end{cases}
\intertext{and}
\idrank(M\wr\Sing_n) &= \begin{cases}
2|M| &\text{if $n=2$ and $|G|=1$}\\
(2|M|-|G|)\times\binom n2 &\text{otherwise.}
\end{cases}
\end{align*}
\end{thm}

\pf The case in which $M=G$ is covered in Proposition \ref{prop:rank2}, so suppose $M\not=G$.  In particular, $|M|\geq2$.  We break the proof up into three cases.

\pfcase{1}
Suppose first that $n\geq3$.  Because of Proposition \ref{prop:rank1}, the proof will be complete (in this case) if we can produce an idempotent generating set for $M\wr\Sing_n$ of size $(2|M|-|G|)\times\binom n2$.  With this in mind, let~$\Ga$ be a strongly connected tournament on vertex set $\bn$, and put $\Lam=\Lam_1\cup\Lam_2$, where
\[
\Lam_1 = \set{\ve_{ij;a}}{i,j\in\bn,\ i\not=j,\ a\in M\sm G} \AND \Lam_2 = \set{\ve_{ij;a}}{(i,j)\in E_\Ga,\ a\in G}.
\]
The proof of Proposition \ref{prop:rank2} gives $\la\Lam_2\ra=G\wr\Sing_n$.  Consequently, $\la\Lam\ra$ contains $\X_1$, and we are done, by Theorem~\ref{thm:IG}.

\pfcase{2}
Next suppose $n=2$ and $|G|\geq2$.  This time, put $\Lam=\Lam_1\cup\Lam_2$, where
\[
\Lam_1 = \set{\ve_{12;a},\ve_{21;a}}{a\in M\sm G} \AND \Lam_2 = \set{\ve_{12;a}}{a\in G\sm\{1\}}\cup\{\ve_{21;1}\}.
\]
Again, the proof of Proposition \ref{prop:rank2} gives $\la\Lam_2\ra=G\wr\Sing_2$, and it quickly follows that $\Lam$ is an idempotent generating set of the required size, completing the proof in this case.

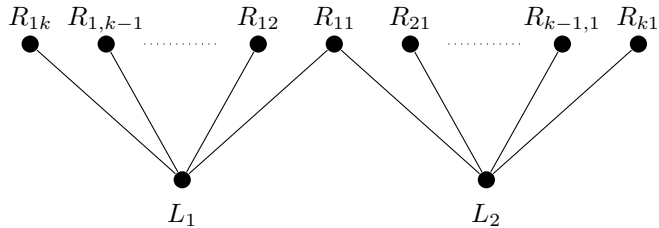
\begin{figure}[h]
\begin{center}
\begin{tikzpicture}[yscale=0.9]
\tikzstyle{vertex}=[circle,fill=black, inner sep = 0.08cm]
\foreach \x in {1,2,4,5,6,8,9} {\node[vertex] (U\x) at (\x,2) {};}
\foreach \x in {3,7} {\node[vertex] (L\x) at (\x,0) {};}
\foreach \x in {1,2,4,5} {\draw (U\x)--(L3);}
\foreach \x in {5,6,8,9} {\draw (U\x)--(L7);}
\foreach \x/\y in {1/k,{2/,k-1},4/2,5/1} {\node[below] () at (\x,2.7) {\small $R_{1\y}$};}
\foreach \x/\y in {9/k,{8/k-1,},6/2} {\node[below] () at (\x,2.7) {\small $R_{\y1}$};}
\foreach \x/\y in {3/1,7/2} {\node[below] () at (\x,-.2) {\small $L_{\y}$};}
\draw[dotted] (2.5,2)--(3.5,2) (6.5,2)--(7.5,2);
\end{tikzpicture}
\end{center}
\vspace{-5mm}
\caption{The graph $\De_{E(J)}$; see the proof of Theorem \ref{thm:rank2} for more details.}
\label{fig:De}
\end{figure}

\pfcase{3}
Finally, suppose $n=2$ and $|G|=1$, noting that $M\wr\Sing_2=\X_2=\set{\ve_{12;ab},\ve_{21;ab}}{a,b\in M}$.  It is easy to check that for any $a\in M$, and for $(i,j)=(1,2)$ or~$(2,1)$,
\[
R_{\ve_{ij;1a}}=R_{\ve_{ji;a1}} = \{\ve_{ij;1a},\ve_{ji;a1}\} \AND
L_{\ve_{ij;1a}} = L_{\ve_{ij;a1}} = \set{\ve_{ij;1b},\ve_{ij;b1}}{b\in M}.
\]
It follows that the set 
\[
J = \set{\ve_{12;ab},\ve_{21;ab}}{a,b\in M,\ 1\in\{a,b\}}
\]
is a ${\D}={\J}$-class of $M\wr\Sing_2$.  It is clear that $\X_1=E(J)$, so it follows that $M\wr\Sing_2=\la J\ra=\la E(J)\ra$.  It then follows from \cite[Lemma 3.2]{EastGray} that 
\[
\rank(M\wr\Sing_2)=\rank(J^\circ) \AND \idrank(M\wr\Sing_2)=\idrank(J^\circ),
\]
where $J^\circ$ is the \emph{principal factor} of $J$: i.e., the semigroup with underlying set $J\cup\{0\}$ and operation $\circ$ defined, for $x,y\in J\cup\{0\}$, by
\[
x\circ y = \begin{cases}
xy &\text{if $x,y,xy\in J$}\\
0 &\text{otherwise.}
\end{cases}
\]
Since $M\wr\Sing_2$ is generated by $\X_1=E(J)$, it follows that $J^\circ$ is idempotent-generated.  It then follows from \cite[Lemma 2.3]{Gray2008} that $\rank(J^\circ)$ is equal to the larger of $|J/{\R}|$ and $|J/{\L}|$.  Write $M=\{a_1,\ldots,a_k\}$, where $k=|M|$ and $a_1=1$.  From the above calculations, we see that the $\L$- and $\R$-classes contained in $J$ are
\[
L_1=L_{\ve_{12;11}}\COMMA L_2=L_{\ve_{21;11}} \AND
R_{1s}=R_{\ve_{12;1a_s}} \COMMA R_{s1}=R_{\ve_{12;a_s1}} \qquad\text{for each $1\leq s\leq k$,}
\]
so that $|J/{\R}|=2k-1$ and $|J/{\L}|=2$, giving $\rank(J^\circ)=2k-1$.
To complete the proof, we need to show that $\idrank(J^\circ)=2k$.  To do this, we use some more ideas from \cite{Gray2008}.  For a subset $F\sub \X_1=E(J)$, define a (bipartite) graph $\De_F$ as follows: the vertex set of $\De_F$ is $(J/{\R})\cup(J/{\L})$, and there is an (undirected) edge between $R\in J/{\R}$ and $L\in J/{\L}$ if and only if $(R\cap L)\cap F\not=\emptyset$.  Such a graph $\De_F$ is a subgraph of $\De_{E(J)}$, which is pictured in Figure \ref{fig:De}.
By \cite[Lemma 2.5]{Gray2008}, $J^\circ=\la F\ra$ if and only if $\De_F$ is connected.
Since the removal of any edge from $\De_{E(J)}$ disconnects the graph, it follows that no proper subset $F$ of $E(J)$ generates~$J^\circ$.
So $\idrank(J^\circ)=|E(J)|=2k$, and the proof is complete. \epf

\begin{rem}
Theorem \ref{thm:rank2} shows that when $M\wr\Sing_n$ is idempotent-generated, its rank and idempotent rank are equal unless $n=2$ and $|G|=1\not=|M|$, in which case the rank and idempotent rank are equal to~$2|M|-1$ and $2|M|$, respectively.
It would be interesting to calculate exact values for $\rank(M\wr\Sing_n)$ for various classes of finite monoids $M$ where $M/{\L}$ is not a chain.  
\end{rem}

We now give some information about the (idempotent) rank of $M\wr\Sing_n$ for infinite $M$.
A monoid $M$ has a natural (right) \emph{diagonal} action on the set $M\times M$, defined by $(a,b)\cdot c=(ac,bc)$, for all $a,b,c\in M$.  We say this action is \emph{finitely generated} if there exists a finite subset $\Omega\sub M\times M$ such that $\Omega\cdot M=M\times M$.  It is easy to check that the diagonal action of an infinite group is never finitely generated.

\begin{thm}\label{thm:infinite}
\begin{itemize}
\itemit{i} If $M$ is an infinite monoid, then $M\wr\Sing_n$ is finitely generated if and only if $M$ is finitely generated (as a semigroup) and the diagonal action of $M$ on $M\times M$ is finitely generated.
\itemit{ii} If $M\wr\Sing_n$ is not finitely generated (e.g., if $M$ is an infinite group), then $\rank(M\wr\Sing_n)=|M|$; if $M/{\L}$ is also a chain, then $\idrank(M\wr\Sing_n)=|M|$.
\end{itemize}
\end{thm}

\pf (i).  Suppose first that $M\wr\Sing_n$ is generated by a finite set $\Lam=\set{(\ba_s,\al_s)}{1\leq s\leq k}$.  For $1\leq s\leq k$, write $\ba_s=(a_{s1},\ldots,a_{sn})$, and put
\[
A=\set{a_{st}}{1\leq s\leq k,\ 1\leq t\leq n} \AND \Omega = \set{(a_{s1},a_{s2})}{1\leq s\leq k},
\]
noting that 
both sets are finite.
Let $a\in M$ be arbitrary.  Consideration of $\ve_{12;a}$ as a product of elements from $\Lam$ shows that $a\in\la A\ra$, giving $M=\la A\ra$.  
To show that $M\times M$ is finitely generated as an $M$-act, let $a,b\in M$ be arbitrary, and consider a product 
\[
\ve_{12;ab} = (\ba_s,\al_s)(\ba_{s_1},\al_{s_1})\cdots(\ba_{s_r},\al_{s_r}) 
\qquad\text{where $1\leq s,s_1,\ldots,s_r\leq k$.}
\]
From $\ve_{12}=\al_s\al_{s_1}\cdots\al_{s_r}$, we deduce $\ker(\al_s)=\eta_{12}$.
Writing $(\bb,\be)=(\ba_{s_1},\al_{s_1})\cdots(\ba_{s_r},\al_{s_r})$, and examining the first and second coordinates of $\ve_{12;ab}=(\ba_s,\al_s)(\bb,\be)=(\ba_s(\al_s\cdot\bb),\al_s\be)$, we see that 
\[
(a,b)=(a_{s1}b_{1\al_s},a_{s2}b_{2\al_s}) =(a_{s1}b_{1\al_s},a_{s2}b_{1\al_s}) = (a_{s1},a_{s2})\cdot b_{1\al_s},
\]
showing that $M\times M=\Omega\cdot M$.

Conversely, suppose $M=\la A\ra$ and $M\times M=\Omega\cdot M$ for finite sets $A\sub M$ and $\Omega\sub M\times M$.  Put
\[
\Lam = \set{\ve_{ij;ab}}{i,j\in\bn,\ i\not=j,\ (a,b)\in\Omega} \cup \set{\ve_{ij;a1}}{i,j\in\bn,\ i\not=j,\ a\in A}.
\]
Since $\Lam$ is finite, the proof of (i) will be complete (by Theorem \ref{thm:IG}(i)) if we can show that $\la\Lam\ra$ contains~$\X_2$.  So let $a,b\in M$, and let $i,j\in\bn$ with $i\not=j$.  Now, $(a,b)=(c,d)\cdot f=(cf,df)$ for some $(c,d)\in\Omega$ and $f\in M$, and we may write $f=f_1\cdots f_r$ for some $f_1,\ldots,f_r\in A$.  
But then
\[
\ve_{ij;ab} = \ve_{ij;cd}\ve_{ij;f1} = \ve_{ij;cd}\ve_{ij;f_11}\cdots\ve_{ij;f_r1}\in\la\Lam\ra,
\]
as required.  

\pfitem{ii}  This follows from the fact that $\rank(S)=|S|$ for any non-finitely generated semigroup $S$ (and the corresponding statement for idempotent ranks of idempotent-generated semigroups).~\epf

\begin{rem}
We have not attempted to give the value (or any estimate) for the rank (or idempotent rank, if applicable) of $M\wr\Sing_n$ in the case that $M\wr\Sing_n$ is infinite but finitely generated.  Although the generating set $\Lam$ constructed in the second paragraph of the proof of Theorem \ref{thm:infinite}(i) leads to an upper bound for $\rank(M\wr\Sing_n)$, $\Lam$ is certainly not of minimal size, in general.
%
It would be interesting to calculate exact values for the rank (and idempotent rank, if appropriate) of $M\wr\Sing_n$ for various classes of finitely generated infinite monoids $M$ with finitely generated diagonal action.
\end{rem}

\begin{rem}
Finite generation of the diagonal action plays an important role in \cite{RRT2002}.  However, the wreath products studied there are different to those studied here, so Theorem \ref{thm:infinite}(i) appears to be independent of the results from \cite{RRT2002}.  See also \cite{RRT2001,Gallagher2006,GR2005}, where diagonal actions are studied in detail, and several interesting examples discussed.
\end{rem}

\boldsection{Presentations for $M\wr\Sing_n$}\label{sect:pres}

We now turn to the main topic of the paper: namely, the task of finding presentations for the singular wreath products $M\wr\Sing_n$.  In Section~\ref{subsect:X2}, we give a presentation in terms of the generating set $\X_2$ (see Theorem \ref{thm:X2}).  In Section \ref{subsect:X1}, we give a presentation in terms of the idempotent generating set~$\X_1$ in the case that $M/{\L}$ is a chain (see Theorem~\ref{thm:X1}).  Our strategy is to begin with the general presentation for semidirect products (Theorem \ref{thm:semidirect}) in order to deduce a presentation for $M\wr\Sing_n$ in terms of the (very large) generating set $\X_n$ (Corollary \ref{cor:Xn}), and to then reduce this to the above-mentioned simpler presentations, using Tietze transformations.  Each of the presentations we give in this section extends the presentation for $\Sing_n$ stated in Theorem \ref{thm:X}.

\boldsubsection{A presentation for $M\wr\Sing_n$}\label{subsect:X2}

Since $M\wr\Sing_n=M^n\rtimes\Sing_n$ is a semidirect product, Theorem \ref{thm:semidirect} allows us to write down a presentation for $M\wr\Sing_n$ in terms of the presentation $\pres{X\!}{\!R}$ for $\Sing_n$ from Theorem \ref{thm:X}.  In order to state this presentation (in Corollary \ref{cor:Xn} below), define an alphabet
\[
X_n = \set{e_{ij;\ba}}{i,j\in\bn,\ i\not=j,\ \ba\in M^n},
\]
an epimorphism
\[
\phi_n:X_n^+\to M\wr\Sing_n:e_{ij;\ba}\mt\ve_{ij;\ba},
\]
and let $R_n$ be the set of relations (identifying a letter $e_{ij}\in X$ with $e_{ij;(1,\ldots,1)}\in X_n$, as in Section \ref{sect:semidirect})
\begin{align}
\label{R1n}\tag*{(R1)$_n$} e_{ij;\ba}e_{ij} = e_{ij;\ba} &= e_{ji;\ba}e_{ij} &&\text{for $\ba\in M^n$ and distinct $i,j$}\\
\label{R2n}\tag*{(R2)$_n$} e_{ij;\ba}e_{kl} &= e_{kl;\ba}e_{ij} &&\text{for $\ba\in M^n$ and distinct $i,j,k,l$}\\
\label{R3n}\tag*{(R3)$_n$} e_{ik;\ba}e_{jk} &= e_{ik;\ba} &&\text{for $\ba\in M^n$ and  distinct $i,j,k$}\\
\label{R4n}\tag*{(R4)$_n$} e_{ij;\ba}e_{ik} = e_{ik;\ba}e_{ij} &= e_{jk;\ba}e_{ij} &&\text{for $\ba\in M^n$ and  distinct $i,j,k$}\\
\label{R5n}\tag*{(R5)$_n$} e_{ki;\ba}e_{ij}e_{jk} &= e_{ik;\ba}e_{kj}e_{ji}e_{ik} &&\text{for $\ba\in M^n$ and  distinct $i,j,k$}\\
\label{R6n}\tag*{(R6)$_n$} e_{ki;\ba}e_{ij}e_{jk}e_{kl} &= e_{ik;\ba}e_{kl}e_{li}e_{ij}e_{jl} &&\text{for $\ba\in M^n$ and  distinct $i,j,k,l$}\\
\label{R7n}\tag*{(R7)$_n$} e_{ij;\ba}e_{kl;\bb} &= e_{ij;\bc}e_{kl} &&\text{for $\ba,\bb\in M^n$ and any $i,j,k,l$,}
\end{align}
where in \ref{R7n}, $\bc=\ba(\ve_{ij}\cdot\bb)=(c_1,\ldots,c_n)$ satisfies $c_j=a_jb_i$ and $c_k=a_kb_k$ for $k\not=j$.  Note that $i,j,k,l$ are not assumed to be distinct (apart from $i\not=j$ and $k\not=l$) in \ref{R7n}.  As noted above, the next result is a special case of Theorem \ref{thm:semidirect}.

\begin{cor}\label{cor:Xn}
The semigroup $M\wr\Sing_n$ has presentation $\pres{X_n}{R_n}$ via $\phi_n$. \epfres
\end{cor}

The presentation $\pres{X_n}{R_n}$ utilises the very large generating set $\X_n$.  Consequently, our next goal is to simplify this presentation to obtain a presentation utilising the considerably smaller generating set $\X_2\sub\X_n$.  With this in mind, define an alphabet
\[
X_2 = \set{e_{ij;ab}}{i,j\in\bn,\ i\not=j,\ a,b\in M},
\]
an epimorphism
\[
\phi_2:X_2^+\to M\wr\Sing_n:e_{ij;ab}\mt\ve_{ij;ab},
\]
and let $R_2$ be the set of relations
\begin{align}
\label{R12}\tag*{(R1)$_2$} e_{ij;ab}e_{ij;cd} = e_{ij;ac,bc} &= e_{ji;ba}e_{ij;dc} &&\text{for $a,b,c,d\in M$ and distinct $i,j$}\\
\label{R22}\tag*{(R2)$_2$} e_{ij;ab}e_{kl;cd} &= e_{kl;cd}e_{ij;ab} &&\text{for $a,b,c,d\in M$ and distinct $i,j,k,l$}\\
\label{R3a2}\tag*{(R3a)$_2$} e_{ik;ab}e_{jk;1c} &= e_{ik;ab} &&\text{for $a,b,c\in M$ and distinct $i,j,k$}\\
\label{R3b2}\tag*{(R3b)$_2$} e_{ik;ab}e_{jk;c1} &= e_{ki;ba}e_{ji;c1}e_{ik;11} &&\text{for $a,b,c\in M$ and distinct $i,j,k$}\\
\label{R3c2}\tag*{(R3c)$_2$} e_{ik;aa}e_{jk;b1} &= e_{ik;11}e_{jk;b1}e_{ik;a1} &&\text{for $a,b\in M$ and distinct $i,j,k$}\\
\label{R4a2}\tag*{(R4a)$_2$} e_{ij;ab}e_{ik;cd} = e_{ik;ac,d}e_{ij;1,bc} &= e_{jk;bc,d}e_{ij;ac,1} &&\text{for $a,b,c,d\in M$ and distinct $i,j,k$}\\
\label{R4b2}\tag*{(R4b)$_2$} e_{ij;c,ad}e_{ik;1,bd} = e_{ik;c,bd}e_{ij;1,ad} &= e_{jk;ab}e_{ij;cd} &&\text{for $a,b,c,d\in M$ and distinct $i,j,k$}\\
\label{R52}\tag*{(R5)$_2$} e_{ki}e_{ij}e_{jk} &= e_{ik}e_{kj}e_{ji}e_{ik} &&\text{for distinct $i,j,k$}\\
\label{R62}\tag*{(R6)$_2$} e_{ki}e_{ij}e_{jk}e_{kl} &= e_{ik}e_{kl}e_{li}e_{ij}e_{jl} &&\text{for distinct $i,j,k,l$.}
\end{align}
Note that the labelling of relations is chosen to reflect the labels of the relations from $R$ (as stated in Section~\ref{sect:prelim}), some of which have been split up into separate relations in $R_2$.
Note also that in relations \ref{R52} and \ref{R62}, we identify $X$ as a subset of $X_2$ via $e_{ij}\equiv e_{ij;11}$; in this way, relations \ref{R52} and \ref{R62} are really just \eqref{R5} and \eqref{R6}.
Finally, note also that in several relations, we have separated the two monoid subscripts on a letter from $X_2$ with a comma when at least one of the subscripts involves a product; for example, the two monoid subscripts on ``$e_{ij;ac,bc}$'' in relation \ref{R12} are the products $ac$ and $bc$ from $M$.
Our goal in this section is to prove the following.

\begin{thm}\label{thm:X2}
The semigroup $M\wr\Sing_n$ has presentation $\pres{X_2}{R_2}$ via $\phi_2$.  
\end{thm}

\begin{rem}\label{rem:R4}
The ``$e_{ik;ac,d}e_{ij;1,bc} = e_{jk;bc,d}e_{ij;ac,1}$'' part of \ref{R4a2} follows from the ``$e_{ik;c,bd}e_{ij;1,ad} = e_{jk;ab}e_{ij;cd}$'' part of \ref{R4b2}, upon making the substitution $(a,b,c,d)\to(bc,d,ac,1)$.  Similarly, the ``$e_{ij;c,ad}e_{ik;1,bd} = e_{ik;c,bd}e_{ij;1,ad}$'' part of \ref{R4b2} follows from the ``$e_{ij;ab}e_{ik;cd} = e_{ik;ac,d}e_{ij;1,bc}$'' part of \ref{R4a2}.  As such, we could replace \ref{R4a2} and \ref{R4b2} by just
\[
e_{ij;ab}e_{ik;cd} = e_{ik;ac,d}e_{ij;1,bc} \AND e_{ik;c,bd}e_{ij;1,ad} = e_{jk;ab}e_{ij;cd} .
\]
However, in the calculations that follow, it is convenient to leave both relations as they are.
\end{rem}

As noted earlier, we will prove Theorem \ref{thm:X2} by performing a sequence of \emph{Tietze transformations}, beginning with the presentation $\pres{X_n}{R_n}$ from Corollary \ref{cor:Xn}.  We will write ${\sim_n}=R_n^\sharp=\ker(\phi_n)$ for the congruence on $X_n^+$ generated by~$R_n$, and we think of $X_2$ as a subset of $X_n$ by identifying $e_{ij;ab}\in X_2$ with $e_{ij;\ba}\in X_n$, where $\ba=(a_1,\ldots,a_n)$ satisfies $a_i=a$, $a_j=b$ and $a_k=1$ if $k\not\in\{i,j\}$.

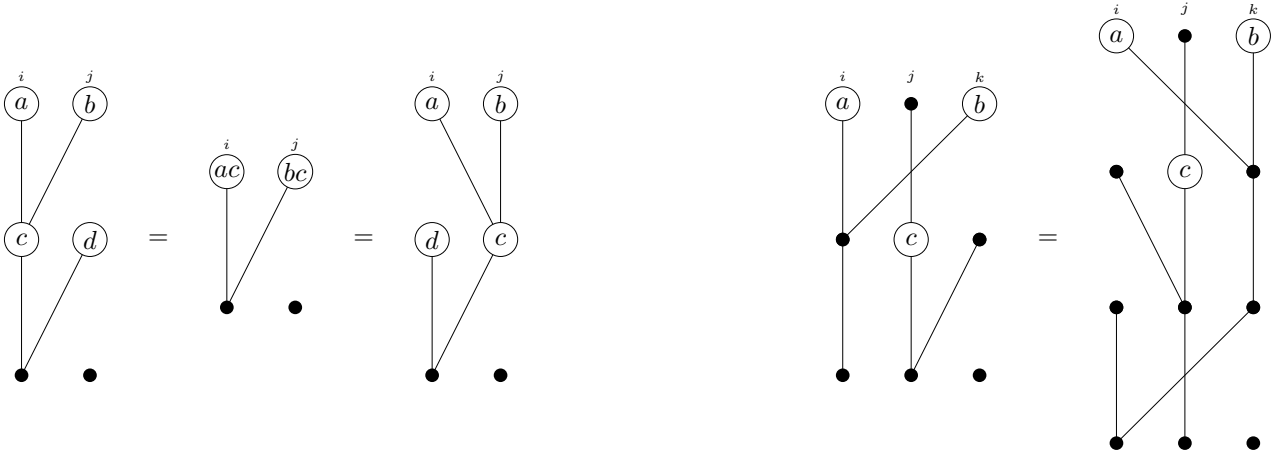
\begin{figure}[h]
\begin{center}
\scalebox{.9}{
\begin{tikzpicture}[scale=1]
\tikzstyle{vertex}=[circle,draw=black, fill=white,minimum size=5mm, inner sep = 0.03cm]
\begin{scope}[shift={(0,0)}]	
\uvs{1,2}
\lvs{1,2}
\stlines{1/1,2/1}
\vertlabelsh{1/i,2/j}
\node[vertex] () at (1,2){$a$};
\node[vertex] () at (2,2){$b$};
\end{scope}
\begin{scope}[shift={(0,-2)}]	
\uvs{1,2}
\lvs{1,2}
\stlines{1/1,2/1}
\node[vertex] () at (1,2){$c$};
\node[vertex] () at (2,2){$d$};
\end{scope}
\begin{scope}[shift={(3,-1)}]	
\vertlabelsh{1/i,2/j}
\uvs{1,2}
\lvs{1,2}
\stlines{1/1,2/1}
\node[vertex] () at (1,2){$ac$};
\node[vertex] () at (2,2){$bc$};
\draw(0,1)node{$=$};
\draw(3,1)node{$=$};
\end{scope}
\begin{scope}[shift={(6,0)}]	
\uvs{1,2}
\lvs{1,2}
\stlines{1/2,2/2}
\vertlabelsh{1/i,2/j}
\node[vertex] () at (1,2){$a$};
\node[vertex] () at (2,2){$b$};
\end{scope}
\begin{scope}[shift={(6,-2)}]	
\uvs{1,2}
\lvs{1,2}
\stlines{1/1,2/1}
\node[vertex] () at (1,2){$d$};
\node[vertex] () at (2,2){$c$};
\end{scope}
%
\begin{scope}[shift={(12+0,0)}]	
\uvs{1,2,3}
\lvs{1,2,3}
\stlines{1/1,2/2,3/1}
\vertlabelsh{1/i,2/j,3/k}
\node[vertex] () at (1,2){$a$};
\node[vertex] () at (3,2){$b$};
\end{scope}
\begin{scope}[shift={(12+0,-2)}]	
\uvs{1,2,3}
\lvs{1,2,3}
\stlines{1/1,2/2,3/2}
\node[vertex] () at (2,2){$c$};
\end{scope}
\begin{scope}[shift={(12+4,1)}]	
\uvs{1,2,3}
\lvs{1,2,3}
\stlines{1/3,2/2,3/3}
\vertlabelsh{1/i,2/j,3/k}
\node[vertex] () at (1,2){$a$};
\node[vertex] () at (3,2){$b$};
\end{scope}
\begin{scope}[shift={(12+4,-1)}]	
\uvs{1,2,3}
\lvs{1,2,3}
\stlines{1/2,2/2,3/3}
\node[vertex] () at (2,2){$c$};
\draw(0,1)node{$=$};
\end{scope}
\begin{scope}[shift={(12+4,-3)}]	
\uvs{1,2,3}
\lvs{1,2,3}
\stlines{1/1,2/2,3/1}
\end{scope}
\end{tikzpicture}
}
\end{center}
\vspace{-5mm}
\caption{Diagrammatic verification of relations \ref{R12}, left, and \ref{R3b2}, right.  Only the relevant parts of the diagrams have been included, and no ordering on $i,j$ (or $i,j,k$) is to be implied.}
\label{fig:R2rels}
\end{figure}

One may check (diagrammatically) that $u\phi_n=v\phi_n$ (in $M\wr\Sing_n$) for each relation $(u,v)$ from $R_2$.  We do this for relations \ref{R12} and \ref{R3b2} in Figure \ref{fig:R2rels}, and leave the reader to check the rest.  In particular, we may add the relations $R_2$ to the presentation $\pres{X_n}{R_n}$ to obtain $\pres{X_n}{R_2\cup R_n}$.  

Out next goal is to show that the generators from $X_n\sm X_2$ may be removed from the presentation.  To do this, we define words $E_{ij;\ba}\in X_2^+$, for each 
$i,j\in\bn$ with $i\not=j$ and each $\ba=(a_1,\ldots,a_n)\in M^n$, as follows.  If it happens that $a_k=1$ for all $k\in\bn\sm\{i,j\}$, then we simply define $E_{ij;\ba}=e_{ij;a_ia_j}$.
Otherwise, we  first write $\bn\sm\{i,j\}=\{i_1,\ldots,i_{n-2}\}$, where $i_1<\cdots<i_{n-2}$, and define
\[
E_{ij;\ba} = e_{ij;a_ia_j} e_{i_1j;a_{i_1}1}\cdots e_{i_{n-2}j;a_{i_{n-2}}1}.
\]
As shown in Figure \ref{fig:eijba}, we have $E_{ij;\ba}\phi_n=\ve_{ij;\ba}=e_{ij;\ba}\phi_n$.  In particular, because ${\sim_n}=\ker(\phi_n)$, it follows that $e_{ij;\ba}\sim_n E_{ij;\ba}$.  So we may remove each generator $e_{ij;\ba}\in X_n\sm X_2$ from the presentation $\pres{X_n}{R_2\cup R_n}$, replacing every occurrence of such a generator in the relations by the word $E_{ij;\ba}\in X_2^+$.  Since $R_2\sub X_2^+\times X_2^+$, the only relations modified in this way are those from $R_n$.  
We label the relations modified in this way as
\begin{align}
\label{R1n'}\tag*{(R1)$'_n$} E_{ij;\ba}e_{ij} = E_{ij;\ba} &= E_{ji;\ba}e_{ij} &&\text{for $\ba\in M^n$ and distinct $i,j$}\\
\label{R2n'}\tag*{(R2)$'_n$} E_{ij;\ba}e_{kl} &= E_{kl;\ba}e_{ij} &&\text{for $\ba\in M^n$ and distinct $i,j,k,l$}\\
\label{R3n'}\tag*{(R3)$'_n$} E_{ik;\ba}e_{jk} &= E_{ik;\ba} &&\text{for $\ba\in M^n$ and  distinct $i,j,k$}\\
\label{R4n'}\tag*{(R4)$'_n$} E_{ij;\ba}e_{ik} = E_{ik;\ba}e_{ij} &= E_{jk;\ba}e_{ij} &&\text{for $\ba\in M^n$ and  distinct $i,j,k$}\\
\label{R5n'}\tag*{(R5)$'_n$} E_{ki;\ba}e_{ij}e_{jk} &= E_{ik;\ba}e_{kj}e_{ji}e_{ik} &&\text{for $\ba\in M^n$ and  distinct $i,j,k$}\\
\label{R6n'}\tag*{(R6)$'_n$} E_{ki;\ba}e_{ij}e_{jk}e_{kl} &= E_{ik;\ba}e_{kl}e_{li}e_{ij}e_{jl} &&\text{for $\ba\in M^n$ and  distinct $i,j,k,l$}\\
\label{R7n'}\tag*{(R7)$'_n$} E_{ij;\ba}E_{kl;\bb} &= E_{ij;\bc}e_{kl} &&\text{for $\ba,\bb\in M^n$ and any $i,j,k,l$,}
\end{align}
and denote the whole set of relations modified in this way by $R_n'$.
So the presentation has now become $\pres{X_2}{R_2\cup R_n'}$.  To complete the proof of Theorem \ref{thm:X2}, it remains to show that the relations $R_n'$ can be removed.  That is, we need to show that $u\sim_2 v$ for each relation $(u,v)\in R_n'$, where ${\sim_2}=R_2^\sharp$ denotes the congruence on $X_2^+$ generated by $R_2$.
We begin this task with some technical results that will be of use on a number of occasions.

\begin{lemma}\label{lem:order}
For any $i,j,k,l\in\bn$ with $i\not=j$ and with $j,k,l$ distinct, and for any $a,b,c,d\in M$, we have $e_{ij;ab}e_{kj;c1}e_{lj;d1} \sim_2 e_{ij;ab}e_{lj;d1}e_{kj;c1}$.  
\end{lemma}

\pf We have
\begin{align*}
e_{ij;ab}e_{kj;c1}e_{lj;d1} &\sim_2 e_{ij;ab}e_{kj;11}e_{kj;c1}e_{lj;d1} &&\text{by \ref{R12} if $k=i$, or \ref{R3a2} if $k\not=i$}\\
&\sim_2 e_{ij;ab}e_{kj;cc}e_{lj;d1}  &&\text{by \ref{R12}}\\
&\sim_2 e_{ij;ab}e_{kj;11}e_{lj;d1}e_{kj;c1}  &&\text{by \ref{R3c2}}\\
&\sim_2 e_{ij;ab}e_{lj;d1}e_{kj;c1}  &&\text{by \ref{R12} or \ref{R3a2}.}  \qedhere
\end{align*}
\epf

\begin{cor}\label{cor:order}
For any $i,j,k_1,\ldots,k_t\in\bn$ with $i\not=j$ and $j,k_1,\ldots,k_t$ distinct, for any ${a,b,c_1,\ldots,c_t\in M}$, and for any permutation $\pi\in\S_t$, we have
$
e_{ij;ab}e_{k_1j;c_11}\cdots e_{k_tj;c_t1} \sim_2 e_{ij;ab}e_{k_{1\pi}j;c_{1\pi}1}\cdots e_{k_{t\pi}j;c_{t\pi}1}.
$
\end{cor}

\pf Consider the subword $e_{k_sj;c_s1}e_{k_{s+1}j;c_{s+1}1}$ in the product $e_{k_1j;c_11}\cdots e_{k_tj;c_t1}$.  The letter immediately before this subword is $e_{lj;uv}$ for some $l\in\bn\sm\{j\}$ and some $u,v\in M$.  In particular, Lemma \ref{lem:order} gives $e_{lj;uv}e_{k_sj;c_s1}e_{k_{s+1}j;c_{s+1}1} \sim_2 e_{lj;uv}e_{k_{s+1}j;c_{s+1}1}e_{k_sj;c_s1}$.  So the result is true if $\pi$ is a simple transposition $(s,s+1)$.  Since $\pi$ is the product of such simple transpositions, the result follows. \epf

As a special case of the previous result, we may immediately deduce the following.

\begin{cor}\label{cor:Eij}
Let $i,j\in\bn$ with $i\not=j$, and let $\ba\in M^n$.  If $j_1,\ldots,j_{n-2}$ is any ordering of $\bn\sm\{i,j\}$, then
\[
\tag*{\qed}
E_{ij;\ba} \sim_2 e_{ij;a_ia_j} e_{j_1j;a_{j_1}1}\cdots a_{j_{n-2}j;a_{j_{n-2}}1}.  
\]
\end{cor}

\begin{lemma}\label{lem:chain}
Suppose $i,j,k_1,\ldots,k_t\in\bn$ are distinct, and let $a,b,c_1,\ldots,c_t\in M$, where $t\geq0$.  Then
\[
e_{ij;ab}  e_{k_1j;c_11}\cdots e_{k_tj;c_t1} 
\sim_2 
e_{ji;ba}  e_{k_1i;c_11}\cdots e_{k_ti;c_t1} e_{ij;11}.
\]
\end{lemma}

\pf We use induction on $t$.  If $t=0$, then the words $e_{k_1j;c_11}\cdots e_{k_tj;c_t1}$ and $e_{k_1i;c_11}\cdots e_{k_ti;c_t1}$ are empty, and we have $e_{ij;ab} \sim e_{ji;ba}e_{ij;11}$, by \ref{R12}.  If $t\geq1$, then
\begin{align*}
e_{ij;ab}  e_{k_1j;c_11}&\cdots e_{k_{t-1}j;c_{t-1}1}e_{k_tj;c_t1} 
 \sim_2
e_{ji;ba}  e_{k_1i;c_11}\cdots e_{k_{t-1}i;c_{t-1}1} e_{ij;11} e_{k_tj;c_t1}  &&\text{by an induction hypothesis}\\
& \sim_2
e_{ji;ba}  e_{k_1i;c_11}\cdots e_{k_{t-1}i;c_{t-1}1} e_{ji;11} e_{k_ti;c_t1}e_{ij;11}  &&\text{by \ref{R3b2}}\\
& \sim_2
e_{ji;ba}  e_{k_1i;c_11}\cdots e_{k_{t-1}i;c_{t-1}1} e_{k_ti;c_t1}e_{ij;11}  &&\hspace{-1.8cm}\text{by \ref{R3a2} if $t\geq2$, or \ref{R12} if $t=1$.} \qedhere
\end{align*}
\epf

\begin{lemma}\label{lem:d1}
Let $i,j,k\in\bn$ with $i\not=j$ and $k\not=j$, and let $a,b,c,d\in M$.  Then 
\bit
\itemit{i} $e_{ij;ab}e_{kj;cd}\sim_2e_{ij;ab}e_{kj;c1}$\emph{;}
\itemit{ii} $e_{ij;ab}e_{kj;c1}e_{kj;d1}\sim_2 e_{ij;ab}e_{kj;cd,1}$.
\eit
\end{lemma}

\pf (i).  We have $e_{ij;ab}e_{kj;cd} \sim_2 e_{ij;ab}e_{kj;11}e_{kj;cd} 
\sim_2 e_{ij;ab}e_{kj;cc} 
\sim_2 e_{ij;ab}e_{kj;11}e_{kj;c1} 
\sim_2 e_{ij;ab}e_{kj;c1}$,
using~\ref{R12} in each step and also \ref{R3a2} in the first and last.

\pfitem{ii}  Here we have $e_{ij;ab}e_{kj;c1}e_{kj;d1} \sim_2 e_{ij;ab}e_{kj;cd,d} \sim_2 e_{ij;ab}e_{kj;cd,1}$, by \ref{R12} and Part (i). \epf

\pf[\bf Proof of Theorem \ref{thm:X2}.]
It remains to show that the relations from $R_2$ imply those from $R_n'$.
The relations from $R_n'$ all involve at least one word of the form $E_{ij;\ba}$.  Recall that if $a_k=1$ for all $k\in\bn\sm\{i,j\}$, the word $E_{ij;\ba}$ is simply defined to be $e_{ij;a_ia_j}$.  But in this case, if $\bn\sm\{i,j\}=\{i_1,\ldots,i_{n-2}\}$, with $i_1<\cdots<i_{n-2}$, then
\[
E_{ij;\ba}=e_{ij;a_ia_j} \sim_2 e_{ij;a_ia_j}e_{i_1j;11}\cdots e_{i_{n-2}j;11} = e_{ij;a_ia_j}e_{i_1j;a_{i_1}1}\cdots e_{i_{n-2}j;a_{i_{n-2}}1},
\]
by \ref{R3a2}.  So, for uniformity, it will be convenient to assume that $E_{ij;\ba}$ is always given by the longer expression.  We now consider the relations from $R_n'$ one by one.

\bigskip\noindent
{\bf\boldmath \ref{R1n'}:} Here we have
$
E_{ij;\ba}= e_{ij;a_ia_j}e_{i_1j;a_{i_1}1}\cdots e_{i_{n-2}j;a_{i_{n-2}}1} \sim_2 e_{ij;a_ia_j}e_{i_1j;a_{i_1}1}\cdots e_{i_{n-2}j;a_{i_{n-2}}1} e_{ij;11} = E_{ij;\ba}e_{ij},
$
by \ref{R3a2}, and
$
E_{ij;\ba}= e_{ij;a_ia_j}e_{i_1j;a_{i_1}1}\cdots e_{i_{n-2}j;a_{i_{n-2}}1} \sim_2 
e_{ji;a_ja_i}e_{i_1i;a_{i_1}1}\cdots e_{i_{n-2}i;a_{i_{n-2}}1}e_{ij;11} = E_{ji;\ba}e_{ij},
$
by Lemma \ref{lem:chain}.

\bigskip\noindent
Now that we know \ref{R1n'} holds (modulo $R_2$), we can quickly derive \ref{R3n'}, \ref{R5n'} and \ref{R6n'}.

\bigskip\noindent
{\bf\boldmath \ref{R3n'}:}  Here we have $E_{ik:\ba}e_{jk} \sim_2 E_{ik:\ba}e_{ik}e_{jk} \sim E_{ik:\ba}e_{ik} \sim E_{ik:\ba}$, using \ref{R1n'} and \ref{R3a2}.

\bigskip\noindent
{\bf\boldmath \ref{R5n'}:} Here we have
$
E_{ki;\ba}e_{ij}e_{jk} \sim_2 E_{ki;\ba}e_{ki}e_{ij}e_{jk} \sim_2 E_{ki;\ba}e_{ik}e_{kj}e_{ji}e_{ik} \sim_2 E_{ik;\ba}e_{kj}e_{ji}e_{ik},
$
by \ref{R1n'} and~\ref{R52}.

\bigskip\noindent
{\bf\boldmath \ref{R6n'}:} This is almost identical to \ref{R5n'}.

\bigskip\noindent
{\bf\boldmath \ref{R2n'}:}  Relabelling the elements of $\bn$ if necessary, and using Corollary \ref{cor:Eij}, it suffices to assume that ${(i,j,k,l)=(1,2,3,4)}$, and to show that $E_{12;\ba}e_{34}\sim_2 E_{34;\ba}  e_{12}$.
Here we have
\begin{align*}
E_{12;\ba}e_{34} &= e_{12;a_1a_2}e_{32;a_31}e_{42;a_41}\cdot e_{52;a_51}\cdots e_{n2;a_n1} \cdot e_{34;11} \\
&\sim_2 e_{12;a_1a_2}e_{32;a_31}e_{42;a_41} e_{34;11} \cdot e_{52;a_51}\cdots e_{n2;a_n1} &&\text{by \ref{R22}}\\
&\sim_2 e_{12;a_1a_2}e_{32;a_31}e_{42;a_41} e_{34;11}  e_{14;11}e_{14;11} \cdot e_{52;a_51}\cdots e_{n2;a_n1} &&\text{by \ref{R3a2}}\\
&\sim_2 e_{12;a_1a_2}e_{32;a_31} e_{34;1a_4}e_{32;11}   e_{14;11}e_{14;11} \cdot e_{52;a_51}\cdots e_{n2;a_n1} &&\text{by \ref{R4b2}}\\
&\sim_2 e_{12;a_1a_2}e_{34;a_3a_4}e_{32;11} e_{32;11}   e_{14;11}e_{14;11} \cdot e_{52;a_51}\cdots e_{n2;a_n1} &&\text{by \ref{R4a2}}\\
&\sim_2 e_{34;a_3a_4}e_{12;a_1a_2}e_{32;11} e_{32;11}   e_{14;11}e_{14;11} \cdot e_{52;a_51}\cdots e_{n2;a_n1} &&\text{by \ref{R22}}\\
&\sim_2 e_{34;a_3a_4}e_{12;a_1a_2}   e_{14;11}e_{14;11} \cdot e_{52;a_51}\cdots e_{n2;a_n1} &&\text{by \ref{R3a2}}\\
&\sim_2 e_{34;a_3a_4}e_{14;a_11}e_{12;1a_2}   e_{14;11} \cdot e_{52;a_51}\cdots e_{n2;a_n1} &&\text{by \ref{R4a2}}\\
&\sim_2 e_{34;a_3a_4}e_{14;a_11}e_{24;a_21} e_{12;11}    \cdot e_{52;a_51}\cdots e_{n2;a_n1} &&\text{by \ref{R4a2}}\\
&\sim_2 e_{34;a_3a_4}e_{14;a_11}e_{24;a_21} e_{12;11} \cdot e_{42;11}   \cdot e_{52;a_51}\cdots e_{n2;a_n1} &&\text{by \ref{R3a2}}\\
&\sim_2 e_{34;a_3a_4}e_{14;a_11}e_{24;a_21} e_{12;11} \cdot e_{24;11}   \cdot e_{54;a_51}\cdots e_{n4;a_n1} \cdot e_{42;11} &&\text{by Lemma \ref{lem:chain}}\\
&\sim_2 e_{34;a_3a_4}e_{14;a_11}e_{24;a_21}e_{34;11} e_{12;11} \cdot e_{24;11}   \cdot e_{54;a_51}\cdots e_{n4;a_n1} \cdot e_{42;11} &&\text{by \ref{R3a2}}\\
&\sim_2 e_{34;a_3a_4}e_{14;a_11}e_{24;a_21} e_{12;11} e_{34;11}\cdot e_{24;11}   \cdot e_{54;a_51}\cdots e_{n4;a_n1} \cdot e_{42;11} &&\text{by \ref{R22}}\\
&\sim_2 e_{34;a_3a_4}e_{14;a_11}e_{24;a_21} e_{12;11} e_{34;11}   \cdot e_{54;a_51}\cdots e_{n4;a_n1} \cdot e_{42;11} &&\text{by \ref{R3a2}}\\
&\sim_2 e_{34;a_3a_4}e_{14;a_11}e_{24;a_21}  e_{34;11} e_{12;11}   \cdot e_{54;a_51}\cdots e_{n4;a_n1} \cdot e_{42;11} &&\text{by \ref{R22}}\\
&\sim_2 e_{34;a_3a_4}e_{14;a_11}e_{24;a_21}   e_{12;11}   \cdot e_{54;a_51}\cdots e_{n4;a_n1} \cdot e_{42;11} &&\text{by \ref{R3a2}}\\
&\sim_2 e_{34;a_3a_4}e_{14;a_11}e_{24;a_21}    \cdot e_{54;a_51}\cdots e_{n4;a_n1} \cdot  e_{12;11} e_{42;11} &&\text{by \ref{R22}}\\
&\sim_2 e_{34;a_3a_4}e_{14;a_11}e_{24;a_21}    \cdot e_{54;a_51}\cdots e_{n4;a_n1} \cdot  e_{12;11} = E_{34;\ba}  e_{12} &&\text{by \ref{R3a2}.}
\end{align*}
{\bf\boldmath \ref{R4n'}:} As in the previous calculation, we may assume $(i,j,k)=(1,2,3)$.  First,
\begin{align*}
E_{12;\ba}e_{13} &= e_{12;a_1a_2}e_{32;a_31}\cdot e_{42;a_41}\cdots e_{n2;a_n1}\cdot e_{13;11} \\
&\sim_2 e_{12;a_1a_2}e_{32;a_31} e_{13;11}\cdot e_{42;a_41}\cdots e_{n2;a_n1} &&\text{by \ref{R22}}\\
&\sim_2 e_{12;a_1a_2}e_{13;1a_3}e_{12;11} \cdot e_{42;a_41}\cdots e_{n2;a_n1} &&\text{by \ref{R4b2}}\\
&\sim_2 e_{13;a_1a_3}e_{12;1a_2}e_{12;11} \cdot e_{42;a_41}\cdots e_{n2;a_n1} &&\text{by \ref{R4a2}}\\
&\sim_2 e_{13;a_1a_3}e_{13;11}e_{21;a_21}e_{12;11} \cdot e_{42;a_41}\cdots e_{n2;a_n1} &&\text{by \ref{R12} twice}\\
&\sim_2 e_{13;a_1a_3}e_{23;a_21}e_{21;11}e_{12;11} \cdot e_{42;a_41}\cdots e_{n2;a_n1} &&\text{by \ref{R4b2}}\\
&\sim_2 e_{13;a_1a_3}e_{23;a_21}e_{12;11} \cdot e_{42;a_41}\cdots e_{n2;a_n1} &&\text{by \ref{R12}}\\
&\sim_2 e_{13;a_1a_3}e_{23;a_21}e_{12;11} \cdot e_{32;11} \cdot e_{42;a_41}\cdots e_{n2;a_n1} &&\text{by \ref{R3a2}}\\
&\sim_2 e_{13;a_1a_3}e_{23;a_21}e_{12;11} \cdot e_{23;11} \cdot e_{43;a_41}\cdots e_{n3;a_n1} \cdot e_{32;11} &&\text{by Lemma \ref{lem:chain}}\\
&\sim_2 e_{13;a_1a_3}e_{23;a_21}e_{43;11}e_{12;11} \cdot e_{23;11} \cdot e_{43;a_41}\cdots e_{n3;a_n1} \cdot e_{32;11} &&\text{by \ref{R3a2}}\\
&\sim_2 e_{13;a_1a_3}e_{23;a_21}e_{12;11}e_{43;11} \cdot e_{23;11} \cdot e_{43;a_41}\cdots e_{n3;a_n1} \cdot e_{32;11} &&\text{by \ref{R22}}\\
&\sim_2 e_{13;a_1a_3}e_{23;a_21}e_{12;11}e_{43;11}  \cdot e_{43;a_41}\cdots e_{n3;a_n1} \cdot e_{32;11} &&\text{by \ref{R3a2}}\\
&\sim_2 e_{13;a_1a_3}e_{23;a_21}e_{43;11}e_{12;11}  \cdot e_{43;a_41}\cdots e_{n3;a_n1} \cdot e_{32;11} &&\text{by \ref{R22}}\\
&\sim_2 e_{13;a_1a_3}e_{23;a_21}e_{12;11}  \cdot e_{43;a_41}\cdots e_{n3;a_n1} \cdot e_{32;11} &&\text{by \ref{R3a2}}\\
&\sim_2 e_{13;a_1a_3}e_{23;a_21}  \cdot e_{43;a_41}\cdots e_{n3;a_n1} \cdot e_{12;11}e_{32;11} &&\text{by \ref{R22}}\\
&\sim_2 e_{13;a_1a_3}e_{23;a_21}  \cdot e_{43;a_41}\cdots e_{n3;a_n1} \cdot e_{12;11} = E_{13;\ba}e_{12} &&\text{by \ref{R3a2},}
\intertext{establishing the first part of \ref{R4n'}.  For the second part,}
E_{13;\ba}e_{12} &= e_{13;a_1a_3}e_{23;a_21}\cdot e_{43;a_41}\cdots e_{n3;a_n1}\cdot e_{12;11} \\
&\sim_2 e_{13;a_1a_3}e_{23;a_21}e_{12;11}\cdot e_{43;a_41}\cdots e_{n3;a_n1} &&\text{by \ref{R22}}\\
&\sim_2 e_{13;a_1a_3}e_{12;1a_2}e_{13;11}\cdot e_{43;a_41}\cdots e_{n3;a_n1} &&\text{by \ref{R4b2}}\\
&\sim_2 e_{23;a_2a_3}e_{12;a_11}e_{13;11}\cdot e_{43;a_41}\cdots e_{n3;a_n1} &&\text{by \ref{R4b2}}\\
&\sim_2 e_{23;a_2a_3}e_{13;a_11}e_{12;11}\cdot e_{43;a_41}\cdots e_{n3;a_n1} &&\text{by \ref{R4a2}}\\
&\sim_2 e_{23;a_2a_3}e_{13;a_11}\cdot e_{43;a_41}\cdots e_{n3;a_n1} \cdot e_{12;11} = E_{23;\ba}e_{12} &&\text{by \ref{R22}.}
\end{align*}
{\bf\boldmath \ref{R7n'}:}  Here we may assume that $(i,j)=(1,2)$, but (because $i,j,k,l$ are not necessarily distinct), we need to consider several possibilities for the values of $k,l$:
\begin{itemize}
\begin{multicols}{4}
\item[(i)] $(k,l)=(1,2)$;
\item[(ii)] $(k,l)=(2,1)$;
\item[(iii)] $(k,l)=(1,3)$;
\item[(iv)] $(k,l)=(3,1)$;
\item[(v)] $(k,l)=(3,2)$;
\item[(vi)] $(k,l)=(2,3)$;
\item[(vii)] $(k,l)=(3,4)$.
\end{multicols}
\end{itemize}
In each case, we have $\bc = (c_1,\ldots,c_n) = \ba(\ve_{12}\cdot\bb) = (a_1b_1,a_2b_1,a_3b_3,\ldots,a_nb_n)$.  
In Case (i), we have
\begin{align*}
E_{12;\ba}E_{12;\bb} &= e_{12;a_1a_2}e_{32;a_31}\cdots e_{n2;a_n1} e_{12;b_1b_2}e_{32;b_31}\cdots e_{n2;b_n1} \\
&\sim_2 e_{12;a_1a_2}e_{32;a_31}\cdots e_{n2;a_n1} e_{12;b_11}e_{32;b_31}\cdots e_{n2;b_n1} &&\text{by Lemma \ref{lem:d1}(i)}\\
&\sim_2 e_{12;a_1a_2}e_{12;b_11} (e_{32;a_31}e_{32;b_31})\cdots (e_{n2;a_n1}e_{n2;b_n1})  &&\text{by Corollary \ref{cor:order}, repeatedly}\\
&\sim_2 e_{12;a_1b_1,a_2b_1} e_{32;a_3b_3,1}\cdots e_{n2;a_nb_n,1}  &&\text{by \ref{R12} and Lemma \ref{lem:d1}(ii)}\\
&\sim_2 e_{12;a_1b_1,a_2b_1} e_{32;a_3b_3,1}\cdots e_{n2;a_nb_n,1} e_{12;11}  = E_{12;\bc}e_{12} &&\text{by \ref{R3a2}.}
\end{align*}
For Case (ii), we use \ref{R1n'}, Case (i) of \ref{R7n'}, and \ref{R12}:
\[
E_{12;\ba}E_{21;\bb} \sim_2 E_{12;\ba}E_{12;\bb}e_{21} \sim_2 E_{12;\bc}e_{12}e_{21}  \sim_2 E_{12;\bc}e_{21}.
\]
In Case (iii), we have
\begin{align*}
E_{12;\ba}&E_{13;\bb} 
= e_{12;a_1a_2}e_{32;a_31}e_{42;a_41}\cdots e_{n2;a_n1} e_{13;b_1b_3}e_{23;b_21}e_{43;b_41}\cdots e_{n3;b_n1} \\
&\sim_2 e_{12;a_1a_2}e_{32;a_31}e_{13;b_1b_3}e_{42;a_41}\cdots e_{n2;a_n1} e_{23;b_21}e_{43;b_41}\cdots e_{n3;b_n1} &&\text{by \ref{R22}}\\
&\sim_2 e_{12;a_1a_2}e_{32;a_31}e_{13;b_1b_3}e_{42;a_41}\cdots e_{n2;a_n1} e_{32;1b_2}e_{42;b_41}\cdots e_{n2;b_n1}e_{23;11} &&\text{by Lemma \ref{lem:chain}}\\
&\sim_2 e_{12;a_1a_2}e_{12;b_1b_3}e_{13;1,a_3b_3}e_{42;a_41}\cdots e_{n2;a_n1} e_{32;1b_2}e_{42;b_41}\cdots e_{n2;b_n1}e_{23;11} &&\text{by \ref{R4b2}}\\
&\sim_2 e_{12;a_1b_1,a_2b_1}e_{13;1,a_3b_3}e_{42;a_41}\cdots e_{n2;a_n1} e_{32;1b_2}e_{42;b_41}\cdots e_{n2;b_n1}e_{23;11} &&\text{by \ref{R12}}\\
\tag{$*$}\label{eq:*} &\sim_2 e_{12;a_1b_1,a_2b_1}e_{13;1,a_3b_3}e_{42;a_41}\cdots e_{n2;a_n1} e_{42;b_41}\cdots e_{n2;b_n1}e_{23;11} &&\text{by \ref{R3a2} unless $n=3$ }\\
&\sim_2 e_{12;a_1b_1,a_2b_1}e_{42;a_41}\cdots e_{n2;a_n1} e_{42;b_41}\cdots e_{n2;b_n1}e_{13;1,a_3b_3}e_{23;11} &&\text{by \ref{R22}}\\
&\sim_2 e_{12;a_1b_1,a_2b_1}(e_{42;a_41}e_{42;b_41})\cdots (e_{n2;a_n1}e_{n2;b_n1}) e_{13;1,a_3b_3}e_{23;11} &&\text{by Corollary \ref{cor:order}}\\
&\sim_2 e_{12;a_1b_1,a_2b_1}e_{42;a_4b_4,1}\cdots e_{n2;a_nb_n,1} e_{13;1,a_3b_3}e_{23;11} &&\text{by Lemma \ref{lem:d1}(ii)}\\
&\sim_2 e_{12;a_1b_1,a_2b_1}e_{42;a_4b_4,1}\cdots e_{n2;a_nb_n,1} e_{12;11}e_{13;1,a_3b_3}e_{23;11} &&\text{by \ref{R3a2}}\\
&\sim_2 e_{12;a_1b_1,a_2b_1}e_{42;a_4b_4,1}\cdots e_{n2;a_nb_n,1} e_{32;a_3b_3,1}e_{13;11}e_{23;11} &&\text{by \ref{R4b2}}\\
&\sim_2 e_{12;a_1b_1,a_2b_1}e_{32;a_3b_3,1}e_{42;a_4b_4,1}\cdots e_{n2;a_nb_n,1} e_{13;11} = E_{12;\bc}e_{13} &&\text{by Corollary \ref{cor:order} and \ref{R3a2}.}
\end{align*}
Note that at the step labelled \eqref{eq:*}, relation \ref{R3a2} does not apply if $n=3$, since then the word $e_{42;a_41}\cdots e_{n2;a_n1}$ is empty.  However, this step can still be accomplished, albeit by using more relations:
\begin{align*}
e_{12;a_1b_1,a_2b_1}e_{13;1,a_3b_3}e_{42;a_41}\cdots e_{n2;a_n1} e_{32;1b_2} &= e_{12;a_1b_1,a_2b_1}e_{13;1,a_3b_3} e_{32;1b_2} \\
&\sim_2 e_{13;a_1b_1,a_3b_3}e_{12;1,a_2b_1} e_{32;1b_2} &&\text{by \ref{R4a2}}\\
&\sim_2 e_{13;a_1b_1,a_3b_3}e_{12;1,a_2b_1}  &&\text{by \ref{R3a2}}\\
&\sim_2 e_{12;a_1b_1,a_2b_1}e_{13;1,a_3b_3}  &&\text{by \ref{R4a2}}\\
&= e_{12;a_1b_1,a_2b_1}e_{13;1,a_3b_3}e_{42;a_41}\cdots e_{n2;a_n1}.
\end{align*}
Again, we may deduce Case (iv) from Case (iii), together with \ref{R1n'} and \ref{R12}:
\[
E_{12;\ba}E_{31;\bb} \sim_2 E_{12;\ba}E_{13;\bb}e_{31} \sim_2 E_{12;\bc}e_{13}e_{31} \sim_2 E_{12;\bc}e_{31}.
\]
In Case (v), we have
\begin{align*}
E_{12;\ba}E_{32;\bb} &= e_{12;a_1a_2}e_{32;a_31}e_{42;a_41}\cdots e_{n2;a_n1} e_{32;b_3b_2}e_{12;b_11}e_{42;b_41}\cdots e_{n2;b_n1} \\
&\sim_2 e_{12;a_1a_2}e_{32;a_31}e_{42;a_41}\cdots e_{n2;a_n1} e_{32;b_31}e_{12;b_11}e_{42;b_41}\cdots e_{n2;b_n1} &&\text{by Lemma \ref{lem:d1}(i)}\\
&\sim_2 e_{12;a_1a_2}e_{12;b_11} (e_{32;a_31}e_{32;b_31})(e_{42;a_41}e_{42;b_41})\cdots (e_{n2;a_n1}e_{n2;b_n1})   &&\text{by Corollary \ref{cor:order}}\\
&\sim_2 e_{12;a_1b_1,a_2b_1} e_{32;a_3b_3,1}e_{42;a_4b_4,1}\cdots e_{n2;a_nb_n,1}   &&\text{by \ref{R12} and Lemma \ref{lem:d1}(ii)}\\
&\sim_2 e_{12;a_1b_1,a_2b_1} e_{32;a_3b_3,1}e_{42;a_4b_4,1}\cdots e_{n2;a_nb_n,1}e_{32;11} = E_{12;\bc}e_{32}   &&\text{by \ref{R3a2}.}
\end{align*}
Case (vi) follows quickly from Case (v), \ref{R1n'} and \ref{R12}.  Finally, for Case (vii), we first observe that for any $u,v,x,y\in M$,
\begin{align*}
e_{12;uv}e_{32;x1}e_{42;y1}e_{34;11} &\sim_2 e_{12;uv}e_{32;x1}e_{32;11}e_{34;1y} &&\text{by \ref{R4b2}}\\
&\sim_2 e_{12;uv}e_{32;x1}e_{34;1y} &&\text{by \ref{R12}}\\
&\sim_2 e_{12;uv}e_{34;xy}e_{32;11} &&\text{by \ref{R4a2}}\\
&\sim_2 e_{34;xy}e_{12;uv}e_{32;11} &&\text{by \ref{R22}}\\
&\sim_2 e_{34;xy}e_{12;uv} &&\text{by \ref{R3a2}}\\
&\sim_2 e_{12;uv}e_{34;xy} &&\text{by \ref{R22}.}
\end{align*}
We then calculate
\begin{align*}
&E_{12;\ba}E_{34;\bb} = e_{12;a_1a_2}e_{32;a_31}e_{42;a_41}e_{52;a_51}\cdots e_{n2;a_n1} e_{34;b_3b_4}e_{14;b_11}e_{24;b_21}e_{54;b_51}\cdots e_{n4;b_n1} \\
&\sim_2 e_{12;a_1a_2}e_{32;a_31}e_{42;a_41}e_{34;b_3b_4}e_{14;b_11}e_{52;a_51}\cdots e_{n2;a_n1} e_{24;b_21}e_{54;b_51}\cdots e_{n4;b_n1}
&&\text{by \ref{R22}}\\
&\sim_2 e_{12;a_1a_2}e_{32;a_31}e_{42;a_41}e_{34;b_3b_4}e_{14;b_11}e_{52;a_51}\cdots e_{n2;a_n1} e_{42;1b_2}e_{52;b_51}\cdots e_{n2;b_n1}e_{24;11} 
&&\text{by Lemma \ref{lem:chain}}\\
%
%
&\sim_2 e_{12;a_1a_2}e_{32;a_31}e_{34;b_3,a_4b_4}e_{32;1b_4}e_{14;b_11}e_{52;a_51}\cdots e_{n2;a_n1} e_{42;1b_2}e_{52;b_51}\cdots e_{n2;b_n1}e_{24;11} 
&&\text{by \ref{R4b2}}\\
&\sim_2 e_{12;a_1a_2}e_{34;a_3b_3,a_4b_4}e_{32;1,b_3}e_{14;b_11}e_{32;1b_4}e_{52;a_51}\cdots e_{n2;a_n1} e_{42;1b_2}e_{52;b_51}\cdots e_{n2;b_n1}e_{24;11} 
&&\text{by \ref{R22} and \ref{R4a2}}\\
&\sim_2 e_{34;a_3b_3,a_4b_4}e_{12;a_1a_2}e_{32;1,b_3}e_{14;b_11}e_{32;1b_4}e_{52;a_51}\cdots e_{n2;a_n1} e_{42;1b_2}e_{52;b_51}\cdots e_{n2;b_n1}e_{24;11} 
&&\text{by \ref{R22}}\\
&\sim_2 e_{34;a_3b_3,a_4b_4}e_{12;a_1a_2}e_{14;b_11}e_{32;1b_4}e_{52;a_51}\cdots e_{n2;a_n1} e_{42;1b_2}e_{52;b_51}\cdots e_{n2;b_n1}e_{24;11} 
&&\text{by \ref{R3a2}}\\
&\sim_2 e_{34;a_3b_3,a_4b_4}e_{14;a_1b_1,1}e_{12;1,a_2b_1}e_{32;1b_4}e_{52;a_51}\cdots e_{n2;a_n1} e_{42;1b_2}e_{52;b_51}\cdots e_{n2;b_n1}e_{24;11} 
&&\text{by \ref{R4a2}}\\
&\sim_2 e_{34;a_3b_3,a_4b_4}e_{14;a_1b_1,1}e_{12;1,a_2b_1}e_{52;a_51}\cdots e_{n2;a_n1} e_{52;b_51}\cdots e_{n2;b_n1}e_{24;11} 
&&\text{by \ref{R3a2} twice}\\
&\sim_2 e_{34;a_3b_3,a_4b_4}e_{12;a_1b_1,a_2b_1}e_{14;11}e_{52;a_51}\cdots e_{n2;a_n1} e_{52;b_51}\cdots e_{n2;b_n1}e_{24;11} 
&&\text{by \ref{R4a2}}\\
&\sim_2 e_{12;a_1b_1,a_2b_1}e_{34;a_3b_3,a_4b_4}e_{14;11}e_{52;a_51}\cdots e_{n2;a_n1} e_{52;b_51}\cdots e_{n2;b_n1}e_{24;11} 
&&\text{by \ref{R22}}\\
&\sim_2 e_{12;a_1b_1,a_2b_1}e_{34;a_3b_3,a_4b_4}e_{52;a_51}\cdots e_{n2;a_n1} e_{52;b_51}\cdots e_{n2;b_n1}e_{24;11} 
&&\text{by \ref{R3a2}}\\
&\sim_2 e_{34;a_3b_3,a_4b_4}e_{12;a_1b_1,a_2b_1}e_{52;a_51}\cdots e_{n2;a_n1} e_{52;b_51}\cdots e_{n2;b_n1}e_{24;11} 
&&\text{by \ref{R22}}\\
&\sim_2 e_{34;a_3b_3,a_4b_4}e_{12;a_1b_1,a_2b_1}(e_{52;a_51}e_{52;b_51})\cdots (e_{n2;a_n1}e_{n2;b_n1}) e_{24;11} 
&&\text{by Corollary \ref{cor:order}}\\
%
&\sim_2 e_{34;a_3b_3,a_4b_4}e_{12;a_1b_1,a_2b_1}e_{52;a_5b_5,1}\cdots e_{n2;a_nb_n,1} e_{24;11} 
&&\text{by Lemma \ref{lem:d1}(ii)}\\
&\sim_2e_{12;a_1b_1,a_2b_1}e_{52;a_5b_5,1}\cdots e_{n2;a_nb_n,1}  e_{34;a_3b_3,a_4b_4}e_{24;11} 
&&\text{by \ref{R22}}\\
&\sim_2e_{12;a_1b_1,a_2b_1}e_{52;a_5b_5,1}\cdots e_{n2;a_nb_n,1}  e_{34;a_3b_3,a_4b_4}
&&\text{by \ref{R3a2}}\\
&\sim_2e_{12;a_1b_1,a_2b_1}e_{34;a_3b_3,a_4b_4}e_{52;a_5b_5,1}\cdots e_{n2;a_nb_n,1}  
&&\text{by \ref{R22}}\\
&\sim_2e_{12;a_1b_1,a_2b_1}e_{32;a_3b_3,1}e_{42;a_4b_4,1}e_{34;11}e_{52;a_5b_5,1}\cdots e_{n2;a_nb_n,1}  
&&\text{by the observation}\\
&\sim_2e_{12;a_1b_1,a_2b_1}e_{32;a_3b_3,1}e_{42;a_4b_4,1}e_{52;a_5b_5,1}\cdots e_{n2;a_nb_n,1}  e_{34;11}
= E_{12;\bc}e_{34}
&&\text{by \ref{R22}.}
\end{align*}
We have finally shown that all the relations from $R_n'$ may be removed, and the proof is therefore complete.~\epf

\boldsubsection{An idempotent-generated  presentation for $M\wr\Sing_n$ with $M/{\L}$ a chain}\label{subsect:X1}

For the duration of this section, we fix a monoid $M$ for which $M/{\L}$ is a chain.  Recall from Theorem \ref{thm:IG} that, for such a monoid $M$, the singular wreath product $M\wr\Sing_n$ is generated by its idempotents: indeed, by the idempotents from the set $\X_1$.  It is therefore our goal in this section to obtain a presentation (see Theorem \ref{thm:X1}) for $M\wr\Sing_n$ in terms of the idempotent generating set $\X_1$.  In the special case that $M$ is a group, some of the relations take on a simpler form (see Theorem \ref{thm:X1group}).

With the stated goal in mind, define an alphabet
\[
X_1 = \set{e_{ij;a}}{i,j\in\bn,\ i\not=j,\ a\in M},
\]
an epimorphism
\[
\phi_1:X_1^+\to M\wr\Sing_n:e_{ij;a}\mt\ve_{ij;a},
\]
and let $R_1$ be the set of relations
\begin{align}
\label{R1a1}\tag*{(R1a)$_1$} e_{ij;a}e_{ij;b} &= e_{ij;a}  &&\text{for $a,b\in M$ and distinct $i,j$}\\
\label{R1b1}\tag*{(R1b)$_1$} e_{ij;1}e_{ji;a}e_{ij;b} &= e_{ji;1}e_{ij;ab} &&\text{for $a,b\in M$ and distinct $i,j$}\\
\label{R1c1}\tag*{(R1c)$_1$} e_{ji;a}e_{ij;c} &= e_{ji;b}e_{ij;c} &&\text{for $a,b,c\in M$ and distinct $i,j$ with $ac=bc$}\\
\label{R1d1}\tag*{(R1d)$_1$} e_{ij;b}e_{ji;c}e_{ij;1} &= e_{ji;a}e_{ij;bc} &&\text{for $a,b,c\in M$ and distinct $i,j$ with $abc=c$}\\
\label{R1e1}\tag*{(R1e)$_1$} e_{ji;1}e_{ij;1} &= e_{ij;1} &&\text{for distinct $i,j$}\\
\label{R21}\tag*{(R2)$_1$} e_{ij;a}e_{kl;b} &= e_{kl;b}e_{ij;a} &&\text{for $a,b\in M$ and distinct $i,j,k,l$}\\
\label{R3a1}\tag*{(R3a)$_1$} e_{ik;a}e_{jk;b} &= e_{ik;a} &&\text{for $a,b\in M$ and distinct $i,j,k$}\\
\label{R3b1}\tag*{(R3b)$_1$} e_{ij;1}e_{jk;a}e_{kj;1} &= e_{ji;1}e_{ik;a}e_{ki;1}e_{ij;1} &&\text{for $a\in M$ and distinct $i,j,k$}\\
\label{R3c1}\tag*{(R3c)$_1$} e_{ij;1}e_{ji;a}e_{ik;b} &= e_{ji;1}e_{ik;b}e_{kj;a}e_{jk;1} &&\text{for $a,b\in M$ and distinct $i,j,k$}\\
\label{R41}\tag*{(R4)$_1$} e_{ij;b}e_{ik;ab} = e_{ik;ab}e_{ij;b} &= e_{jk;a}e_{ij;b} &&\text{for $a,b\in M$ and distinct $i,j,k$}\\
\label{R51}\tag*{(R5)$_1$} e_{ki}e_{ij}e_{jk} &= e_{ik}e_{kj}e_{ji}e_{ik} &&\text{for distinct $i,j,k$}\\
\label{R61}\tag*{(R6)$_1$} e_{ki}e_{ij}e_{jk}e_{kl} &= e_{ik}e_{kl}e_{li}e_{ij}e_{jl} &&\text{for distinct $i,j,k,l$.}
\end{align}
Again, we have identified $X$ with a subset of $X_1$ in relations \ref{R51} and \ref{R61}.  
It is important to note that some relations (such as \ref{R1b1}) involve a letter of the form $e_{ij;ab}$ from $X_1$, where ``$ab$'' denotes a single subscript (the product of $a$ and $b$ in $M$): in particular, $e_{ij;ab}$ does not represent the letter from $X_2$ where $a$ and $b$ are separate subscripts.  In order to avoid any potential confusion, we will write \emph{all} letters from $X_2$ as $e_{ij;a,b}$ (not just those for which the monoid subscripts are expressed as products) for the duration of this section.
%
We may now state the main result of this section.

\begin{thm}\label{thm:X1}
If $M/{\L}$ is a chain, then the semigroup $M\wr\Sing_n$ has presentation $\pres{X_1}{R_1}$ via $\phi_1$.  
\end{thm}

To prove Theorem \ref{thm:X1}, we begin with the presentation $\pres{X_2}{R_2}$ from Theorem~\ref{thm:X2}.  Again, we think of~$X_1$ as a subset of $X_2$ by identifying $e_{ij;a}\in X_1$ with $e_{ij;1,a}\in X_2$.  Since $u\phi_2=v\phi_2$ for all $(u,v)\in R_1$, as may easily be checked diagrammatically, we may add relations $R_1$ to the presentation to obtain $\pres{X_2}{R_1\cup R_2}$.  Now write ${\sim_1}=R_1^\sharp$ for the congruence on $X_1^+$ generated by $R_1$.  


Since $M/{\L}$ is a chain, we may fix some subset $\Om\sub M\times M$ such that the following conditions are satisfied:~
\bit
\item[(i)] for all $a,b\in M$, $\Om$ contains exactly one of $(a,b)$ or $(b,a)$; and
\item[(ii)] for all $(a,b)\in\Om$, $a\leq_{\L} b$.
\eit
Note that $\Om$ (regarded as a binary relation) is reflexive and anti-symmetric, but need not be transitive (though $\Om$ could be chosen to have this additional property, in which case $\Om$ would then be a total order on~$M$ that refines the preorder $\leq_{\L}$).  
For each $(a,b)\in\Om$, we choose some $x_{ab}\in M$
 such that $a=x_{ab}b$.  
 For each $i,j\in\bn$ with $i\not=j$, and for each $a,b\in M$, we define the word
\[
E_{ij;ab} = \begin{cases}
e_{ji;x_{ab}}e_{ij;b} &\text{if $(a,b)\in\Om$}\\
e_{ij;x_{ba}}e_{ji;a}e_{ij;1} &\text{if $(a,b)\not\in\Om$.}
\end{cases}
\]
(Since there is no chance of confusion, we will not generally need to write $E_{ij;a,b}$ for these words.)  We begin with a simple lemma.

\begin{lemma}\label{lem:Eijab}
Let $i,j\in\bn$ with $i\not=j$, and let $a,b\in M$.  Then
\bit
\itemit{i} $E_{ij;ab}\phi_2=e_{ij;a,b}\phi_2$\emph{;}
\itemit{ii} $E_{ij;ab}\sim_1 e_{ji;x}e_{ij;b}$ for any $x\in M$ with $a=xb$\emph{;}
\itemit{iii} $E_{ij;ab}\sim_1 e_{ij;x}e_{ji;a}e_{ij;1}$ for any $x\in M$ with $b=xa$\emph{;}
\itemit{iv} $E_{ij;1a}\sim_1e_{ij;a}$.
\eit
\end{lemma}

\pf  (i).  This is easily checked diagrammatically.  

\pfitem{ii}  We must consider two cases.  If $(a,b)\in\Om$, then \ref{R1c1} gives $E_{ij;ab}=e_{ji;x_{ab}}e_{ij;b}\sim_1 e_{ji;x}e_{ij;b}$, since $xb=a=x_{ab}b$.  If $(a,b)\not\in\Om$, then $(b,a)\in\Om$, $b=x_{ba}a$ and $a=xb=xx_{ba}a$, so \ref{R1d1} gives $E_{ij;ab}=e_{ij;x_{ba}}e_{ji;a}e_{ij;1}\sim_1 e_{ji;x}e_{ij;x_{ba}a}=e_{ji;x}e_{ij;b}$.  

\pfitem{iii}  If $(a,b)\in\Om$, then $a=x_{ab}b=x_{ab}xa$, and \ref{R1d1} gives $E_{ij;ab}=e_{ji;x_{ab}}e_{ij;b}=e_{ji;x_{ab}}e_{ij;xa}\sim_1e_{ij;x}e_{ji;a}e_{ij;1}$.  If $(a,b)\not\in\Om$, then \ref{R1c1} gives $E_{ij;ab}=e_{ij;x_{ba}}e_{ji;a}e_{ij;1}\sim_1e_{ij;x}e_{ji;a}e_{ij;1}$. 

\pfitem{iv}  We have $E_{ij;1a} \sim_1 e_{ij;a}e_{ji;1}e_{ij;1} \sim_1 e_{ij;a}e_{ij;1} \sim_1 e_{ij;a}$, using Part (iii) and then \ref{R1e1} and \ref{R1a1}. \epf



By Lemma \ref{lem:Eijab}(i) and Theorem \ref{thm:X2}, it follows that $e_{ij;a,b}\sim_2 E_{ij;ab}$ for each $i,j,a,b$, so we may remove each generator $e_{ij;a,b}\in X_2\sm X_1$ from the presentation $\pres{X_2}{R_1\cup R_2}$, replacing every occurence of such a generator in the relations by the word $E_{ij;ab}\in X_1^+$.  The only relations that are affected in this way are those from $R_2$, excluding \ref{R52} and \ref{R62}, which are just \ref{R51} and \ref{R61}.  We denote the relations modified in this way by:
\begin{align}
\label{R12'}\tag*{(R1)$'_2$} E_{ij;ab}E_{ij;cd} = E_{ij;ac,bc} &= E_{ji;ba}E_{ij;dc} &&\text{for $a,b,c,d\in M$ and distinct $i,j$}\\
\label{R22'}\tag*{(R2)$'_2$} E_{ij;ab}E_{kl;cd} &= E_{kl;cd}E_{ij;ab} &&\text{for $a,b,c,d\in M$ and distinct $i,j,k,l$}\\
\label{R3a2'}\tag*{(R3a)$'_2$} E_{ik;ab}E_{jk;1c} &= E_{ik;ab} &&\text{for $a,b,c\in M$ and distinct $i,j,k$}\\
\label{R3b2'}\tag*{(R3b)$'_2$} E_{ik;ab}E_{jk;c1} &= E_{ki;ba}E_{ji;c1}E_{ik;11} &&\text{for $a,b,c\in M$ and distinct $i,j,k$}\\
\label{R3c2'}\tag*{(R3c)$'_2$} E_{ik;aa}E_{jk;b1} &= E_{ik;11}E_{jk;b1}E_{ik;a1} &&\text{for $a,b\in M$ and distinct $i,j,k$}\\
\label{R4a2'}\tag*{(R4a)$'_2$} E_{ij;ab}E_{ik;cd} = E_{ik;ac,d}E_{ij;1,bc} &= E_{jk;bc,d}E_{ij;ac,1} &&\text{for $a,b,c,d\in M$ and distinct $i,j,k$}\\
\label{R4b2'}\tag*{(R4b)$'_2$} E_{ij;c,ad}E_{ik;1,bd} = E_{ik;c,bd}E_{ij;1,ad} &= E_{jk;ab}E_{ij;cd} &&\text{for $a,b,c,d\in M$ and distinct $i,j,k$,}
\end{align}
and denote the whole set of relations modified in this way by $R_2'$.  So the presentation has now become $\pres{X_1}{R_1\cup R_2'}$, and we must show that the relations from $R_2'$ may be removed.  

We now pause to prove a technical lemma that will be useful on a number of occasions.  For its proof, and for many calculations in this section, we observe that the first part of \ref{R41} implies
\[
e_{ij;a}e_{ik;b} \sim_1 e_{ik;b}e_{ij;a} \qquad\text{for $a,b\in M$ and distinct $i,j,k$,}
\]
since we either have $a=xb$ or $b=xa$ for some $x\in M$, as $M/{\L}$ is a chain.

\begin{lemma}\label{lem:X1tech}
Let $i,j,k\in\bn$ be distinct, and let $a,b,c,d\in M$.  Then
\bit
\itemit{i} $e_{jk;a}e_{kj;b}e_{ki;c}e_{ik;d} \sim_1 e_{jk;a}e_{kj;bd}e_{ki;cd}e_{ik;1}$\emph{;}
\itemit{ii} $e_{ki;1}e_{ik;b}e_{ij;ab}\sim_1e_{ik;1}e_{ij;a}e_{ji;b}e_{ij;1}$\emph{;}
\itemit{iii} $e_{ki;1}e_{ik;ab}e_{ij;b}\sim_1e_{ik;1}e_{ji;a}e_{ij;b}$\emph{;}
\itemit{iv} $e_{ki;a}e_{ik;c}e_{ij;d}\sim_1e_{ji;b}e_{ik;c}e_{ij;d}$ if $ac=bd$.
\eit
\end{lemma}

\pf 
(i).  We have
\begin{align*}
e_{jk;a}e_{kj;b}e_{ki;c}e_{ik;d} &\sim_1 e_{jk;a}e_{ki;c}e_{kj;b}e_{ik;d} &&\text{by \ref{R41}}\\
&\sim_1 e_{jk;a}e_{ik;1}e_{ki;c}e_{ik;d}e_{ij;bd} &&\text{by \ref{R3a1} and \ref{R41}}\\
&\sim_1 e_{jk;a}e_{ki;1}e_{ik;cd}e_{ij;bd} &&\text{by \ref{R1b1}}\\
&\sim_1 e_{jk;a}e_{ik;1}e_{ki;cd}e_{ik;1}e_{ij;bd} &&\text{by \ref{R1b1}}\\
&\sim_1 e_{jk;a}e_{ki;cd}e_{ik;1}e_{ij;bd} &&\text{by \ref{R3a1}}\\
&\sim_1 e_{jk;a}e_{ki;cd}e_{kj;bd}e_{ik;1} &&\text{by \ref{R41}}\\
&\sim_1 e_{jk;a}e_{kj;bd}e_{ki;cd}e_{ik;1} &&\text{by \ref{R41}.}
\intertext{(ii).  We have}
e_{ik;1}e_{ij;a}e_{ji;b}e_{ij;1} &\sim_1 e_{kj;a}e_{ik;1}e_{ji;b}e_{ij;1} &&\text{by \ref{R41}}\\
&\sim_1 e_{kj;a}e_{ji;b}e_{jk;b}e_{ij;1} &&\text{by \ref{R41}}\\
&\sim_1 e_{kj;a}e_{jk;b}e_{ji;b}e_{ij;1} &&\text{by \ref{R41}}\\
&\sim_1 e_{kj;a}e_{ki;1}e_{jk;b}e_{ij;1} &&\text{by \ref{R41}}\\
&\sim_1 e_{ki;1}e_{kj;a}e_{ij;1}e_{ik;b} &&\text{by \ref{R41} twice}\\
&\sim_1 e_{ki;1}e_{kj;a}e_{ik;b} &&\text{by \ref{R3a1}}\\
&\sim_1 e_{ki;1}e_{ik;b}e_{ij;ab} &&\text{by \ref{R41}.}
\intertext{(iii).  We have $e_{ki;1}e_{ik;ab}e_{ij;b} \sim_1 e_{ki;1}e_{jk;a}e_{ij;b} \sim_1 e_{jk;a}e_{ji;a}e_{ij;b} \sim_1 e_{ik;1}e_{ji;a}e_{ij;b}$, by three applications of~\ref{R41}. 
\endgraf\bigskip\noindent
(iv). We must consider two cases.  If $c=xd$ for some $x\in M$, then}
e_{ki;a}e_{ik;c}e_{ij;d} = e_{ki;a}e_{ik;xd}e_{ij;d} &\sim_1 e_{ki;a}e_{jk;x}e_{ij;d} &&\text{by \ref{R41}}\\
&\sim_1 e_{ji;ax}e_{jk;x}e_{ij;d} &&\text{by \ref{R41}}\\
&\sim_1 e_{ji;ax}e_{ik;xd}e_{ij;d} &&\text{by \ref{R41}}\\
&\sim_1 e_{ji;ax}e_{ij;d}e_{ik;xd} &&\text{by \ref{R41}}\\
&\sim_1 e_{ji;b}e_{ij;d}e_{ik;c} &&\text{by \ref{R1c1}, as $axd=ac=bd$ and $xd=c$}\\
&\sim_1 e_{ji;b}e_{ik;c}e_{ij;d} &&\text{by \ref{R41}}
\intertext{while if $d=yc$ for some $y\in M$, then}
e_{ji;b}e_{ik;c}e_{ij;d} = e_{ji;b}e_{ik;c}e_{ij;yc} &\sim_1 e_{ji;b}e_{ij;yc}e_{ik;c} &&\text{by \ref{R41}}\\
&\sim_1 e_{ji;b}e_{kj;y}e_{ik;c} &&\text{by \ref{R41}}\\
&\sim_1 e_{ki;by}e_{kj;y}e_{ik;c} &&\text{by \ref{R41}}\\
&\sim_1 e_{ki;by}e_{ij;yc}e_{ik;c} &&\text{by \ref{R41}}\\
&\sim_1 e_{ki;by}e_{ik;c}e_{ij;yc} &&\text{by \ref{R41}}\\
&\sim_1 e_{ki;a}e_{ik;c}e_{ij;d} &&\text{by \ref{R1c1}, as $byc=bd=ac$ and $yc=d$.} \qedhere
\end{align*}
\epf

We now have all we need to begin the proof of the main result of this section.

\pf[\bf Proof of Theorem \ref{thm:X1}.]
We consider the relations from $R_2'$ one at a time; each splits into several cases, depending on the relationship between the monoid subscripts in the $\leq_{\L}$ order.

\bigskip\noindent
{\bf\boldmath \ref{R12'}:}  There are four cases to consider:
\begin{itemize}
\begin{multicols}{2}
\item[(i)] $a\leq_{\L} b$ and $c\leq_{\L} d$;
\item[(ii)] $a\leq_{\L} b$ and $d\leq_{\L} c$;
\item[(iii)] $b\leq_{\L} a$ and $c\leq_{\L} d$;
\item[(iv)] $b\leq_{\L} a$ and $d\leq_{\L} c$.
\end{multicols}
\end{itemize}
In Case (i), writing $a=xb$ and $c=yd$, we have
\begin{align*}
E_{ij;ab}E_{ij;cd} &\sim_1 e_{ji;x}e_{ij;b} e_{ji;y}e_{ij;d} &&\text{by Lemma \ref{lem:Eijab}(ii)}\\
&\sim_1 e_{ji;x}e_{ji;1}e_{ij;b} e_{ij;1}e_{ji;y}e_{ij;d} &&\text{by \ref{R1a1} twice}\\
&\sim_1 e_{ji;x}e_{ji;1}e_{ij;b} e_{ji;1}e_{ij;c} &&\text{by \ref{R1b1}, noting that $c=yd$}\\
&\sim_1 e_{ji;x}e_{ij;1}e_{ji;b} e_{ij;c} &&\text{by \ref{R1b1}}\\
&\sim_1 e_{ji;x}e_{ji;1}e_{ij;bc} &&\text{by \ref{R1b1}}\\
&\sim_1 e_{ji;x}e_{ij;bc} &&\text{by \ref{R1a1}}\\
&\sim_1 E_{ij;ac,bc} &&\text{by Lemma \ref{lem:Eijab}(ii), noting that $xbc=ac$,}
\intertext{while}
E_{ji;ba}E_{ij;dc} &\sim_1 e_{ji;x}e_{ij;b}e_{ji;1} e_{ij;y}e_{ji;d}e_{ij;1} &&\text{by Lemma \ref{lem:Eijab}(iii)}\\
&\sim_1 e_{ji;x}e_{ij;b}e_{ij;1} e_{ji;c}e_{ij;1} &&\text{by \ref{R1b1}, noting that $c=yd$}\\
&\sim_1 e_{ji;x}e_{ji;1}e_{ij;b} e_{ji;c}e_{ij;1} &&\text{by \ref{R1a1} twice}\\
&\sim_1 e_{ji;x}e_{ij;1}e_{ji;bc}e_{ij;1} &&\text{by \ref{R1b1}}\\
&\sim_1 e_{ji;x}e_{ji;1}e_{ij;bc} &&\text{by \ref{R1b1}}\\
&\sim_1 e_{ji;x}e_{ij;bc} &&\text{by \ref{R1a1}}\\
&\sim_1 E_{ij;ac,bc} &&\text{by Lemma \ref{lem:Eijab}(ii).}
\intertext{In Case (ii), writing $a=xb$ and $d=yc$, we have}
 E_{ij;ab}E_{ij;cd} &\sim_1 e_{ji;x}e_{ij;b}e_{ij;y}e_{ji;c}e_{ij;1} &&\text{by Lemma \ref{lem:Eijab}(ii) and (iii)}\\
 &\sim_1 e_{ji;x}e_{ji;1}e_{ij;b}e_{ji;c}e_{ij;1} && \text{by \ref{R1a1} twice}\\
 &\sim_1 e_{ji;x}e_{ij;1}e_{ji;bc}e_{ij;1} && \text{by \ref{R1b1}}\\
 &\sim_1 e_{ji;x}e_{ji;1}e_{ij;bc} && \text{by \ref{R1b1}}\\
 &\sim_1 e_{ji;x}e_{ij;bc}, && \text{by \ref{R1a1}}\\
 &\sim_1 E_{ij;ac,bc} &&\text{by Lemma \ref{lem:Eijab}(ii),}\\
\intertext{while}
 E_{ji;ba}E_{ij;dc} &\sim_1 e_{ji;x}e_{ij;b}e_{ji;1}e_{ji;y}e_{ij;c} &&\text{by Lemma \ref{lem:Eijab}(ii) and (iii)}\\
 &\sim_1 e_{ji;x}e_{ji;1}e_{ij;b}e_{ji;1}e_{ij;c} && \text{by \ref{R1a1} twice}\\
 &\sim_1 e_{ji;x}e_{ji;1}e_{ij;bc} && \text{by \ref{R1b1} twice}\\
 &\sim_1 e_{ji;x}e_{ij;bc} && \text{by \ref{R1a1}}\\
 &\sim_1 E_{ij;ac,bc} &&\text{by Lemma \ref{lem:Eijab}(ii).}
 \intertext{In Case (iii), writing $b=xa$ and $c=yd$, we have}
 E_{ij;ab}E_{ij;cd}&\sim_1 e_{ij;x}e_{ji;a}e_{ij;1}e_{ji;y}e_{ij;d} &&\text{by Lemma \ref{lem:Eijab}(ii) and (iii)}\\
 &\sim_1 e_{ij;x}e_{ij;1}e_{ji;a}e_{ij;1}e_{ji;y}e_{ij;d} && \text{by \ref{R1a1}}\\
 &\sim_1 e_{ij;x}e_{ji;1}e_{ij;ayd} && \text{by \ref{R1b1} three times,}\\
\intertext{while}
 E_{ij;ac,bc}&\sim_1 e_{ij;x}e_{ji;ac}e_{ij;1} &&\text{by Lemma \ref{lem:Eijab}(iii)}\\
 &\sim_1 e_{ij;x}e_{ij;1}e_{ji;ac}e_{ij;1} && \text{by \ref{R1a1}}\\
 &\sim_1 e_{ij;x}e_{ji;1}e_{ij;ayd} && \text{by \ref{R1b1} and $c=yd$,}\\
\intertext{and}
 E_{ji;ba}E_{ij;dc}&\sim_1 e_{ij;x}e_{ji;a}e_{ij;y}e_{ji;d}e_{ij;1} &&\text{by Lemma \ref{lem:Eijab}(ii) and (iii)}\\
 &\sim_1 e_{ij;x}e_{ij;1}e_{ji;a}e_{ij;y}e_{ji;d}e_{ij;1} && \text{by \ref{R1a1}}\\
 &\sim_1 e_{ij;x}e_{ji;1}e_{ij;ayd} && \text{by \ref{R1b1} three times.}\\
\intertext{In Case (iv), writing $b=xa$ and $d=yc$, we have}
 E_{ij;ab}E_{ij;cd}&\sim_1 e_{ij;x}e_{ji;a}e_{ij;1}e_{ij;y}e_{ji;c}e_{ij;1} &&\text{by Lemma \ref{lem:Eijab}(iii)}\\
 &\sim_1 e_{ij;x}e_{ij;1}e_{ji;a}e_{ij;1}e_{ji;c}e_{ij;1} && \text{by \ref{R1a1} twice}\\
 &\sim_1 e_{ij;x}e_{ji;1}e_{ij;ac} && \text{by \ref{R1b1} three times,}\\
\intertext{while}
 E_{ij;ac,bc}&\sim_1 e_{ij;x}e_{ji;ac}e_{ij;1} &&\text{by Lemma \ref{lem:Eijab}(iii)}\\
 &\sim_1 e_{ij;x}e_{ij;1}e_{ji;ac}e_{ij;1} && \text{by \ref{R1a1}}\\
 &\sim_1 e_{ij;x}e_{ji;1}e_{ij;ac} && \text{by \ref{R1b1},}\\
\intertext{and}
 E_{ji;ba}E_{ij;dc}&\sim_1 e_{ij;x}e_{ji;a}e_{ji;y}e_{ij;c} &&\text{by Lemma \ref{lem:Eijab}(ii)}\\
 &\sim_1 e_{ij;x}e_{ij;1}e_{ji;a}e_{ij;c} && \text{by \ref{R1a1} twice}\\
 &\sim_1 e_{ij;x}e_{ji;1}e_{ij;ac} && \text{by \ref{R1b1}.}
\end{align*}
{\bf\boldmath \ref{R22'}:}  This relation follows immediately from \ref{R21}.  

\bigskip\noindent
{\bf\boldmath \ref{R3a2'}:}  If $a=xb$, then Lemma \ref{lem:Eijab} and \ref{R3a1} give $E_{ik;ab}E_{jk;1c} \sim_1 e_{ki;x}e_{ik;b}e_{jk;c} \sim_1 e_{ki;x}e_{ik;b} \sim_1 E_{ik;ab}$.  An almost identical calculation deals with the case in which $b=xa$.

\bigskip\noindent
{\bf\boldmath \ref{R3b2'}:}  
If $a=xb$, then 
\begin{align*}
 E_{ki;ba}E_{ji;c1}E_{ik;11} &\sim_1 e_{ki;x}e_{ik;b}e_{ki;1}e_{ij;c}e_{ji;1}e_{ik;1} &&\text{by Lemma \ref{lem:Eijab}(ii), (iii) and (iv)}\\
 &\sim_1 e_{ki;x}e_{ik;b}e_{ik;1}e_{kj;c}e_{jk;1} &&\text{by \ref{R3b1}}\\
 &\sim_1 e_{ki;x}e_{ik;b}e_{kj;c}e_{jk;1} &&\text{by \ref{R1a1}}\\
 &\sim_1 E_{ik;ab}E_{jk;c1} &&\text{by Lemma \ref{lem:Eijab}(ii).}
\intertext{If $b=xa$, then}
 E_{ik;ab}E_{jk;c1}&\sim_1 e_{ik;x}e_{ki;a}e_{ik;1}e_{kj;c}e_{jk;1} &&\text{by Lemma \ref{lem:Eijab}(ii) and (iii)}\\
 &\sim_1 e_{ik;x}e_{ki;a}e_{ki;1}e_{ij;c}e_{ji;1}e_{ik;1} &&\text{by \ref{R3b1}}\\
 &\sim_1 e_{ik;x}e_{ki;a}e_{ij;c}e_{ji;1}e_{ik;1} &&\text{by \ref{R1a1}}\\
 &\sim_1 E_{ki;ba}E_{ji;c1}E_{ik;11} &&\text{by Lemma \ref{lem:Eijab}(ii) and (iv).}
\intertext{{\bf\boldmath \ref{R3c2'}:}  Here we have}
 E_{ik;aa}E_{jk;b1} &\sim_1 e_{ki;1}e_{ik;a}e_{kj;b}e_{jk;1} &&\text{by Lemma \ref{lem:Eijab}(ii)}\\
 &\sim_1 e_{ik;1}e_{kj;b}e_{ji;a}e_{ij;1}e_{jk;1} &&\text{by \ref{R3c1}}\\
 &\sim_1 e_{ik;1}e_{kj;b}e_{kj;1}e_{ji;a}e_{ij;1}e_{jk;1} &&\text{by \ref{R1a1}}\\
 &\sim_1 e_{ik;1}e_{kj;b}e_{jk;1}e_{ki;a}e_{ik;1}e_{kj;1}e_{jk;1} &&\text{by \ref{R3b1}}\\
 &\sim_1 e_{ik;1}e_{kj;b}e_{jk;1}e_{ki;a}e_{ik;1}e_{jk;1} &&\text{by \ref{R1e1}}\\
 &\sim_1 e_{ik;1}e_{kj;b}e_{jk;1}e_{ki;a}e_{ik;1} &&\text{by \ref{R3a1}}\\
 &\sim_1 E_{ik;11}E_{jk;b1}E_{ik;a1} &&\text{by Lemma \ref{lem:Eijab}(ii) and (iv).}
\end{align*}
{\bf\boldmath \ref{R4a2'}:}  As explained in Remark \ref{rem:R4}, we only need to show that $E_{ij;ab}E_{ik;cd} \sim_1 E_{ik;ac,d}E_{ij;1,bc}$.  However, to do this, we must consider six cases:
\begin{itemize}
\begin{multicols}{3}
\item[(i)] $a\leq_{\L} b$ and $c\leq_{\L} d$;
\item[(ii)] $b\leq_{\L} a$ and $c\leq_{\L} d$;
\item[(iii)] $a\leq_{\L} b$ and $d\leq_{\L} ac$;
\item[(iv)] $b\leq_{\L} a$ and $d\leq_{\L} ac$;
\item[(v)] $a\leq_{\L} b$ and $ac\leq_{\L} d\leq_{\L} c$;
\item[(vi)] $b\leq_{\L} a$ and $ac\leq_{\L} d\leq_{\L} c$.
\end{multicols}
\end{itemize}
(Note that $ac\leq_{\L} d$ in Cases (i) and (ii), and that $d\leq_{\L} c$ in Cases (iii) and (iv).)
For Cases (i) and (ii), write $c=yd$, and note that $ac=(ay)d$, so that
\begin{align*}
E_{ik;ac,d}E_{ij;1,bc} &\sim e_{ki;ay}e_{ik;d}e_{ij;bc} &&\text{by Lemma \ref{lem:Eijab}(ii) and (iv)}\\
&= e_{ki;ay}e_{ik;d}e_{ij;byd} \\
&\sim e_{ki;ay}e_{kj;by}e_{ik;d} &&\text{by \ref{R41}.}\\
\intertext{In Case (i), writing $a=xb$, we have}
E_{ij;ab}E_{ik;cd} &\sim_1 e_{ji;x}e_{ij;b}e_{ki;y}e_{ik;d} &&\text{by Lemma \ref{lem:Eijab}(ii)}\\
&\sim_1 e_{ji;x}e_{ki;y}e_{kj;by}e_{ik;d} &&\text{by \ref{R41}}\\
&\sim_1 e_{ji;x}e_{kj;by}e_{ik;d} &&\text{by \ref{R3a1}}\\
&\sim_1 e_{ki;ay}e_{kj;by}e_{ik;d} &&\text{by \ref{R41}, noting that $xby=ay$.}
\intertext{In Case (ii), writing $b=xa$, we have}
E_{ij;ab}E_{ik;cd} 
&\sim_1 e_{ij;x}e_{ji;a}e_{ij;1}e_{ki;y}e_{ik;d} &&\text{by Lemma \ref{lem:Eijab}(ii) and (iii)}\\
&\sim_1 e_{ij;x}e_{ji;a}e_{ki;y}e_{kj;y}e_{ik;d} &&\text{by \ref{R41}}\\
&\sim_1 e_{ij;x}e_{ji;a}e_{kj;y}e_{ik;d} &&\text{by \ref{R3a1}}\\
&\sim_1 e_{ij;x}e_{ki;ay}e_{kj;y}e_{ik;d} &&\text{by \ref{R41}}\\
&\sim_1 e_{ki;ay}e_{kj;by}e_{kj;y}e_{ik;d} &&\text{by \ref{R41}, noting that $xay=by$}\\
&\sim_1 e_{ki;ay}e_{kj;by}e_{ik;d} &&\text{by \ref{R1a1}.}\\
\intertext{For Case (iii), write $a=xb$ and $d=zac$.  Then}
E_{ij;ab}E_{ik;cd} &\sim_1 e_{ji;x}e_{ij;b}e_{ik;za}e_{ki;c}e_{ik;1} &&\text{by Lemma \ref{lem:Eijab}(ii) and (iii)}\\
&= e_{ji;x}e_{ij;b}e_{ik;zxb}e_{ki;c}e_{ik;1} \\
&\sim_1 e_{ji;x}e_{jk;zx}e_{ij;b}e_{ki;c}e_{ik;1} &&\text{by \ref{R41}}\\
&\sim_1 e_{jk;zx}e_{ji;x}e_{ki;c}e_{kj;bc}e_{ik;1} &&\text{by \ref{R41} twice}\\
&\sim_1 e_{jk;zx}e_{ji;x}e_{kj;bc}e_{ik;1} &&\text{by \ref{R3a1}}\\
&\sim_1 e_{ji;x}e_{jk;zx}e_{ik;1}e_{ij;bc} &&\text{by \ref{R41} twice,}\\
\intertext{and}
E_{ik;ac,d}E_{ij;1,bc} &\sim e_{ik;z}e_{ki;ac}e_{ik;1}e_{ij;bc} &&\text{by Lemma \ref{lem:Eijab}(iii) and (iv)}\\
&\sim e_{ik;z}e_{ki;xbc}e_{kj;bc}e_{ik;1} &&\text{by \ref{R41} and $a=xb$}\\
&\sim e_{ik;z}e_{ji;x}e_{kj;bc}e_{ik;1} &&\text{by \ref{R41}}\\
&\sim e_{ji;x}e_{jk;zx}e_{ik;1}e_{ij;bc} &&\text{by \ref{R41} twice.}\\
\intertext{For Case (iv), write $b=xa$ and $d=zac$.  Then}
E_{ij;ab}E_{ik;cd} &\sim_1 e_{ij;x}e_{ji;a}e_{ij;1}e_{ik;za}e_{ki;c}e_{ik;1} &&\text{by Lemma \ref{lem:Eijab}(iii)}\\
&\sim_1 e_{ij;x}e_{ji;a}e_{jk;za}e_{ij;1}e_{ki;c}e_{ik;1} &&\text{by \ref{R41}}\\
&\sim_1 e_{ij;x}e_{jk;za}e_{ji;a}e_{ki;c}e_{kj;c}e_{ik;1} &&\text{by \ref{R41} twice}\\
&\sim_1 e_{ij;x}e_{jk;za}e_{ji;a}e_{kj;c}e_{ik;1} &&\text{by \ref{R3a1}}\\
&\sim_1 e_{ij;x}e_{jk;za}e_{ki;ac}e_{kj;c}e_{ik;1} &&\text{by \ref{R41}}\\
&\sim_1 e_{ij;x}e_{kj;1}e_{jk;za}e_{kj;c}e_{ki;ac}e_{ik;1} &&\text{by \ref{R3a1} and \ref{R41}}\\
&\sim_1 e_{ij;x}e_{jk;1}e_{kj;zac}e_{ki;ac}e_{ik;1} &&\text{by \ref{R1b1}}\\
&\sim_1 e_{ij;x}e_{kj;1}e_{ik;z}e_{ki;ac}e_{ik;1} &&\text{by Lemma \ref{lem:X1tech}(iii)}\\
&\sim_1 e_{ij;x}e_{ik;z}e_{ki;ac}e_{ik;1} &&\text{by \ref{R3a1},}\\
\intertext{and}
E_{ik;ac,d}E_{ij;1,bc} &\sim e_{ik;z}e_{ki;ac}e_{ik;1}e_{ij;bc} &&\text{by Lemma \ref{lem:Eijab}(iii) and (iv)}\\
&\sim e_{ik;z}e_{ki;ac}e_{kj;xac}e_{ik;1} &&\text{by \ref{R41} and $b=xa$}\\
&\sim e_{ik;z}e_{ij;x}e_{ki;ac}e_{ik;1} &&\text{by \ref{R41}}\\
&\sim e_{ij;x}e_{ik;z}e_{ki;ac}e_{ik;1} &&\text{by \ref{R41}.}\\
\intertext{For Case (v), write $a=xb$, $d=yc$ and $ac=zd$.  Then}
E_{ij;ab}E_{ik;cd} &\sim_1 e_{ji;x}e_{ij;b}e_{ik;y}e_{ki;c}e_{ik;1} &&\text{by Lemma \ref{lem:Eijab}(ii) and (iii)}\\
&\sim_1 e_{ji;x}e_{ik;y}e_{ij;b}e_{ki;c}e_{ik;1} &&\text{by \ref{R41}}\\
&\sim_1 e_{ji;x}e_{ik;y}e_{ki;c}e_{kj;bc}e_{ik;1} &&\text{by \ref{R41}}\\
&\sim_1 e_{ji;x}e_{ki;1}e_{ik;y}e_{ki;c}e_{ik;1}e_{ij;bc} &&\text{by \ref{R3a1} and \ref{R41}}\\
&\sim_1 e_{ji;x}e_{ik;1}e_{ki;d}e_{ik;1}e_{ij;bc} &&\text{by \ref{R1b1} and $d=yc$}\\
&\sim_1 e_{ji;x}e_{ki;1}e_{ik;d}e_{ij;bc} &&\text{by \ref{R1b1}}\\
&\sim_1 e_{ji;x}e_{ik;d}e_{ij;bc} &&\text{by \ref{R3a1}}\\
&\sim_1 e_{ki;z}e_{ik;d}e_{ij;bc} &&\text{by Lemma \ref{lem:X1tech}(iv), as $zd=ac=xbc$}\\
&\sim_1 E_{ik;ac,d}E_{ij;1,bc} &&\text{by Lemma \ref{lem:Eijab}(ii) and (iv).}\\
\intertext{For Case (vi), write $b=xa$, $d=yc$ and $ac=zd$.  Then}
E_{ij;ab}E_{ik;cd} &\sim_1 e_{ij;x}e_{ji;a}e_{ij;1}e_{ik;y}e_{ki;c}e_{ik;1} &&\text{by Lemma \ref{lem:Eijab}(iii)}\\
&\sim_1 e_{ij;x}e_{ji;a}e_{jk;y}e_{ij;1}e_{ki;c}e_{ik;1} &&\text{by \ref{R41}}\\
&\sim_1 e_{ij;x}e_{jk;y}e_{ji;a}e_{ki;c}e_{kj;c}e_{ik;1} &&\text{by \ref{R41} twice}\\
&\sim_1 e_{ij;x}e_{jk;y}e_{ji;a}e_{kj;c}e_{ik;1} &&\text{by \ref{R3a1}}\\
&\sim_1 e_{ij;x}e_{kj;1}e_{jk;y}e_{kj;c}e_{ki;ac}e_{ik;1} &&\text{by \ref{R3a1} and \ref{R41}}\\
&\sim_1 e_{ij;x}e_{jk;1}e_{kj;d}e_{ki;zd}e_{ik;1} &&\text{by \ref{R1b1}, $yc=d$ and $ac=zd$}\\
&\sim_1 e_{ij;x}e_{jk;1}e_{ji;z}e_{kj;d}e_{ik;1} &&\text{by \ref{R41}}\\
&\sim_1 e_{ij;x}e_{ki;z}e_{jk;1}e_{kj;d}e_{ik;1} &&\text{by \ref{R41}}\\
&\sim_1 e_{ki;z}e_{kj;xz}e_{jk;1}e_{ik;1}e_{ij;d} &&\text{by \ref{R41} twice}\\
&\sim_1 e_{ki;z}e_{kj;xz}e_{jk;1}e_{ij;d} &&\text{by \ref{R3a1}}\\
&\sim_1 e_{ki;z}e_{kj;xz}e_{ij;d}e_{ik;d} &&\text{by \ref{R41}}\\
&\sim_1 e_{ki;z}e_{kj;xz}e_{ik;d} &&\text{by \ref{R3a1}}\\
&\sim_1 e_{ki;z}e_{ik;d}e_{ij;bc} &&\text{by \ref{R41} and $xzd=xac=bc$}\\
&\sim_1 E_{ik;ac,d}E_{ij;1,bc} &&\text{by Lemma \ref{lem:Eijab}(ii) and (iv).}
\end{align*}
{\bf\boldmath \ref{R4b2'}:}  As in Remark \ref{rem:R4}, we only need to show that $E_{ik;c,bd}E_{ij;1,ad}\sim_1E_{jk;ab}E_{ij;cd}$.  Again, we must consider six cases:
\begin{itemize}
\begin{multicols}{3}
\item[(i)] $a\leq_{\L} b$ and $d\leq_{\L} c$;
\item[(ii)] $b\leq_{\L} a$ and $d\leq_{\L} c$;
\item[(iii)] $a\leq_{\L} b$ and $c\leq_{\L} bd$;
\item[(iv)] $b\leq_{\L} a$ and $c\leq_{\L} bd$;
\item[(v)] $a\leq_{\L} b$ and $bd\leq_{\L} c\leq_{\L} d$;
\item[(vi)] $b\leq_{\L} a$ and $bd\leq_{\L} c\leq_{\L} d$.
\end{multicols}
\end{itemize}
(Note that $bd\leq_{\L} c$ in Cases (i) and (ii), and that $c\leq_{\L} d$ in Cases (iii) and (iv).)
For Cases (i) and (ii), write $d=yc$, and note that $bd=(by)c$, so that
\begin{align*}
E_{ik;c,bd}E_{ij;1,ad} &\sim_1 e_{ik;by}e_{ki;c}e_{ik;1}e_{ij;ad} &&\text{by Lemma \ref{lem:Eijab}(iii) and (iv)}\\
&\sim_1 e_{ik;by}e_{ik;1}e_{ki;c}e_{ik;1}e_{ij;ad} &&\text{by \ref{R1a1}}\\
&\sim_1 e_{ik;by}e_{ki;1}e_{ik;c}e_{ij;ayc} &&\text{by \ref{R1b1} and $d=yc$}\\
&\sim_1 e_{ik;by}e_{ik;1}e_{ij;ay}e_{ji;c}e_{ij;1} &&\text{by Lemma \ref{lem:X1tech}(ii)}\\
&\sim_1 e_{ik;by}e_{ij;ay}e_{ji;c}e_{ij;1} &&\text{by \ref{R1a1}.}\\
\intertext{In Case (i), writing $a=xb$, we have}
E_{jk;ab}E_{ij;cd} &\sim_1 e_{kj;x}e_{jk;b} e_{ij;y}e_{ji;c}e_{ij;1} &&\text{by Lemma \ref{lem:Eijab}(ii) and (iii)}\\
&\sim_1 e_{kj;x}e_{ij;y}e_{ik;by} e_{ji;c}e_{ij;1} &&\text{by \ref{R41}}\\
&\sim_1 e_{kj;x}e_{ik;by} e_{ji;c}e_{ij;1} &&\text{by \ref{R3a1}}\\
&\sim_1 e_{ik;by}e_{ij;ay} e_{ji;c}e_{ij;1} &&\text{by \ref{R41} and $xb=a$.}\\
\intertext{In Case (ii), writing $b=xa$, we have}
E_{jk;ab}E_{ij;cd} &\sim_1 e_{jk;x}e_{kj;a}e_{jk;1} e_{ij;y}e_{ji;c}e_{ij;1} &&\text{by Lemma \ref{lem:Eijab}(iii)}\\
&\sim_1 e_{jk;x}e_{kj;a}e_{ij;y}e_{ik;y} e_{ji;c}e_{ij;1} &&\text{by \ref{R41}}\\
&\sim_1 e_{jk;x}e_{kj;a}e_{ik;y} e_{ji;c}e_{ij;1} &&\text{by \ref{R3a1}}\\
&\sim_1 e_{jk;x}e_{ik;y}e_{ij;ay} e_{ji;c}e_{ij;1} &&\text{by \ref{R41}}\\
&\sim_1 e_{jk;x}e_{ij;ay} e_{ji;c}e_{ij;1} &&\text{by \ref{R3a1}}\\
&\sim_1 e_{ik;by}e_{ij;ay} e_{ji;c}e_{ij;1} &&\text{by \ref{R41} and $xa=b$.}\\
\intertext{For Case (iii), write $a=xb$ and $c=zbd$.  Then}
E_{jk;ab}E_{ij;cd} &\sim_1 e_{kj;x}e_{jk;b} e_{ji;zb}e_{ij;d} &&\text{by Lemma \ref{lem:Eijab}(ii)}\\
&\sim_1 e_{kj;x}e_{ji;zb}e_{jk;b} e_{ij;d} &&\text{by \ref{R41}}\\
&\sim_1 e_{kj;x}e_{ij;1}e_{ji;zb}e_{ij;d}e_{ik;bd}  &&\text{by \ref{R3a1} and \ref{R41}}\\
&\sim_1 e_{kj;x}e_{ji;1}e_{ij;zbd}e_{ik;bd}  &&\text{by \ref{R1b1}}\\
&\sim_1 e_{kj;x}e_{ji;1}e_{kj;z}e_{ik;bd}  &&\text{by \ref{R41}}\\
&\sim_1 e_{kj;x}e_{kj;z}e_{ki;z}e_{ik;bd}  &&\text{by \ref{R41}}\\
&\sim_1 e_{kj;x}e_{ki;z}e_{ik;bd}  &&\text{by \ref{R1a1}}\\
&\sim_1 e_{ki;z}e_{kj;x}e_{ik;bd} &&\text{by \ref{R41}}\\
&\sim_1 e_{ki;z}e_{ik;bd}e_{ij;ad} &&\text{by \ref{R41} and $xb=a$}\\
&\sim_1 E_{ik;c,bd}E_{ij;1,ad} &&\text{by Lemma \ref{lem:Eijab}(ii) and (iv).}\\
\intertext{For Case (iv), write $b=xa$ and $c=zbd$.  Then}
E_{jk;ab}E_{ij;cd} &\sim_1 e_{jk;x}e_{kj;a}e_{jk;1} e_{ji;zb}e_{ij;d} &&\text{by Lemma \ref{lem:Eijab}(ii) and (iii)}\\
&\sim_1 e_{jk;x}e_{kj;a}e_{ji;zb}e_{jk;1} e_{ij;d} &&\text{by \ref{R41}}\\
&\sim_1 e_{jk;x}e_{kj;a}e_{ij;1}e_{ji;zb}e_{ij;d}e_{ik;d}  &&\text{by \ref{R3a1} and \ref{R41}}\\
&\sim_1 e_{jk;x}e_{kj;a}e_{ji;1}e_{ij;zbd}e_{ik;d}  &&\text{by \ref{R1b1}}\\
&\sim_1 e_{jk;x}e_{kj;a}e_{ji;1}e_{kj;zb}e_{ik;d}  &&\text{by \ref{R41}}\\
&\sim_1 e_{jk;x}e_{kj;a}e_{kj;zb}e_{ki;zb}e_{ik;d}  &&\text{by \ref{R41}}\\
&\sim_1 e_{jk;x}e_{kj;a}e_{ki;zxa}e_{ik;d}  &&\text{by \ref{R1a1} and $b=xa$}\\
&\sim_1 e_{jk;x}e_{ji;zx}e_{kj;a}e_{ik;d}  &&\text{by \ref{R41}}\\
&\sim_1 e_{ji;zx}e_{jk;x}e_{ik;d}e_{ij;ad}  &&\text{by \ref{R41} twice}\\
&\sim_1 e_{ji;zx}e_{jk;x}e_{ij;ad}  &&\text{by \ref{R3a1}}\\
&\sim_1 e_{ki;z}e_{jk;x}e_{ij;ad}  &&\text{by \ref{R41}}\\
&\sim_1 e_{ki;z}e_{ik;bd}e_{ij;ad}  &&\text{by \ref{R41} and $xa=b$}\\
&\sim_1 E_{ik;c,bd}E_{ij;1,ad}  &&\text{by  Lemma \ref{lem:Eijab}(ii) and (iv).}\\
\intertext{For Case (v), write $a=xb$, $c=yd$ and $bd=zc$.  Then}
E_{jk;ab}E_{ij;cd} &\sim_1 e_{kj;x}e_{jk;b} e_{ji;y}e_{ij;d} &&\text{by Lemma \ref{lem:Eijab}(ii)}\\
&\sim_1 e_{kj;x}e_{ji;y}e_{jk;b} e_{ij;d} &&\text{by \ref{R41}}\\
&\sim_1 e_{kj;x}e_{ij;1}e_{ji;y}e_{ij;d}e_{ik;bd} &&\text{by \ref{R3a1} and \ref{R41}}\\
&\sim_1 e_{kj;x}e_{ji;1}e_{ij;c}e_{ik;zc} &&\text{by \ref{R1b1}, $yd=c$ and $bd=zc$,}\\
\intertext{and}
E_{ik;c,bd}E_{ij;1,ad} &\sim_1 e_{ik;z}e_{ki;c}e_{ik;1}e_{ij;ad} &&\text{by Lemma \ref{lem:Eijab}(iii) and (iv)}\\
&\sim_1 e_{ik;z}e_{ki;c}e_{kj;xzc}e_{ik;1} &&\text{by \ref{R41} and $ad=xbd=xzc$}\\
&\sim_1 e_{ik;z}e_{ij;xz}e_{ki;c}e_{ik;1} &&\text{by \ref{R41}}\\
&\sim_1 e_{kj;x}e_{ik;z}e_{ki;c}e_{ik;1} &&\text{by \ref{R41}}\\
&\sim_1 e_{kj;x}e_{ij;1}e_{ik;z}e_{ki;c}e_{ik;1} &&\text{by \ref{R3a1}}\\
&\sim_1 e_{kj;x}e_{ji;1}e_{ij;c}e_{ik;zc} &&\text{by Lemma \ref{lem:X1tech}(ii).}\\
\intertext{For Case (vi), write $b=xa$, $c=yd$ and $bd=zc$.  Then}
E_{jk;ab}E_{ij;cd} &\sim_1 e_{jk;x}e_{kj;a}e_{jk;1} e_{ji;y}e_{ij;d} &&\text{by Lemma \ref{lem:Eijab}(ii) and (iii)}\\
&\sim_1 e_{jk;x}e_{kj;a}e_{ki;y}e_{jk;1}e_{ij;d} &&\text{by \ref{R41}}\\
&\sim_1 e_{jk;x}e_{ki;y}e_{kj;a}e_{ij;d}e_{ik;d} &&\text{by \ref{R41} twice}\\
&\sim_1 e_{jk;x}e_{ki;y}e_{kj;a}e_{ik;d} &&\text{by \ref{R3a1}}\\
&\sim_1 e_{jk;x}e_{kj;a}e_{ki;y}e_{ik;d} &&\text{by \ref{R41}}\\
&\sim_1 e_{jk;x}e_{kj;ad}e_{ki;c}e_{ik;1} &&\text{by Lemma \ref{lem:X1tech}(i) and $yd=c$}\\
&\sim_1 e_{jk;x}e_{ki;c}e_{kj;ad}e_{ik;1} &&\text{by \ref{R41}}\\
&\sim_1 e_{ik;z}e_{ki;c}e_{kj;ad}e_{ik;1} &&\text{by Lemma \ref{lem:X1tech}(iv), as $zc=bd=xad$}\\
&\sim_1 e_{ik;z}e_{ki;c}e_{ik;1}e_{ij;ad} &&\text{by \ref{R41}}\\
&\sim_1 E_{ik;c,bd}E_{ij;1,ad} &&\text{by Lemma \ref{lem:Eijab}(iii) and (iv).}
\end{align*}
This completes the proof of the theorem. \epf


As noted at the beginning of this section, some of the relations from $R_1$ may be simplified in the case that $M$ is a group.

\begin{thm}\label{thm:X1group}
If $M$ is a group, then $M\wr\Sing_n$ has presentation $\pres{X_1}{R_1'}$ via $\phi_1$, where $R_1'$ is obtained from $R_1$ by replacing 
\emph{\ref{R1a1}--\ref{R1e1}}
by
\begin{align}
\label{R1a1'}\tag*{(R1a)$'_1$} e_{ij;a}e_{ij;b} = e_{ij;a} &= e_{ji;a^{-1}}e_{ij;a} &&\text{for $a,b\in M$ and distinct $i,j$}\\
\label{R1b1}\tag*{(R1b)$_1$} e_{ij;1}e_{ji;a}e_{ij;b} &= e_{ji;1}e_{ij;ab} &&\text{for $a,b\in M$ and distinct $i,j$.}
\end{align}
\end{thm}

\pf Suppose $M$ is a group.  We start with the presentation $\pres{X_1}{R_1}$ from Therorem \ref{thm:X1}.  Since $(e_{ji;a^{-1}}e_{ij;a})\phi_1=e_{ij;a}\phi_1$, we may augment \ref{R1a1} by replacing it with \ref{R1a1'}.  Relation \ref{R1e1} may be removed, as it is just part of \ref{R1a1'}.  Relation \ref{R1c1} may also be removed, since $ac=bc$ is only possible in $M$ (a group) if $a=b$, in which case the relation is vacuous.  For relation \ref{R1d1}, let $i,j\in\bn$ be distinct, and let $a,b,c\in M$ with $abc=c$, noting that this forces $b=a^{-1}$.  Then, writing ${\sim_1'}=(R_1')^\sharp$, we have
\begin{align*}
e_{ji;a}e_{ij;bc} &\sim_1' e_{ij;b}e_{ji;a}e_{ij;bc} &&\text{by \ref{R1a1'} and $b=a^{-1}$}\\
&\sim_1' e_{ij;b}e_{ij;1}e_{ji;a}e_{ij;bc} &&\text{by \ref{R1a1'}}\\
&\sim_1' e_{ij;b}e_{ji;1}e_{ij;c} &&\text{by \ref{R1b1} and $abc=c$}\\
&\sim_1' e_{ij;b}e_{ij;1}e_{ji;c}e_{ij;1} &&\text{by \ref{R1b1}}\\
 &\sim_1' e_{ij;b}e_{ji;c}e_{ij;1} &&\text{by \ref{R1a1'}.} \qedhere
\end{align*}
\epf

\begin{rem}
Theorem \ref{thm:X1group} can also be deduced directly from Theorem \ref{thm:X2} in a similar way to Theorem~\ref{thm:X1}, though the calculations are far easier under the assumption that $M$ is a group.  Here we may simply define the words $E_{ij;ab}$ as $e_{ji;ab^{-1}}e_{ij;b}$ for all $i,j,a,b$, and we never need to consider multiple cases (according to whether $a\leq_{\L} b$ or $b\leq_{\L} a$, etc.).
\end{rem}

\begin{rem}
Although $M\wr\Sing_n\not=\la\X_1\ra$ in the case that $M/{\L}$ is not a chain, it is still the case (by Theorem \ref{thm:IG}) that $\la E(M\wr\Sing_n)\ra=\la\X_1\ra$.  It would be interesting to give a presentation for $\la E(M\wr\Sing_n)\ra$ in terms of the generating set $\X_1$.
\end{rem}

%

\section{An application to the endomorphism monoid of a uniform partition}\label{sect:EndP}

We now apply the results of previous sections to obtain a presentation for the idempotent-generated subsemigroup of the endomorphism monoid of a uniform partition of a finite set.  These monoids are defined as follows.  Let $X$ be a non-empty finite set, and let $\P=\{C_1,\ldots,C_n\}$ be a uniform partition of $X$: by this we mean that the sets $C_1,\ldots,C_n$ are pairwise disjoint, have a common size ($m$, say), and their union is all of~$X$ (so that $|X|=mn$).  The \emph{endomorphism monoid of $\P$} is defined to be the submonoid
\[
\TXP = \set{\al\in\T_X}{(\forall i\in\bn)(\exists j\in\bn)\ C_i\al\sub C_j}
\]
of $\T_X$, the full transformation semigroup on $X$.  These monoids were introduced by Pei in \cite{Pei1994}, and have subsequently been studied by a number of different authors.  In particular, the rank of $\TXP$ was calculated in \cite{AS2009}, while the rank and idempotent rank of the idempotent-generated subsemigroup of $\TXP$ were calculated in \cite{DE2016}; for the corresponding studies of the non-uniform case, see \cite{DEM2016,ABMS2015}.  
%
%
As noted in \cite{AS2009}, $\TXP$ is isomorphic to the wreath product $\T_m\wr\T_n$.  
So we will concentrate on such a wreath product $\T_m\wr\T_n$, and our goal (as stated above) is to give a monoid presentation for the idempotent-generated subsemigroup $\la E(\T_m\wr\T_n)\ra$.  In fact, we are able to solve a more general problem: namely, in Theorem \ref{thm:IGMwrTn}, we give a monoid presentation for the idempotent-generated subsemigroup $\la E(M\wr\T_n)\ra$ of $M\wr\T_n$, modulo a presentation for $\la E(M)\ra$, in the case that $M$ is a monoid satisfying $\la E(M)\ra=\{1\}\cup(M\sm G)$, where $G$ is the group of units of $M$.
%
%
So for the remainder of this section, we fix a monoid $M$, write $G$ for its group of units, and we assume that $\la E(M)\ra = \{1\}\cup(M\sm G)$.  Examples of such monoids $M$ include the finite full transformation semigroups \cite{Howie1966}, finite dimensional full linear monoids \cite{Erdos1967}, finite partition monoids \cite{JEpnsn}, finite Brauer monoids \cite{MM2007}, the endomorphism monoids of certain finite dimensional free $M$-acts (see Corollary~\ref{cor:IGMwrTn} above), and many more.  Presentations for the idempotent-generated subsemigroups of some (but not all) of these examples are known \cite{JEpnsn,JEptnsn2,JEtnsn2,MM2007}.

To avoid confusion, we will write $1$ and $1_n$ for the identity elements of $M$ and $\T_n$, respectively.  Recall that any subsemigroup $S$ of $\T_n$ leads to a wreath product $K\wr S$, for any monoid $K$.  In particular, when $S=\{1_n\}\sub\T_n$ consists of only the identity transformation, $K\wr\{1_n\}$ is isomorphic to the direct product of $n$ copies of $K$.  The $M=\T_m$ case of the next result is contained in \cite[Proposition 4.1]{DE2016}.  We write $A=B\sqcup C$ to indicate that $A$ is the disjoint union of $B$ and $C$.
To simplify notation throughout this section, if $T$ is any semigroup, we will write $\bbE(T)=\la E(T)\ra$ for the idempotent-generated subsemigroup of~$T$.

\begin{prop}\label{prop:decomp}
Suppose $M$ is a monoid with group of units $G$, and that $\bbE(M) = \{1\}\cup(M\sm G)$.  Then
\bit
\itemit{i} $\bbE(M\wr\T_n) = (\bbE(M)\wr\{1_n\})\sqcup(M\wr\Sing_n)$\emph{;} and
\itemit{ii} $M\wr\Sing_n=(\bbE(M)\wr\{1_n\})(G\wr\Sing_n)$. 
\end{itemize}
\end{prop}

\pf (ii).  If $(\ba,1_n)\in\bbE(M)\wr\{1_n\}$ and $(\bb,\be)\in G\wr\Sing_n$, then ${(\ba,1_n)(\bb,\be)=(\ba\bb,\be)\in M\wr\Sing_n}$.  

Conversely, suppose $(\bc,\ga)\in M\wr\Sing_n$.  For each $i\in\bn$, define
\[
a_i =\begin{cases}
1 &\text{if $c_i\in G$}\\
c_i &\text{if $c_i\in M\sm G$}
\end{cases}
\AND
b_i =\begin{cases}
c_i &\text{if $c_i\in G$}\\
1 &\text{if $c_i\in M\sm G$.}
\end{cases}
\]
Since $a_ib_i=c_i$ for each $i$, it follows that $(\bc,\ga)=(\ba,1_n)(\bb,\ga)$.  It is also clear that $(\bb,\ga)\in G\wr\Sing_n$, while $(\ba,1_n)\in\bbE(M)\wr\{1_n\}$ follows from the fact that $\bbE(M) = \{1\}\cup(M\sm G)$.  This completes the proof of (ii).

\pfitem{i}  First suppose $(\ba_1,\al_1),\ldots,(\ba_k,\al_k)\in E(M\wr\T_n)$, and write $(\ba,\al)=(\ba_1,\al_1)\cdots(\ba_k,\al_k)$.  If any of $\al_1,\ldots,\al_k$ belongs to $\Sing_n$, then so too does $\al=\al_1\cdots\al_k$, so that $(\ba,\al)\in M\wr\Sing_n$.  On the other hand, if $\al_1=\cdots=\al_k=1_n$, then Lemma \ref{lem:idempotents} gives $\ba_i\in E(M)^n$ for all $i$, in which case
\[
(\ba,\al)=(\ba_1\cdots\ba_k,1_n)\in\bbE(M)\wr\{1_n\}.
\]
We have shown that $\bbE(M\wr\T_n) \sub (\bbE(M)\wr\{1_n\})\sqcup(M\wr\Sing_n)$.  

To prove the reverse containment, first suppose $(\bb,1_n)\in\bbE(M)\wr\{1_n\}$.  Since $1\in E(M)$, we may write $\bb=\bb_1\cdots\bb_k$ with $\bb_1,\ldots,\bb_k\in E(M)^n$, and it then follows that $(\bb,1_n)=(\bb_1,1_n)\cdots(\bb_k,1_n)\in\bbE(M\wr\T_n)$.  This shows that $\bbE(M)\wr\{1_n\} \sub \bbE(M\wr\T_n)$.  Together with Part (ii) and the fact that $G\wr\Sing_n$ is idempotent-generated (by Theorem \ref{thm:IG}(iii)), it also follows that $M\wr\Sing_n\sub\bbE(M\wr\T_n)$.  \epf

Suppose now that $\bbE(M)=\la E(M)\ra$ has monoid presentation $\pres YQ$ via $\psi:Y^*\to\bbE(M)$.  Since this is a \emph{monoid} presentation, we may assume that $y\psi\not=1$ for all $y\in Y$, and it will be important to do so in what follows.  Define new alphabets $Y_{(i)}=\set{y_{(i)}}{y\in Y}$ for each $i\in\bn$, and put $\bY=Y_{(1)}\cup\cdots\cup Y_{(n)}$.  For a word $w=y_1\cdots y_k\in Y^*$, and for $i\in\bn$, define $w_{(i)}=(y_1)_{(i)}\cdots(y_k)_{(i)}\in Y_{(i)}^*$.  For each $i\in\bn$, write $Q_{(i)}=\set{(u_{(i)},v_{(i)})}{(u,v)\in Q}$ and put $\bQ=Q_{(1)}\cup\cdots\cup Q_{(n)}$.  We also define
\[
\RC = \set{(x_{(i)}y_{(j)},y_{(j)}x_{(i)})}{x,y\in Y,\ i,j\in\bn,\ i\not=j}.
\]
For $a\in M$ and $i\in\bn$, write $a_{(i)}=((1,\ldots,1,a,1,\ldots,1),1_n)$, where the $a$ is in the $i$th position.  Define an epimorphism
\[
\Psi:\bY^*\to \bbE(M)\wr\{1_n\}: y_{(i)}\mt (y\psi)_{(i)}.
\]
The next result follows from an obvious (and essentially folklore) result on presentations for direct products of monoids.

\begin{lemma}\label{lem:Em}
With the above notation, the monoid $\bbE(M)\wr\{1_n\}$ has monoid presentation $\pres{\bY}{\bQ\cup \RC}$ via~$\Psi$.~\epfres
\end{lemma}

We fix the semigroup presentation $\pres{X_1}{R_1'}$ for $G\wr\Sing_n$ 
(via $\phi_1:X_1^+\to G\wr\Sing_n$) 
from Theorem \ref{thm:X1group}, where $X_1=\set{e_{ij;a}}{i,j\in\bn,\ i\not=j,\ a\in G}$, and so on.  We now explain how to stitch this together with the monoid presentation $\pres{\bY}{\bQ\cup \RC}$ for $\bbE(M)\wr\{1_n\}$ in order to yield a monoid presentation for $\bbE(M\wr\T_n)$.  In what follows, it will be convenient to write $\wb=w\psi\in\bbE(M)$ for any word $w\in Y^*$.  By Proposition \ref{prop:decomp}, we may define an epimorphism
\[
\Theta:(\bY\cup X_1)^*\to\bbE(M\wr\T_n):y_{(i)}\mt\yb_{(i)},\ e_{ij;a}\mt\ve_{ij;a}.
\]
We will also choose (and fix for the remainder of the section) a set of words $\set{h_a}{a\in\bbE(M)}\sub Y^*$ such that $h_a\psi=a$ for all $a\in\bbE(M)$.  For $a\in\bbE(M)$ and $i\in\bn$, define $h_{a;i}=(h_a)_{(i)}\in Y_{(i)}^*$, noting that $h_{a;i}\Theta=h_{a;i}\Psi=a_{(i)}$.  In practice, we might like to choose the words $h_a$ to be as short as possible, but this is not a requirement.  Note that if $w\in \bY^*$ is such that $w\Theta = ((a_1,\ldots,a_n),1_n)\in\bbE(M)\wr\{1_n\}$, then $w$ may be transformed into $h_{a_1;1}\cdots h_{a_n;n}$ using relations~$\bQ\cup\RC$, by Lemma \ref{lem:Em}.

Now let $R_\nabla$ denote the set of relations
\allowdisplaybreaks[0]
\begin{alignat}{2}
\tag{$\nabla$1a}\label{N1a}     &\mathrel{\phantom{=}}\hphantom{\hack} y_{(i)}h_{a\yb;j}e_{ij;1} &\quad&\qquad\text{if $k=i$} \\
\tag{$\nabla$1b}\label{N1b} e_{ij;a}y_{(k)} &=\hack e_{ij;a} &&\qquad\text{if $k=j$} \\
\tag{$\nabla$1c}\label{N1c}     &\mathrel{\phantom{=}}\hphantom{\hack} y_{(k)}e_{ij;a} 	&&\qquad\text{otherwise} \\[1ex]
\tag*{($\nabla$2)}\label{N2}  y_{(j)}e_{ij;a} &= h_{\yb a;j}e_{ij;1} \\
\tag*{($\nabla$3)}\label{N3}  y_{(i)}e_{ji;a}e_{ij;b} &= h_{\yb ab;i}e_{ij;b},
\end{alignat}			
\allowdisplaybreaks
where $y\in Y$ in each relation, and $i,j,k,a,b$ range over all allowable values, subject to the stated constraints.  Note that the assumption that $\yb=y\psi\not=1$ for all $y\in Y$ (and the assumption that $\bbE(M)=\{1\}\cup(M\sm G)$) implies that $a\yb,\yb a\in M\sm G\sub\bbE(M)$ for all $y\in Y$ and $a\in G$ (so that the words $h_{a\yb},h_{\yb a},h_{\yb ab}$ appearing in the above relations are well defined).  The main result we wish to prove in this section is the following.

\begin{thm}\label{thm:IGMwrTn}
Suppose $M$ is a monoid with group of units $G$, and that $\bbE(M) = \la E(M)\ra = \{1\}\cup(M\sm G)$.  
With the above notation, the idempotent-generated subsemigroup $\bbE(M\wr\T_n)=\la E(M\wr\T_n)\ra$ of $M\wr\T_n$ has monoid presentation $\pres{\bY\cup X_1}{\bQ\cup\RC\cup R_1'\cup R_\nabla}$ via $\Theta$.
\end{thm}

To prove Theorem \ref{thm:IGMwrTn}, we first need some preliminary lemmas.  We will write ${\ssim}=(\bQ\cup\RC\cup R_1'\cup R_\nabla)^\sharp$ for the congruence on $(\bY\cup X_1)^*$ generated by the relations $\bQ\cup\RC\cup R_1'\cup R_\nabla$.  The next result follows by a simple diagrammatic check that the relations $R_\nabla$ are preserved by $\Theta$.

\begin{lemma}\label{lem:rels}
We have ${\ssim}\sub\ker(\Theta)$.  \epfres
\end{lemma}

As usual, proving the reverse containment is more of a challenge.

\begin{lemma}\label{lem:w1w2}
If $w\in(\bY\cup X_1)^*$, then $w\ssim w_1w_2$ for some $w_1\in \bY^*$ and $w_2\in X_1^*$.  If $w\not\in \bY^*$, then $w_2\in X_1^+$.
\end{lemma}

\pf For a word $u\in(\bY\cup X_1)^*$, we write $\xi(u)$ for the number of letters from $X_1$ appearing in $u$.  We prove the lemma by induction on $\xi(w)$.  If $\xi(w)=0$, then we are already done (with $w_1=w$ and $w_2=1$), so suppose $\xi(w)\geq1$, and write $w=ue_{ij;a}v$, where $u\in(\bY\cup X_1)^*$ and $v\in \bY^*$, so $\xi(u)=\xi(w)-1$.  By (\ref{N1a}--\ref{N1c}), we have $e_{ij;a}v\ssim ze_{ij;b}$ for some $z\in \bY^*$ and some $b\in G$.  Since $\xi(uz)=\xi(u)=\xi(w)-1$, an induction hypothesis gives $uz\ssim u_1u_2$ for some $u_1\in \bY^*$ and $u_2\in X_1^*$.  So $w=ue_{ij;a}v\ssim uze_{ij;b}\ssim u_1u_2e_{ij;b}$, and we are done (with $w_1=u_1\in\bY^*$ and $w_2=u_2e_{ij;b}\in X_1^+$).  (Note that the final assertion in the lemma follows from the above argument.)
\epf

We now improve Lemma \ref{lem:w1w2} by showing that the two words $w_1,w_2$ can be chosen to have a very specific form, in the case that $w\not\in \bY^*$.

\begin{lemma}\label{lem:normal}
Let $w\in(\bY\cup X_1)^*\sm \bY^*$, and write $w\Theta=(\ba,\al)$.  For $i\in\bn$, define
\[
b_i = 
\begin{cases}
1 &\text{if $a_i\in G$}\\
a_i &\text{if $a_i\in M\sm G$}
\end{cases}
\AND
c_i = 
\begin{cases}
a_i &\text{if $a_i\in G$}\\
1 &\text{if $a_i\in M\sm G$.}
\end{cases}
\]
Then $w\ssim w_1w_2$ for some $w_1\in \bY^*$ and $w_2\in X_1^+$ with $w_1\Theta=(\bb,1_n)$ and $w_2\Theta=(\bc,\al)$.
\end{lemma}

\pf By Lemma \ref{lem:w1w2}, the set $\Omega = \set{(w_1,w_2)\in \bY^*\times X_1^+}{w\ssim w_1w_2}$ is non-empty.  We define a function $\xi:\Om\to\N$ as follows.  Let $(w_1,w_2)\in\Om$, and write
\[
w_1\Theta = (\bp,1_n) \ANd w_2\Theta=(\bq,\al) \qquad\text{where $\bp\in\bbE(M)^n$ and $\bq\in G^n$.}
\]
(Because the last coordinate of $w_1\Theta$ must be $1_n$, it follows that the last coordinate of $w_2\Theta$ must be $\al$.)  We then define $\xi(w_1,w_2)$ to be the cardinality of the set $\Xi(w_1,w_2)=\set{i\in\bn}{(p_i,q_i)\not=(b_i,c_i)}$.  (Here, we assume that $\bp=(p_1,\ldots,p_n)$ and $\bq=(q_1,\ldots,q_n)$.)

We now choose a pair $(w_1,w_2)\in\Om$ for which $\xi(w_1,w_2)$ is minimal.  The proof will of course be complete if we can show that $\xi(w_1,w_2)=0$.
To do so, suppose to the contrary that $\xi(w_1,w_2)\geq1$.  As above, write $w_1\Theta = (\bp,1_n)$ and $w_2\Theta=(\bq,\al)$, noting that 
\[
(\ba,\al)=w\Theta=(w_1\Theta)(w_2\Theta)
=(\bp\bq,\al).
\]
This gives $p_iq_i=a_i$ for all $i$.  Since $w_2\in X_1^+$, we have $\al\in\Sing_n$, so we may fix some $(i,j)\in\ker(\al)$ with $i\not=j$.  By relabelling the elements of $\bn$ if necessary, we may assume that $(i,j)=(1,2)$.  Define words
\[
u_1 = (e_{21;q_1q_2^{-1}}e_{12;q_2})\cdot(e_{23;q_3}e_{32;1})\cdots(e_{2n;q_n}e_{n2;1})
\ \text{and}\ 
u_2 = (e_{12;q_2q_1^{-1}}e_{21;q_1})\cdot(e_{13;q_3}e_{31;1})\cdots(e_{1n;q_n}e_{n1;1}),
\]
and let $v$ be any word over $X$ (regarded as a subset of $X_1$ as usual) with $v\Theta=((1,\ldots,1),\al)$.
It is easy to check (diagrammatically) that $u_1\Theta=(\bq,\ve_{12})$ and $u_2\Theta=(\bq,\ve_{21})$.  In particular, since $\al=\ve_{12}\al=\ve_{21}\al$, we have
\[
(u_1\Theta)(v\Theta) = (u_2\Theta)(v\Theta) = (\bq,\al)=w_2\Theta.
\]
Since $w_2,u_1v,u_2v$ all belong to $X_1^+$, Theorem \ref{thm:X1group} then gives $w_2\ssim u_1v\ssim u_2v$.  As noted earlier, Lemma \ref{lem:Em} also gives $w_1\ssim h_{p_1;1}\cdots h_{p_n;n}$.  

Since $\xi(w_1,w_2)\geq1$, we may fix some $r\in\Xi(w_1,w_2)$.  Note that we could not have $p_r=1$, or else then $a_r=p_rq_r=q_r\in G$, which would give $(b_r,c_r)=(1,a_r)=(p_r,q_r)$, contradicting our assumption that $r\in\Xi(w_1,w_2)$.  In particular, $h_{p_r;r}\not=1$, so we may write $h_{p_r;r}=(y_1)_{(r)}\cdots(y_k)_{(r)}y_{(r)}$, where ${y_1,\ldots,y_k,y\in Y}$ (and where $p_r$ is therefore equal to $\yb_1\cdots\yb_k\yb$).  
Note that $\RC$ gives $w_1\ssim w_3h_{p_r;r}$, where $w_3=h_{p_1;1}\cdots h_{p_{r-1};r-1}h_{p_{r+1};r+1} \cdots h_{p_n;n}$. 
Note also that
\[
p_r\not=1 \implies p_r\in M\sm G \implies a_r=p_rq_r\in M\sm G \implies (b_r,c_r)=(a_r,1).
\]
We now consider separate cases, depending on the value of $r$.

\pfcase{1}
Suppose first that $r\geq3$.  Note that
\begin{align*}
h_{p_r;r} u_1 &=  (y_1)_{(r)}\cdots(y_k)_{(r)}y_{(r)} (e_{21;q_1q_2^{-1}}e_{12;q_2})\cdot(e_{23;q_3}e_{32;1})\cdots(e_{2,r-1;q_{r-1}}e_{r-1,2;1})\\
&\hspace{3cm} \times (e_{2r;q_r}e_{r2;1})(e_{2,r+1;q_{r+1}}e_{r+1,2;1})\cdots(e_{2n;q_n}e_{n2;1})\\
&\ssim  (y_1)_{(r)}\cdots(y_k)_{(r)} (e_{21;q_1q_2^{-1}}e_{12;q_2})\cdot(e_{23;q_3}e_{32;1})\cdots(e_{2,r-1;q_{r-1}}e_{r-1,2;1})\\
&\hspace{3cm} \times y_{(r)}(e_{2r;q_r}e_{r2;1})(e_{2,r+1;q_{r+1}}e_{r+1,2;1})\cdots(e_{2n;q_n}e_{n2;1}) &&\text{by \eqref{N1c}}\\
&\ssim  (y_1)_{(r)}\cdots(y_k)_{(r)} (e_{21;q_1q_2^{-1}}e_{12;q_2})\cdot(e_{23;q_3}e_{32;1})\cdots(e_{2,r-1;q_{r-1}}e_{r-1,2;1})\\
&\hspace{3cm} \times h_{\yb q_r;r}(e_{2r;1}e_{r2;1})(e_{2,r+1;q_{r+1}}e_{r+1,2;1})\cdots(e_{2n;q_n}e_{n2;1}) &&\text{by \ref{N2}}\\
&\ssim  (y_1)_{(r)}\cdots(y_k)_{(r)}h_{\yb q_r;r} (e_{21;q_1q_2^{-1}}e_{12;q_2})\cdot(e_{23;q_3}e_{32;1})\cdots(e_{2,r-1;q_{r-1}}e_{r-1,2;1})\\
&\hspace{3cm} \times (e_{2r;1}e_{r2;1})(e_{2,r+1;q_{r+1}}e_{r+1,2;1})\cdots(e_{2n;q_n}e_{n2;1}) &&\text{by \eqref{N1c}.}
\end{align*}
(In the last step of the previous calculation, recall that $h_{\yb q_r;r}$ involves only letters from $Y_{(r)}$.)  Note also that
\[
((y_1)_{(r)}\cdots(y_k)_{(r)}h_{\yb q_r;r})\Theta 
= (\yb_1\cdots\yb_k\yb q_{r})_{(r)}
= (p_rq_{r})_{(r)}
= (a_r)_{(r)}.
\]
As seen above, $a_r\in M\sm G$, so Lemma \ref{lem:Em} gives $(y_1)_{(r)}\cdots(y_k)_{(r)}h_{\yb q_r;r}\ssim h_{a_r;r}$.  Now put
\[
u_3 = (e_{21;q_1q_2^{-1}}e_{12;q_2})\cdot(e_{23;q_3}e_{32;1})\cdots(e_{2,r-1;q_{r-1}}e_{r-1,2;1})(e_{2r;1}e_{r2;1})(e_{2,r+1;q_{r+1}}e_{r+1,2;1})\cdots(e_{2n;q_n}e_{n2;1}).
\]
The above calculations show that $h_{p_r;r} u_1 \ssim h_{a_r;r}u_3$, and it follows that
\[
w\ssim w_1w_2\ssim (w_3h_{p_r;r} )(u_1v) \ssim (w_3h_{a_r;r})(u_3v) = v_1v_2,
\]
where $v_1=w_3h_{a_r;r}\in\bY^*$ and $v_2=u_3v\in X_1^+$.  It follows that $(v_1,v_2)\in\Om$.  But one may easily check that
\[
v_1\Theta = ((p_1,\ldots,p_{r-1},a_r,p_{r+1},\ldots,p_n),1_n) \AND
v_2\Theta = ((q_1,\ldots,q_{r-1},1,q_{r+1},\ldots,q_n),\al).
\]
Since $(b_r,c_r)=(a_r,1)$, it follows that $\xi(v_1,v_2)=\xi(w_1,w_2)-1$, contradicting the minimality of $\xi(w_1,w_2)$, and completing the proof in this case.

\pfcase{2}
Next, suppose $r=1$.  So now we have $w_1\ssim w_3h_{p_1;1}$ and $h_{p_1;1}=(y_1)_{(1)}\cdots(y_k)_{(1)}y_{(1)}$.  First note that $y_{(1)}e_{21;q_1q_2^{-1}}e_{12;q_2} \ssim h_{\yb q_1;1}e_{12;q_2}$, by \ref{N3}.  As in the previous case, we have $(y_1)_{(1)}\cdots(y_k)_{(1)}h_{\yb q_1;1}\ssim h_{a_1;1}$.  It quickly follows that $h_{p_1;1}u_1\ssim h_{a_1;1}u_3$, where $u_3=e_{12;q_2}\cdot(e_{23;q_3}e_{32;1})\cdots(e_{2n;q_n}e_{n2;1})$.    So $w\ssim w_1w_2 \ssim (w_3h_{p_1;1})(u_1v) \ssim (w_3h_{a_1;1})(u_3v)=v_1v_2$, where $v_1=w_3h_{a_1;1}\in\bY^*$ and $v_2=u_3v\in X_1^+$.  This time we have
\[
v_1\Theta=((a_1,p_2,\ldots,p_n),1_n) \AND v_2\Theta=((1,q_2,\ldots,q_n),\al),
\]
and again we obtain $\xi(v_1,v_2)=\xi(w_1,w_2)-1$, a contradiction.

\pfcase{3}
The case in which $r=2$ is almost identical to the previous case, but we use the word $u_2$ (defined above) instead of $u_1$.
\epf


We are now ready to tie the loose ends together.

\pf[\bf Proof of Theorem \ref{thm:IGMwrTn}.]  It remains to prove that $\ker(\Theta)\sub{\ssim}$, so suppose $u,v\in(\bY\cup X_1)^*$ are such that $u\Theta=v\Theta$.  If $u\in \bY^*$, then $u\Theta\in\bbE(M)\wr\{1_n\}$, which also gives $v\in \bY^*$; in this case, $u\ssim v$ follows from Lemma~\ref{lem:Em}.  So suppose $u\not\in \bY^*$, noting that this also forces $v\not\in \bY^*$.
Lemma \ref{lem:normal} then gives $u\ssim u_1u_2$ and $v\ssim v_1v_2$ for some $u_1,v_1\in \bY^*$ and $u_2,v_2\in X_1^+$ with $u_1\Theta=v_1\Theta$ and $u_2\Theta=v_2\Theta$.
Lemma \ref{lem:Em} and Theorem~\ref{thm:X1group} (respectively) then give $u_1\ssim v_1$ and $u_2\ssim v_2$.  Putting this all together, we obtain $u\ssim u_1u_2\ssim v_1v_2\ssim v$. \epf

\begin{rem}
As noted above, in the case that $M=\T_m$, Theorem \ref{thm:IGMwrTn} gives a presentation for the idempotent-generated subsemigroup of $\T_m\wr\T_n\cong\TXP$, where $\TXP$ is the endomorphism monoid of a uniform partition $\P$ of the set $X$ into $n$ blocks of size $m$.  Here we have $G=\S_m$, and the monoid presentation $\pres YQ$ for $\bbE(\T_m)=\{1\}\cup(\TmSm)$ is deduced from the semigroup presentation for $\TmSm$ in Theorem \ref{thm:X}.
\end{rem}

\subsection*{Acknowledgements}

The first author is supported by Grant No.~fsyq201408 of the Outstanding Young Teacher Training Programme of Foshan University.
The second author would like to thank the Higher Committee of Education Development (HCED) in Iraq for the financial support to pursue her PhD at the University of York.
The third author is supported by Grant No.\ 174019 of the Ministry of Education, Science, and Technological Development of the Republic of Serbia.
The first and fourth authors thank the University of York for its hospitality during their stay while this research was being conducted.
The {\sc Semigroups} package for GAP \cite{GAP} was useful during the early stages of the project.

\footnotesize
\def\bibspacing{-1.1pt}
\bibliography{biblio}
\bibliographystyle{abbrv}
\end{document}